\documentclass[12pt]{article}

\usepackage[utf8]{inputenc}
\usepackage[a4paper,nomarginpar]{geometry}
\usepackage[english]{babel}
\usepackage{enumitem}
\usepackage{amsmath,amsfonts,amssymb,amsthm,amscd}
\usepackage{tikz,xcolor}
\usepackage{mathrsfs}
\usepackage[bookmarks=true,hyperindex,pdftex,colorlinks,citecolor=red, linkcolor=blue
]{hyperref}
\usepackage{graphicx}
\usepackage{float}
\usepackage{caption}

\allowdisplaybreaks

\theoremstyle{plain}
\newtheorem{theorem}{Theorem}[section]
\newtheorem{proposition}[theorem]{Proposition}

\newtheorem{lemma}[theorem]{Lemma}
\newtheorem{corollary}[theorem]{Corollary}

\theoremstyle{definition}
\newtheorem{definition}[theorem]{Definition}
\newtheorem{remark}[theorem]{Remark}

\newtheorem{example}[theorem]{Example}

\newcommand{\N}{\mathbb{N}}
\newcommand{\Z}{\mathbb{Z}}

\newcommand{\K}{\mathbb{K}}

\DeclareMathOperator{\Chi}{Chi}
\DeclareMathOperator{\Gen}{Gen}
\DeclareMathOperator{\Par}{par}

\DeclareMathOperator{\Span}{span}
\DeclareMathOperator{\card}{card}
\DeclareMathOperator{\supp}{supp}
\DeclareMathOperator{\NS}{NS}
\DeclareMathOperator{\MProx}{MProx}
\DeclareMathOperator{\MAsym}{MAsym}

\def\la{\lambda}

\newcommand{\abs}[1]{\left\lvert#1\right\rvert}

\newcommand{\norm}[1]{\lVert#1\rVert}

\newcommand{\ldens}{\operatorname{\underline{dens}}}
\newcommand{\udens}{\operatorname{\overline{dens}}}

\newcommand{\eps}{\varepsilon}

\def\root{\texttt{r}}

\title{On the dynamics of weighted shifts on Fréchet sequence spaces over trees}

\author{Mateus Azevedo, Fernando Costa Jr.,   Jo\~ao V. A. Pinto}

\date{}

\begin{document}
	
	\maketitle
	
\begin{abstract}
	In this paper, we give characterizations of Li--Yorke chaos, mean Li--Yorke chaos and distributional chaos for weighted backward shifts acting on general Fréchet sequence spaces over directed trees. As corollaries, we obtain explicit characterizations in the setting of Köthe spaces $\lambda_p(A,V)$, where $p\in\{0\}\cup[1,\infty)$ and $V$ is a directed tree. In particular, we obtain criteria for the case of the classical Banach spaces $\ell^p(V)$, $p\in[1,\infty)$, and $c_0(V)$.
\end{abstract}
	
	
\section{Introduction}
    
Linear Dynamics is the subfield of dynamical systems devoted to the study of dynamical properties of continuous linear operators on topological vector spaces, in particular Fréchet and Banach spaces. The study of this kind of system goes back to the works of Birkhoff \cite{Birkhoff}, MacLane \cite{MacLane}, and Rolewicz \cite{Rolewicz}, who found the first hypercyclic operators. Hypercyclicity later became the main notion in the theory of linear dynamics. The first systematic studies of this property were carried out in the works of Kitai \cite{Kitai} and Godefroy and Shapiro \cite{GodefroyShapiro}. In this context, other related properties that have also been extensively studied include weak mixing, mixing, and Devaney chaos. The basic theory can be found in the books \cite{BaMa,GE-P}.

Over the last two decades, other notions of chaos, such as Li--Yorke chaos and distributional chaos, have been extensively studied in Linear Dynamics. These concepts are based on the dynamical behavior of pairs of points. In \cite{BBMP11}, Li--Yorke chaos and distributional chaos were studied for the first time in the general setting of operators on Banach spaces. In \cite{Nilson0}, distributional chaos was studied in the setting of Fréchet spaces, where it was characterized in terms of the existence of a single vector, the so-called distributionally irregular vector, and the useful Distributional Chaos Criterion (DCC) was established. We will make use of these results in our work. An analogous approach to Li--Yorke chaos was developed in \cite{BBMP2}: the property was described through the existence of a single vector, called an irregular vector, and a criterion for Li--Yorke chaos was obtained. Both ingredients will be used later in our paper.

During the past decade, increasing attention has been devoted to average dynamical properties, such as mean equicontinuity and mean sensitivity (see \cite{Mean}). Within this framework, mean Li--Yorke chaos, introduced in \cite{Hu}, has emerged as an important notion. For seminal works on mean Li--Yorke chaos in the setting of operators on Banach spaces, we refer the reader to \cite{BBP,BBPW}. In \cite{JiangLi}, the authors extended the notions and results established in \cite{Nilson0,BBMP2,BBP} to the setting of endomorphisms on metrizable groups. We will use their results in our characterizations of mean Li--Yorke chaos.

Obtaining explicit characterizations of dynamical properties for particular classes of operators is fundamental to the development of Linear Dynamics. Such characterizations facilitate the construction of examples and counterexamples and help distinguish between different dynamical properties. Weighted shifts constitute one of the most important and widely studied classes of operators in the theory. Their simple structure allows their iterates to be computed explicitly, while the freedom to choose the weights and the underlying sequence space makes them remarkably flexible. Consequently, classical weighted shifts on sequence spaces over $\N$ or $\Z$ have long served as a testing ground for general questions in Linear Dynamics. This particular kind of operators also belong to the broader class of weighted composition operators, for which properties such as shadowing and expansivity have been deeply investigated. The characterizations of Li--Yorke, mean Li--Yorke, and distributional chaos obtained in \cite{BerBonPin26,Fernanda,joaojmaa} are of particular relevance to the present paper, since we use some of the strategies developed in these works, especially in Sections~\ref{sec:LY} and~\ref{sec:dist}.

The importance of weighted shifts in the classical theory naturally motivates their study on more general index sets. Replacing $\N$ and $\Z$ by a directed tree introduces the geometry of the tree as an additional element of the dynamics. Given a directed tree $V$ and a sequence space $X\subset\K^V$, a weighted backward shift on $X$ is an operator of the form
\begin{equation}\label{tree}
(B_\lambda\varphi)(v)=\sum_{u\in\Chi(v)}\lambda_u\varphi(u)
\quad\text{for every }\varphi\in X\text{ and }v\in V,
\end{equation}
where $\lambda:=(\lambda_v)_{v\in V}\in(\K\setminus\{0\})^V$ is a suitably chosen weight sequence. It is worth noting that weighted backward shifts on trees are not, in general, weighted composition operators, since the value at a vertex may depend on several of its children. This provides further motivation for studying them as a distinct class of operators.

A first systematic study of weighted shifts on sequence spaces over directed trees was given by Jab{\l}o\'nski, Jung and Stochel in \cite{JablonskiJungStochel}, and the first investigation of their dynamical properties appeared in \cite{MartinezAvendano}. Building on this work, Grosse-Erdmann and Papathanasiou established in the fundamental papers \cite{GrossePapChaotic,Gro1} a broad basis for the dynamical study of these operators, addressing boundedness, hypercyclicity, weak mixing, mixing, and Devaney chaos. More recent work has shown that shifts on trees are also effective in the study of other dynamical questions (see for example \cite{AbakumovAbbar, baranov, chen, ernst}). In \cite{LopezPap}, weighted shifts on unrooted directed trees were used to solve an important problem on chain recurrence posed in \cite{BerPer}. Particularly significantly, the relevant phenomenon can be realized by a comparatively simple shift on a tree, whereas it cannot occur for classical unilateral or bilateral weighted shifts. This example illustrates one of the main advantages of working with shifts on trees. Besides the weights and the underlying sequence space, one can exploit the geometry and branching structure of the tree. This additional flexibility can lead to simpler examples and counterexamples, especially when one seeks operators that satisfy one dynamical property but fail another. Thus, shifts on trees are useful not only as generalizations of classical shifts, but also as tools for clarifying the relations between different notions in Linear Dynamics.

In parallel, recent works have studied Li--Yorke chaos, mean Li--Yorke chaos, and distributional chaos for weighted shifts acting on Banach and Fréchet sequence spaces indexed by $\N$ or $\Z$ (see \cite{BerBonPin26,Fernanda,joaojmaa}). It is therefore natural to extend these results from the classical chains to more complex directed trees. The purpose of the present paper is precisely to study and characterize these three notions of chaos for weighted backward shifts on Fréchet sequence spaces over directed trees.

The paper is organized as follows. In Section~\ref{sec:prelim}, we introduce the necessary preliminaries concerning directed trees and sequence spaces over them. In Section~\ref{sec:LY}, we obtain a general characterization of Li--Yorke chaos for weighted shifts on Fréchet sequence spaces over directed trees. As an immediate corollary, we derive a criterion in the setting of Köthe sequence spaces $\lambda_p(A,V)$, where $p\in\{0\}\cup[1,\infty)$. Under additional assumptions, we obtain further characterizations and recover some previously known results in more specific settings. For example, in Proposition~\ref{p2} we establish an equivalence between Li--Yorke chaos and the non-equicontinuity of the family of operators $\{B_\lambda^n\}_{n\geq1}$. In Section~\ref{sec:meanLY}, we obtain a general characterization of mean Li--Yorke chaos for weighted shifts on Fréchet sequence spaces over directed trees. In the case of Banach sequence spaces over directed trees, we obtain a stronger characterization. In addition, several further results are established. For example, Proposition~\ref{prop:mean-LY-equiv} establishes the equivalence between mean Li--Yorke chaos and the mean sensitivity of the operator $B_\lambda$. In Section~\ref{sec:dist}, we investigate distributional chaos in the same settings as above.

\section{Preliminaries} \label{sec:prelim}

Throughout this paper, $\mathbb{K}$ denotes either the field $\mathbb{R}$ of real numbers or the field $\mathbb{C}$ of complex numbers, $\mathbb{Z}$ denotes the ring of integers, $\mathbb{N}$ denotes the set of all positive integers, and $\mathbb{N}_0 := \mathbb{N} \cup \{0\}$. We start our preliminaries by introducing the basic definitions on trees.
    \begin{definition}
        A \emph{directed tree} $(V,E)$ is a connected directed graph consisting of a countable set $V$ of \emph{vertices} and a set $E \subset (V \times V) \backslash \{(v,v) : v \in V\}$ of \emph{edges} that are required to satisfy the properties below.
	\begin{enumerate}[label=$(\roman*)$]
		\item $(V,E)$ has no cycles, that is, there are no paths leading back to a vertex already visited.
		\item With the possible exception of a single vertex $\root\in V$, for each other vertex $v \in V$ there is a unique $u\in V$ such that $(u,v)\in E$. In this case, $u$ is called the \emph{parent} of $v$ and denoted by $u=\Par(v)$. 
		\item For each $v\in V$, there is at least one $u\in V$ such that $(v,u)\in E$. In other words, the tree has no \emph{leaves}. When $(v,u)\in E$, the vertex $u$ is called a \emph{child} of $v$. The set of all the children of a vertex $v$ is denoted by $\Chi(v)$.
	\end{enumerate}
    	When it exists, the only vertex $\root\in V$ without a parent is called the \emph{root} of $V$, in this case we say that $V$ is a \emph{rooted directed tree}. If there is no such an element, we say that $V$ is \emph{unrooted}.
    \end{definition}
    \begin{figure}[H]
        \centering
        \includegraphics[width=0.8\linewidth]{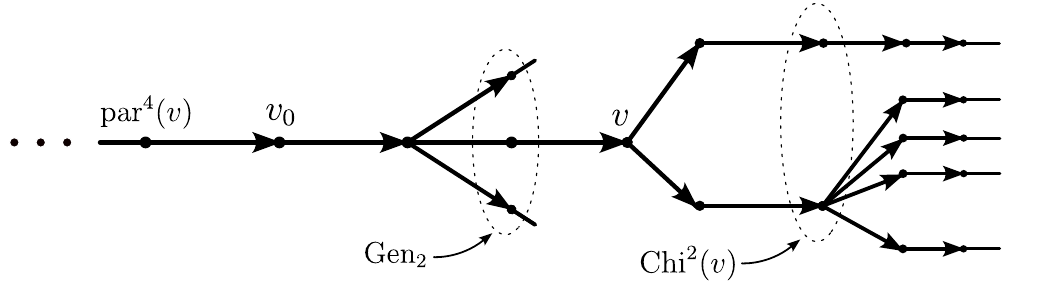}
        \caption{Representation of an unrooted directed tree and some of its elements}
        \label{fig:ex-tree-geomtria}
    \end{figure}
    In the remainder of this article, a directed tree $(V,E)$ will be simply denoted by $V$ and sometimes called only a \emph{tree}. Since we deal with iterates of operators, it will be useful to define the following notations. For each $v\in V$, we define $\Par^0(v)=v$, $\Chi^0(v)=\{v\}$ and, for all $n\geq 1$,
	\[\Par^n(v)=\Par(\Par^{n-1}(v)) \quad \text{and} \quad \Chi^n(v) = \bigcup_{u\in \Chi^{n-1}(v)}\Chi(u).\]
    It can be also useful to associate a generation to the elements of a directed tree $V$. We start by defining the following equivalence relation: we say that $u\sim v$ in $V$ when there is $n\in \N_0$ such that $\Par^n(u)=\Par^n(v)$. The equivalence class $[v]$ of any $v\in V$ is called \emph{the generation} of $v$. In order to count the generations, we must define what is the generation 0. We fix some $v_0\in V$, being $v_0=\root$ when $V$ is rooted, and we define
	\[\Gen_n=\begin{cases}
\Chi^n(\root), & n\geq 0, \\
\emptyset, & n< 0,
\end{cases}\]
	in the rooted case, and
	\[\Gen_n = \Gen_n(v_0) = \{v\in V : \exists k\geq 0 \text { with }\Par^{n+k}(v)=\Par^k(v_0) \}\quad n\in\Z,\]
	when $V$ is unrooted. The vertex $v_0$ used to count the generations will be also called the \emph{center} of the tree $V$. See Figure \ref{fig:ex-tree-geomtria} for a visual exemplification of some of these elements.

    \begin{example}
    The sets of natural numbers $\N$ and $\N_0$ are rooted trees, whereas the set of integers $\Z$ is an unrooted tree. As usual, we consider $\Gen_n:=\{n\}$ for $n\in J$, where $J=\N_0$ or $J=\Z$. Another simply defined and useful tree is the so called \emph{$N$-adic} (or \emph{$N$-ary}) \emph{tree}: it is the rooted or unrooted tree in which every vertex has exactly $N$ children. Both $\N$ and $\Z$ are particular cases of $1$-adic trees.
\end{example}
    
    For the next definition, recall that a vector space $X$ is said to be a {\it Fréchet space} if it is endowed with an increasing sequence $(\| \cdot \|_k)_{k \in \N}$ of seminorms (called a {\it fundamental sequence of seminorms}) that defines a metric 
	\begin{equation}\label{d}
	    d(x,y) := \sum_{k=1}^{\infty} \frac{1}{2^k}\,\min\{1,\|x-y\|_k\} \quad \text{for }x,y \in X,
	\end{equation}
	under which $X$ is complete. When we want to emphasize the underlying family of seminorms, we usually write the Fréchet space as $\left(X,(\|\cdot\|_n)_{n\in\N}\right)$. If $X$ is a Banach space endowed with a norm $\|\cdot\|$, we denote it by $(X,\|\cdot\|)$. Recall that the metric $d$ satisfies the following properties: 
	\begin{enumerate}[label=$(\mathrm{P}\arabic*)$]
		\item \label{metric:P1} $d(x_1+x_2,0)\le d(x_1,0)+d(x_2,0)$ for all $x_1,x_2\in X$.
		\item \label{metric:P2} $d(\lambda x,0)\le (1+|\lambda|)d(x,0)$ for all $\lambda\in\mathbb{K}$ and $x\in X$.
	\end{enumerate}
	
	\begin{definition}
	    Let $V$ be a directed tree. A Fréchet (Banach) space $X$ which is a vector subspace of the product space $\mathbb{K}^{V}:= \{ f : V \to \mathbb{K} : f \text{ is a function} \}$ is a {\it Fréchet (Banach) sequence space over} $V$ if the inclusion map $X \hookrightarrow \mathbb{K}^{V}$ is continuous, i.e., convergence in $X$ implies coordinatewise convergence.
\end{definition}
\begin{definition}
Let $V$ be a directed tree and let $X$ be a Fréchet sequence space over $V$. Fix a sequence $\lambda\in(\K\backslash \{0\})^V$, called a \emph{weight} and denoted by $(\lambda_v)_{v\in V}$. If the function $B_{\lambda}: X \to X$ defined by 
\begin{equation}\label{eqq0}
    (B_\lambda f)(v) = \sum_{u\in \Chi(v)}\lambda_uf(u)
\end{equation}
for all $v\in V$, is well-defined, that is, the series in (\ref{eqq0}) converges unconditionally and $B_{\lambda}f \in X$ for every $f\in X$, then we say that $B_{\lambda}$ is a \textit{weighted backward shift on} $X$. In this case, the closed graph theorem guarantees that $B_{\lambda}$ is continuous.
\end{definition}
 Since it can be useful, the iterates of $B_\lambda$ are given by
	\[(B_\lambda^n f)(v) = \sum_{u\in \Chi^n(v)} \lambda(v\to u)f(u),\]
	where 
    $$\lambda(v\to u) = \lambda_{\Par^{n-1}(u)}\cdots\lambda_{\Par^2(u)}\lambda_{\Par(u)}\lambda_u.$$
Fix a directed tree $V$. Consider the canonical vectors $(e_v)_{v\in V}$ given by
\[
e_v(u)=
\begin{cases}
1, & \text{if } u=v,\\
0, & \text{if } u\neq v.
\end{cases}
\]
Fix $v\in V$. Whenever $\operatorname{par}^n(v)$ does not exist for some $n\in \N$, we adopt the convention
\[
\lambda(\operatorname{par}^n(v)\to v):=0 \quad \text{and}\quad e_{\operatorname{par}^n(v)}:=0.
\]
Moreover, we adopt the usual conventions that empty sums are 0 and empty products are 1. In particular, $\lambda(v\to v):=1$ for any $v\in V$. Additionally, we convention that the supremum over the empty set is equal to zero.

\begin{definition}\label{defAd}
         Let $V$ be a directed tree and let $X$ be a Fréchet sequence space over $V$. Let the generations be
enumerated according to some vertex $v_0$. We say that $X$ is {\it admissible} if the following conditions hold:
         \begin{enumerate}[label=$(\mathrm{B}\arabic*)$, leftmargin=*] \setcounter{enumi}{-1}
            \item \label{frechet:condition- B0}  for any $v \in V$, $e_v \in X$;

\item \label{frechet:condition- B1} for any $f \in X$ and any $n\in \Z$, the series 
\[f\chi_{\Gen_n}=\sum_{v\in \Gen_n} f(v)e_v\]
converges unconditionally in $X$;
\item \label{frechet:condition- B2} for any $f \in X$, the series
\[
f=\sum_{n\in \Z} f\chi_{\Gen_n} 
\]
converges in $X$. 
\end{enumerate}
\end{definition}
Let $V$ be a directed tree. We denote by $c_{00}(V)$ the set of all functions $f:V\to \K$ such that $\{v\in V:f(v)\neq0\}$ is finite. Given an admissible Fréchet sequence space $X$ over $V$, it is easy to see that $c_{00}(V)$ is dense in $X$. By conditions \ref{frechet:condition- B1} and \ref{frechet:condition- B2}, the closed graph theorem, and the Banach-Steinhaus theorem, we obtain the following lemma.
\begin{lemma}\label{LM1}
    Let $V$ be a directed tree with generations enumerated according to some fixed vertex $v_0$. Let $X$ be an admissible Fréchet sequence space over $V$. Then the family of operators $(\Phi_{n})_{n\in \Z}$ given by
        $$\begin{matrix}
\Phi_n: & X & \to  &  X\\
 & g & \mapsto & g\chi_{\Gen_n}
\end{matrix}$$
is equicontinuous.
\end{lemma}
The next lemma follows from straightforward computations.
\begin{lemma}\label{LM2}
    Let $V$ be a directed tree with generations enumerated according to some fixed vertex $v_0$. Let $X$ be an admissible Fréchet sequence space over $V$. Suppose that the weighted backward shift $B_\lambda$, with nonzero weights $\lambda := (\lambda_v)_{v \in V}$, is well-defined and continuous on $X$. Then
     $$B_\lambda^j(f\chi_{\Gen_n})=\chi_{\Gen_{n-j}}B_\lambda^jf \quad \text{for all } \,f\in X,\, j\in\mathbb N_0, \text{ and }\, n\in \Z.$$
\end{lemma}

    Now, we will introduce a concrete class of examples of admissible Fréchet sequence spaces. To do this, we need the following definition.
\begin{definition}
		Let $V$ be a directed tree. A matrix $A = (a_{v,k})_{v\in V,\; k\in\mathbb{N}}$ is called a {\it  K\"othe matrix} if it satisfies:
		\begin{enumerate}[label=$(\roman*)$]
			\item $a_{v,k} \ge 0$ for all $v\in V$ and $k\in\mathbb{N}$;
			\item for each fixed $v\in V$, $ a_{v,k} \le a_{v,k+1}$ for all $k \in \N$;
			\item for each $v\in V$ there exists at least one $k\in \N$ such that $a_{v,k} > 0$.
		\end{enumerate}
	\end{definition}    
    \begin{definition}\label{def1}
		Let $V$ be a directed tree. Consider $p\in \{0\}\cup [1,\infty)$ and a Köthe matrix $A = (a_{v,k})_{v\in V,\; k\in\mathbb{N}}$. The {\it associated Köthe sequence space} (or simply {\it Köthe sequence space}) $\lambda_p(A,V)$ is the Fréchet sequence space defined as follows.
		\begin{itemize}
			\item If $p\in[1,\infty)$, then
			\[
			\lambda_p(A,V) := \Big\{ f \in \K^V : \sum_{v\in V} |a_{v,k}f(v)|^p < \infty, \, \forall k\in\mathbb{N} \Big\},
			\]
			endowed with the seminorms
			\[
			\|f\|_k := \left( \sum_{v\in V} |a_{v,k} f(v)|^p \right)^{1/p} 
			\quad\text{ for }\; f\in \lambda_p(A,V)\;\text{ and } k\in\mathbb{N}.
			\]
			
			\item If $p=0$, then
			\[
			\lambda_0(A,V) := \Big\{ f\in \K^V : \forall \varepsilon>0\text{ and }\forall k\in \N, \; \exists F\subset V\text{ finite s.t. } |a_{v,k}f(v)|< \varepsilon, \, \forall v\in V\setminus F \Big\},
			\]
			endowed with the seminorms
			\[
			\|f\|_k := \sup_{v\in V} |a_{v,k} f(v)|
			\quad \text{ for }\; f\in \lambda_0(A,V)\;\text{ and } k\in\mathbb{N}.
			\]
		\end{itemize}
	\end{definition}
    It is easy to see that the Köthe sequence spaces $\lambda_p(A,V)$, where $V$ is a directed tree and $p \in \{0\}\cup [1,\infty)$, are examples of admissible Fréchet sequence spaces.  
    \begin{example} \label{exlpc0}
         Let $V$ be a directed tree and let $(\nu_v)_{v\in V}\subset (0,\infty)$. Consider a Köthe matrix $A = (a_{v,k})_{v\in V,\; k\in\mathbb{N}}$ defined by $a_{v,k}=\nu_v$ for all $k\in \N$.
\begin{enumerate}[label=$(\alph*)$]
\item If $p=0$, then the Köthe space $\lambda_0(A,V)$ coincides with the Banach space
\[
c_0(\nu,V):=\left\{f\in\K^V:\forall \varepsilon>0,\ \exists F\subset V \text{ finite s.t. }|\nu_v f(v)|<\varepsilon\text{  }\forall v\in V\setminus F\right\},
\]
endowed with the norm $\|f\|=\sup_{v\in V} |\nu_v f(v)|$ for $f\in c_0(\nu,V)$. When $\nu=(1)_{v\in V}$, we simply denote $c_0(\nu,V)$ by $c_0(V)$.
\item If $p\in [1,\infty)$, then the Köthe space $\lambda_p(A,V)$ coincides with the Banach space
             \[
\ell^p(\nu,V):=\left\{f\in\K^V:\sum_{v\in V} |\nu_v f(v)|^p<\infty\right\},
\]
endowed with the norm $\|f\|=\left(\sum_{v\in V} |\nu_v f(v)|^p\right)^{1/p}$ for $f \in\ell^p(\nu,V)$. When $\nu=(1)_{v\in V}$, we simply denote $\ell^p(\nu,V)$ by $\ell^p(V)$.
\end{enumerate}
\end{example}
    
    \begin{example}\label{ExFNZ}
        Note that the usual Fréchet sequence spaces over $J:=\N$ or $J:=\Z$, whose canonical vectors $(e_n)_{n\in J}$ form a Schauder basis are examples of admissible Fréchet sequence spaces. Moreover, the usual weighted backward shift
\[
B_w((x_n)_{n\in J})=(w_{n}x_{n+1})_{n\in J} \quad \text{for }(x_n)_{n\in J}\in X,
\]
where $w:=(w_n)_{n\in J}\subset \K\setminus \{0\}$ is chosen so that $B_w$ is well defined on $X$, coincides with the weighted backward shift introduced in \eqref{eqq0} after reindexing the weights by $w_n=\lambda_{n+1}$.
    \end{example}	

\section{Li--Yorke chaos} \label{sec:LY}
    
    The term `chaos' was first introduced into the mathematical literature by Li and Yorke \cite{16} in their seminal study of the dynamics of interval maps. Intuitively, Li--Yorke chaos describes the existence of a large set of points such that, pairwise, their orbits repeatedly come arbitrarily close to each other and also separate from each other over time. Below, we recall the corresponding formal definition in the setting of metric spaces.
    
\begin{definition}
    Let $M$ be a metric space with metric $d$. A continuous map
$f\colon M\to M$ is said to be ({\it densely}) \textit{Li--Yorke chaotic} if there is an 
uncountable ({\it dense}) set $S\subset M$ such that
each pair $(x,y)\in S\times S$ of distinct points is a \textit{Li--Yorke pair} for
$f$ in the sense that
\[
\liminf_{n\to\infty} d(f^n(x),f^n(y))=0
\quad \text{and} \quad
\limsup_{n\to\infty} d(f^n(x),f^n(y))>0.
\]
In this case, $S$ is said to be a \emph{scrambled set}.
\end{definition}
In the linear setting, \cite{BBMP2} provides the following result, which we shall use in our propositions.
\begin{lemma}[{\cite[Theorems 9 and 15]{BBMP2}}] \label{LLY1}
Let $T:X\to X$ be a continuous linear operator on a separable Fréchet space $X$. Then, the following assertions are equivalent:
\begin{enumerate}[label=$(\roman*)$]
\item \label{LLY1:i} $T$ is Li--Yorke chaotic;
\item \label{LLY1:ii} $T$ admits a \emph{semi-irregular vector}, that is, a vector $x \in X$ such that the sequence $(T^n x)_{n\in\mathbb{N}}$ does not converge to zero but has a subsequence converging to zero;
\item \label{LLY1:iii} $T$ admits an \emph{irregular vector}, that is, a vector $x \in X$ such that the sequence $(T^n x)_{n\in\mathbb{N}}$ is unbounded but has a subsequence converging to zero;
\item \label{LLY1:iv} $T$ satisfies the \emph{Li--Yorke chaos criterion}, that is, there exists a subset $X_0$ of $X$ with the following properties:
    \begin{enumerate}[label=$(\alph*)$]
    \item \label{LLY1:iv:a} for every $x \in X_0$, $(T^n x)_{n\in \N}$ has a subsequence converging to zero;
    \item \label{LLY1:iv:b} there is a bounded sequence $(a_n)_{n\in \N}$ in $\overline{\Span(X_0)}$ such that the sequence $(T^n a_n)_{n\in \N}$ is unbounded.
    \end{enumerate}
\end{enumerate}
\end{lemma}
In the context of dense Li--Yorke chaos, we will use the following result.
\begin{lemma}[{\cite[Theorem 17]{BBMP2}}] \label{l01}
    Let $T:X\to X$ be a continuous linear operator on a separable Fréchet space $X$. Then $T$ is densely Li--Yorke chaotic if and only if there exists a dense subset $X_0\subset X$ with properties \ref{LLY1:iv:a} and \ref{LLY1:iv:b} of the Li--Yorke chaos criterion. In particular, if $T$ admits a dense set of irregular vectors, then $T$ is densely Li--Yorke chaotic.
\end{lemma}
The next theorem establishes a characterization of Li--Yorke chaos for weighted backward shifts on admissible Fréchet sequence spaces over directed trees.
\begin{theorem}\label{teo:carac:li-yorke}
    Let $V$ be a directed tree with generations 
enumerated according to some fixed vertex $v_0$. Let $\left( X, (\|\cdot\|_n)_{n\in \N} \right)$ be an admissible Fréchet sequence space over $V$. Suppose that the weighted backward shift $B_\lambda$, with nonzero weights $\lambda := (\lambda_v)_{v \in V}$, is well-defined and continuous on $X$. Then $B_\lambda$ is Li--Yorke chaotic if and only if there exist $S\subset \{n\in\mathbb Z:\Gen_n\neq\varnothing\}$, and, for each $n\in S$, a vector $f_n \in X$ supported in $\Gen_n$ such that the following conditions hold:
\begin{enumerate}[label=$(\Alph*)$]
    \item \label{teo:carac:li-yorke:A} there is a subsequence $(m_k)_{k \in \N}\subset \N$ such that, for all $n \in S$, we have $\lim_{k\to \infty}B_{\lambda}^{m_k}f_n=0$;
    \item \label{teo:carac:li-yorke:B} there exist $m \in \N$ and an increasing sequence $(N_k)_{k \in \N}$ of positive integers such that, for each $k \in \N$, there is a finite subset $F_k\subset S$ with $\left\| \sum_{j\in F_k}f_j\right\|_{\rho(k)}$ $\ne 0$ and
    $$\frac{\left\| \sum_{j\in F_k} B_{\lambda}^{N_k}f_j \right\|_m}{\left\| \sum_{j\in F_k}^{}f_j\right\|_{\rho(k)}}>k,$$
    where $\rho(k):= \max\{k,m\}$.
\end{enumerate}
\end{theorem}
	\begin{proof}
	    $(\Rightarrow):$ Suppose that $B_{\lambda}$ is Li--Yorke chaotic. Then, by Lemma \ref{LLY1}, $B_{\lambda}$ admits an irregular vector $f\in X$, that is, there exist $m \in \N$ and subsequences $(m_k)_{k\in \N}\subset \N$ and $(M_k)_{k\in \N}\subset \N$ such that 
        \begin{equation}\label{eq1}
            \lim_{k\to \infty}B_{\lambda}^{m_k}f=0 \quad \text{and} \quad \lim_{k\to \infty}\big\|  B_{\lambda}^{M_k}f \big\|_m=\infty.
        \end{equation}
        It is easy to see that we can take $m\in \N$ in (\ref{eq1}) sufficiently large such that $\left\| f \right\|_m \neq 0$. Fix such an $m\in\N$ and consider the definition of $\rho(k)$ given in the statement for this $m$. Define $S:=\left\{ n\in \Z:\Gen_n \neq \emptyset \right\}.$ By Lemma \ref{LM1}, the family of operators $(\Phi_{n})_{n\in \Z}$ given by 
        $$\Phi_{n}: g\in X \to g\chi_{\operatorname{Gen}_n}\in X \quad \text{for } n\in \Z,$$
        is equicontinuous. Thus, given a neighborhood $O\subset X$ of $0$, there exists a neighborhood $U\subset X$ of $0$  that satisfies 
 \begin{equation}\label{eq2}
     g\in U \quad \Rightarrow \quad g\chi_{\operatorname{Gen}_n}\in O\quad \text{ for all }n\in \Z.
 \end{equation}
 By (\ref{eq1}) and (\ref{eq2}) there exists $k_0 \in \N$ such that 
 $$k\ge k_0 \quad \Rightarrow \quad B_\lambda^{m_k}f\in U \quad \Rightarrow \quad (B_\lambda^{m_k}f)\chi_{\operatorname{Gen}_n}\in O$$
for all $n\in \Z$. By Lemma \ref{LM2}, $B_\lambda^{m_k}(f\chi_{\Gen_n})=\chi_{\Gen_{n-m_k}}B_\lambda^{m_k}f$  for all $k\in \N$ and $n\in S$. Therefore, item \ref{teo:carac:li-yorke:A} holds with $f_n:= f\chi_{\Gen_{n}}$ for $n\in S$. On the other hand, by (\ref{eq1}), for each $k\in \N$ there exists $N_k\in \N$ such that 
 \begin{equation}\label{eq3}
     \big\|  B_{\lambda}^{N_k}f \big\|_m> k (\|f\|_{\rho(k)}+1).
 \end{equation}
 Take the integers $N_k$, so that $(N_k)_{k\in\N}$ forms a subsequence of $\N$. Since $f=\sum_{n\in S} f_n$, and $B_{\lambda}$ is continuous, there exists a finite subset $F_k \subset S$ such that 
 \begin{equation}\label{eq4}
     \left| \frac{\left\|B_{\lambda}^{N_k}\left( \sum_{j\in F_k}^{}f_j \right) \right\|_m}{\left\| \sum_{j\in F_k}^{}f_j \right\|_{\rho(k)}} -\frac{\big\| B_{\lambda}^{N_k}f \big\|_m}{\| f \|_{\rho(k)}} \right|<\frac{1}{2\| f\|_{\rho(k)}}.
 \end{equation}
 Therefore, by (\ref{eq3}) and (\ref{eq4}), we obtain item \ref{teo:carac:li-yorke:B}. 

 $(\Leftarrow):$ By item \ref{teo:carac:li-yorke:A}, taking $X_0:=\left\{ f_n: n \in S \right\}$, item \ref{LLY1:iv:a} of the Li--Yorke chaos criterion is satisfied. Fix $m\in \N$ provided in item \ref{teo:carac:li-yorke:A}. For each $k \in \N$, consider the finite set $F_k\subset S$ given in item \ref{teo:carac:li-yorke:B} and define 
 $$g_{N_k}:= \frac{\sum_{j\in F_k}^{}f_j}{\left\| \sum_{j\in F_k}^{}f_j\right\|_{\rho(k)}}.$$
 If $n\notin \left\{ N_k:k \in \N \right\}$, then we put $g_n:=0$. The sequence $(g_k)_{k \in \mathbb{N}}$ is bounded. In fact, fix $n\in \N$. Thus, for all $k > \max \left\{ m,n \right\}$, we have 
 $$\left\| g_{N_k} \right\|_n =\frac{\left\| \sum_{j\in F_k}^{}f_j \right\|_n}{\left\| \sum_{j\in F_k}^{}f_j\right\|_{\rho(k)}}\le 1.$$
 By item \ref{teo:carac:li-yorke:B}, we have $\left\| B_{\lambda}^{N_k}g_{N_k} \right\|_m>k$ for all $k\in \N$. Thus, the sequence $(B_{\lambda}^{n}g_{n})_{n\in \N}$ is unbounded. Therefore, item \ref{LLY1:iv:b} of the Li--Yorke chaos criterion holds.
	\end{proof}
\begin{remark} \label{rmkLY} \
\begin{enumerate}[label=$(\alph*)$]
    \item In Theorem \ref{teo:carac:li-yorke}, for each $n\in S$, we chose a single vector $f_n \in X$ supported in $\Gen_n$, such that the family $(f_n)_{n\in S}$ satisfies conditions \ref{teo:carac:li-yorke:A} and \ref{teo:carac:li-yorke:B}. From the proof, we see that this was not necessary, that is, we could restate the theorem by taking, for each $n\in S$, not only a single vector, but rather a sequence $(f_{n,k})_{k\in \N}\subset X$ of vectors supported in $\Gen_n$, such that the family of vectors $\left\{ f_{n,k}:n\in S\text{ and }k\in \N \right\}$ satisfies conditions \ref{teo:carac:li-yorke:A} and \ref{teo:carac:li-yorke:B} of the theorem. 
    \item Item \ref{teo:carac:li-yorke:A} of Theorem \ref{teo:carac:li-yorke} is used in the proof of the implication $(\Leftarrow)$ in order to guarantee condition \ref{LLY1:iv:a} of the Li--Yorke chaos criterion. In this criterion, for each $x\in X_0$, the subsequence along which $(T^n x)_{n\in\N}$ converges to zero need not be the same; that is, it may depend on $x\in X_0$. Therefore, we may weaken condition \ref{teo:carac:li-yorke:A} of Theorem \ref{teo:carac:li-yorke} as follows:
\begin{enumerate}[label=$(\Alph*')$]
\item \label{teo:carac:li-yorke:A'} for each $n\in S$, the sequence $(B_{\lambda}^{k}f_n)_{k\in\N}$ admits a subsequence converging to zero.
\end{enumerate}
\end{enumerate}
\end{remark}
    As a consequence of Theorem \ref{teo:carac:li-yorke}, we obtain the following corollary, which characterizes Li--Yorke chaos in the context of weighted shifts on Köthe sequence spaces over trees. This corollary is obtained by directly applying Theorem \ref{teo:carac:li-yorke} to the seminorms from Definition \ref{def1}, and writing each $f_j$ (from Theorem \ref{teo:carac:li-yorke}) in terms of the canonical vectors $(e_v)_{v\in V}$. We observe that the next corollary can also be stated in the form mentioned in Remark~\ref{rmkLY}.
    
\begin{corollary} \label{cor:particular:Kothe}
Let $V$ be a directed tree with generations enumerated according to some fixed vertex $v_0$. Consider a Köthe sequence space $X := \lambda_p(A,V)$, where $A := (a_{v,k})_{v \in V, k \in \N}$ is a Köthe matrix and $p \in \{0\} \cup [1,\infty)$. Suppose that the weighted backward shift $B_\lambda$, with nonzero weights $\lambda := (\lambda_v)_{v \in V}$, is well-defined and continuous on $X$. Then $B_\lambda$ is Li--Yorke chaotic if and only if there exist a subset $S\subset \{n\in\mathbb Z:\Gen_n\neq\varnothing\}$, and a family of scalars $\left\{ c_{v}\in \K : v \in \Gen_n \text{for some } n\in S \right\}$ with $\sum_{v\in \operatorname{Gen}_n}^{}c_{v}e_{v}\in X$ for all $n\in S$, such that the following conditions hold:
\begin{enumerate}[label=$(\Alph*)$]
    \item \label{cor:particular:Kothe : (A)} there exists a subsequence $(m_k)_{k \in \N}\subset \N$ such that for all $n \in S$ we have
    $$\begin{cases}
\displaystyle\lim_{k\to \infty}\underset{v \in \operatorname{Gen}_{n-m_k}}{\sup} \Bigg|\sum_{u\in \operatorname{Chi}^{m_k}(v)}^{} \lambda(v\to u)c_{u} a_{v, \ell} \Bigg|=0\quad \text{ for all }\ell\in \N, & \text{if }p=0,\\[1mm]
\displaystyle\lim_{k\to \infty}\underset{v \in \operatorname{Gen}_{n-m_k}}{\sum} \Bigg|\sum_{u\in \operatorname{Chi}^{m_k}(v)}^{} \lambda(v\to u)c_{u} a_{v, \ell} \Bigg|^p=0\quad \text{ for all }\ell\in \N, & \text{if }p\in[1,\infty);
\end{cases}$$
\item \label{cor:particular:Kothe: (B)} there exist $m \in \N$ and an increasing sequence $(N_k)_{k \in \N}$ of positive integers such that for each $k \in \N$ there is a finite subset $F_k\subset S$ with 
    $$\underset{j\in F_k}{\max} \underset{u \in\operatorname{Gen}_{j}}{\sup} \left| c_{u}a_{u,\rho(k)} \right|\ne 0\,\,\text{ if }\,p=0,\quad \text{or}\quad \underset{j\in F_k}{\sum}\underset{u \in \operatorname{Gen}_{j}}{\sum} \left| c_{u}a_{u,\rho(k)} \right|^p \ne 0\,\,\text{ if }\,p\in[1,\infty),$$
and, taking $\rho(k):= \max\{k,m\}$,
    $$\begin{cases}
\displaystyle\frac{\underset{j\in F_k}{\max}\, \underset{v \in\operatorname{Gen}_{j-N_k}}{\sup}\left| \underset{u\in \operatorname{Chi}^{N_k}(v)}{\sum}\lambda(v\to u) c_{u}a_{v,m}\right|}{ \underset{j\in F_k}{\max} \underset{u \in\operatorname{Gen}_{j}}{\sup} \left| c_{u}a_{u,\rho(k)} \right|}>k, & \text{if }p=0,\\
\displaystyle\frac{ \underset{j\in F_k}{\sum}\,\underset{v \in \operatorname{Gen}_{j-N_k}}{\sum}\left| \underset{u\in \operatorname{Chi}^{N_k}(v)}{\sum}\lambda(v\to u) c_{u}a_{v,m}\right|^p}{\underset{j\in F_k}{\sum}\underset{u \in \operatorname{Gen}_{j}}{\sum} \left| c_{u}a_{u,\rho(k)} \right|^p}>k^p, & \text{if }p\in[1,\infty).
\end{cases}$$  
\end{enumerate}
\end{corollary}

\begin{example}\label{ex1}
Fix $N\in \N$ with $N\ge 2$. Consider the unrooted directed $N$-adic tree $V$ centered at the vertex $v_0$. Let $A=(a_{v,k})_{v\in V,\,k\in\N}$ be a Köthe matrix whose weights are constant on each generation, that is, assume that there exists a family $(b_{n,k})_{n\in\Z,\,k\in\N}$ of positive numbers such that
\begin{equation}\label{eqndia1}
    a_{v,k}=b_{n,k}\quad\text{for every }v\in\Gen_n.
\end{equation}
Let $p\in\{0\}\cup[1,\infty)$ and let $X:=\lambda_p(A,V)$. Fix $\alpha\in\K\setminus\{0\}$ and consider the weighted backward shift $B_\lambda$ with weights $\lambda_v=\alpha$ for every $v\in V$. We assume that, for every $k\in\N$, there exist $\ell\in\N$ and $C>0$ such that
\begin{equation}\label{eqndia2}
    b_{n-1,k}\leq Cb_{n,\ell}\quad\text{for every }n\in\Z.
\end{equation}
Under this assumption, $B_\lambda$ is well-defined and continuous on $X$. Define
\begin{equation}\label{eqndia3}
    C_p=
\begin{cases}
N|\alpha|, & p=0,\\
|\alpha|, & p=1,\\
N^{\frac{1}{p^*}}|\alpha|, & 1<p<\infty,
\end{cases}
\end{equation}
where $p^*=\frac{p}{p-1}$. We will prove that the following assertions are equivalent:
\begin{enumerate}[label=$(\roman*)$]
    \item \label{LY: Cp condition: (i)} $B_{\lambda}$ is Li--Yorke chaotic;
    \item \label{LY: Cp condition: (ii)} there exists $m\in \N$ such that, for all $\ell\in \N$,
    $$\limsup_{n\to\infty} \sup_{j\in\Z} C_p^n\frac{b_{j-n,m}}{b_{j,\ell}}=\infty.$$
\end{enumerate}

\ref{LY: Cp condition: (i)}$\Rightarrow$\ref{LY: Cp condition: (ii)}: We first note that if $f\in X$ is a nonzero vector supported on $\text{Gen}_n$ for some $n\in\Z$, then, using Hölder's inequality if $p\in [1,\infty)$, and triangle inequality if $p=0$, we obtain 
\begin{equation}\label{eqex1}
\frac{\big\|B_{\lambda}^{k}f\big\|_{\ell_1}}{\|f\|_{\ell_2}} \le C_p^{k}\frac{b_{n-k,\ell_1}}{b_{n,\ell_2}},\quad \text{for every }\ell_1,\ell_2,k\in \mathbb{N}.
\end{equation}
By condition \ref{cor:particular:Kothe: (B)} of Corollary \ref{cor:particular:Kothe}, there exist $S\subset \Z$, $m\in \N$, $(f_n)_{n\in S}\subset X$ where $f_n$ is supported on $\text{Gen}_n$ for every $n\in S$, and an increasing sequence $(N_k)_{k \in \N}$ of positive integers such that for each $k \in \N$ there is a finite subset $F_k\subset S$ with $\sum_{j\in F_k}^{}f_j\ne 0$ and 
\begin{equation}\label{eqex2}
    \begin{cases}
\displaystyle\frac{\max_{j\in F_k}\big\| B_{\lambda}^{N_k}f_j \big\|_m}{\max_{j\in F_k}\|f_j\|_{\rho(k)}}>k,& \text{if }p=0,\\
\displaystyle\frac{\sum_{j\in F_k}\big\| B_{\lambda}^{N_k}f_j \big\|_{m}^p}{\sum_{j\in F_k}\|f_j\|_{\rho(k)}^p}>k^p,& \text{if }p\in[1,\infty).
\end{cases}
\end{equation}
By (\ref{eqex2}), for each $k\in \N$ there exists $j_k\in F_k$ such that $\frac{\big\| B_{\lambda}^{N_k}f_{j_k}\big\|_{m}}{\| f_{j_k}\|_{\rho(k)}}> k$. Fix $\ell\in\N$ and take $k>\max\{\ell, m\}$. Then, by (\ref{eqex1}), we obtain  
$$C_p^{N_k}\frac{b_{j_k-N_k,m}}{b_{j_k,\ell}}>C_p^{N_k}\frac{b_{j_k-N_k,m}}{b_{j_k,\rho(k)}}>k.$$
Therefore, item \ref{LY: Cp condition: (ii)} holds.

\ref{LY: Cp condition: (ii)}$\Rightarrow$\ref{LY: Cp condition: (i)} : By item \ref{LY: Cp condition: (ii)} , there exist $m\in \N$ and a subsequence $(N_k)_{k\in \N}\subset \N$ such that, for each $k\in\N$, there is $j_k\in \Z$ such that 
\begin{equation}\label{eqex3}
    C_p^{N_k}\frac{b_{j_k-N_k,m}}{b_{j_k,\rho(k)}}>k.
\end{equation}
For this implication we use the version of Corollary \ref{cor:particular:Kothe} commented in Remark \ref{rmkLY}. For each $k\in \N$ choose $v_{k}\in \Gen_{j_k-N_k-1}$ and take two distinct vertices $v_{1,k},v_{2,k}\in\Chi(v_k).$ For the generation $\Gen_{j_k}$, define $c_u=1$ for $u\in \Chi^{N_k}(v_{1,k})$, $c_u=-1$ for $u\in \Chi^{N_k}(v_{2,k})$, and $c_u=0$ otherwise. Define $f_{j_k}:=\sum_{v\in \operatorname{Gen}_{j_k}}^{}c_{v}e_{v}$ for each $k\in \N$. Since $B_{\lambda}^{n}f_{j_k}=0$ for all $k\in \N$ and $n>N_k$, item \ref{cor:particular:Kothe : (A)} of Corollary \ref{cor:particular:Kothe} holds. By (\ref{eqex3}) we have 
$$\frac{\big\| B_{\lambda}^{N_k}f_{j_k} \big\|_m}{\|f_{j_k}  \|_{\rho(k)}}=C_p^{N_k}\frac{b_{j_k-N_k,m}}{b_{j_k,\rho(k)}}>k \quad \text{for every }k\in\N.$$
Therefore, item \ref{cor:particular:Kothe: (B)} of Corollary \ref{cor:particular:Kothe} holds with $F_k=\{j_k\}$. Thus, $B_\lambda$ is Li--Yorke chaotic.
\end{example}
For a continuous linear operator  $T$ on a Fréchet space $X$, we define
\[
\NS(T):=\{x\in X:(T^n x)_{n\geq 0} \text{ has a subsequence converging to }0\}.
\]
Bernardes and Vasconcellos proved that, if $\NS(T)$ is dense, then either $T^nx\to0$ for every $x\in X$ or $T$ has a residual set of irregular vectors (see \cite[Theorem 17]{Fernanda}). We use this result in the next proposition.
\begin{proposition} \label{prop:Fernanda17}
Let $V$ be a directed tree and $X$ an admissible Fréchet sequence space over $V$. Suppose that the weighted backward shift $B_\lambda$, with nonzero weights $\lambda := (\lambda_v)_{v \in V}$, is well-defined and continuous on $X$, and assume that the following condition holds.
\begin{enumerate}[label=$(\mathrm{C}1)$] 
    \item \label{cond:(C):li-yorke} There are $v_0\in V$ and a subsequence $(n_j)_{j\in \N}\subset \N$ such that
$$B_\lambda^{n_j}e_{v_0} =\lambda\bigl(\Par^{n_j}(v_0)\to v_0\bigr) e_{\Par^{n_j}(v_0)} \xrightarrow{j\to\infty}0.$$ 
\end{enumerate}
Then $\NS(B_\lambda)$ is dense in $X$, and exactly one of the following
properties holds:
\begin{enumerate}[label=$(\alph*)$]
\item \label{cor:cond:irregular-residual} $B_\lambda$ has a residual set of irregular vectors (and is, by Lemma \ref{l01}, densely Li--Yorke chaotic);
\item \label{cor:cond:Bfto0} $B_\lambda^n f\xrightarrow{n\to\infty}0$ for every $f\in X$.
\end{enumerate}
\end{proposition}
\begin{proof}
We show that $\NS(B_\lambda)$ is dense by verifying that $c_{00}(V)\subset\NS(B_\lambda)$. Let $F\subset V$ be finite. For each $v\in F$, consider $p_v,q_v\in\N_0$ such that $\Par^{p_v}(v_0)=\Par^{q_v}(v)$. For all $n\in \N$ large enough, we can write
\begin{align*}
B_\lambda^{n-p_v+q_v}e_v & = B_\lambda^{n-p_v}B_{\lambda}^{q_v}e_v\\
&=\lambda(\operatorname{par}^{q_v}(v)\to v)B_\lambda^{n-p_v}e_{\operatorname{par}^{p_v}(v_0)}\\
&=\overset{=:c_v}{\overbrace{\frac{\lambda(\operatorname{par}^{q_v}(v)\to v)}{\lambda(\operatorname{par}^{p_v}(v_0)\to v_0)}}}B_\lambda^{n-p_v}B_\lambda^{p_v}e_{v_0}\\
&=c_vB_{\lambda}^n e_{v_0}.
\end{align*}
Fix $N_0\in\N$ such that $N_0 + p_v-q_v\geq 0$ for all $v\in F$. Then, for all $j\in\N$ large enough and for all $v\in F$, we obtain
\begin{align*}
B_\lambda^{n_j+N_0}e_v &=B_\lambda^{(n_j-p_v+q_v)+(N_0+p_v-q_v)}e_v \\
&=B_\lambda^{N_0+p_v-q_v}B_\lambda^{n_j-p_v+q_v}e_v\\
&=c_vB_\lambda^{N_0+p_v-q_v}B_\lambda^{n_j}e_{v_0}
\end{align*}
Since, by condition \ref{cond:(C):li-yorke}, the expression in the last equality tends to zero as $j\to\infty$, we conclude that $c_{00}(V)\subset\NS(B_\lambda)$. Since $c_{00}(V)$ is dense in $X$, it follows that $\NS(B_\lambda)$ is dense. Therefore, the result follows from \cite[Theorem 17]{Fernanda}.
\end{proof}

In our next result, we establish a characterization of Li–Yorke chaos for weighted backward
shifts on admissible Fréchet sequence spaces over directed trees, under the assumption condition of \ref{cond:(C):li-yorke}.

\begin{proposition}\label{p2}
    Let $V$ be a directed tree. Let $\left( X, (\|\cdot\|_n)_{n\in \N} \right)$ be an admissible Fréchet sequence space over $V$. Suppose that the weighted backward shift $B_\lambda$, with nonzero weights $\lambda := (\lambda_v)_{v \in V}$, is well-defined and continuous on $X$. Suppose that condition \ref{cond:(C):li-yorke} holds. Then the following assertions are equivalent:
\begin{enumerate}[label=$(\roman*)$]
    \item \label{p2 : (i)} $B_{\lambda}$ is Li--Yorke chaotic;
    \item\label{p2 : (ii)} $B_{\lambda}$ is densely Li--Yorke chaotic;
    \item\label{p2 : (iii)} the family $\{B_\lambda^n\}_{n\geq 1}$ is not equicontinuous;
    \item \label{p2 : (iv)} there exists $k\in\N$ such that, for every $\ell\in\N$,
\[
\sup\{\|B_\lambda^n f\|_k : n\geq1, f\in X, \|f\|_\ell \leq 1 \}=\infty;
\]
    \item\label{p2 : (v)} $B_\lambda^n f\not\rightarrow 0$ for some $f\in X$;
    \item\label{p2 : (vi)} there exist $m \in \N$ and an increasing sequence $(N_k)_{k \in \N}$ of positive integers such that, for each $k \in \N$, there is a finite set $F_k\subset V$ and scalars $(c_v)_{v\in F_k}$ with $\left\| \sum_{v\in F_k} c_{v}e_v\right\|_{\rho(k)}$ $\ne 0$ and
    \begin{equation}\label{eqply}
        \frac{\left\| \sum_{v\in F_k}\lambda(\Par^{N_k}(v)\to v) c_{v}e_{\Par^{N_k}(v)} \right\|_m}{\left\| \sum_{v\in F_k} c_{v}e_v\right\|_{\rho(k)}}>k,
    \end{equation}
    where $\rho(k):= \max\{m,k\}$.
\end{enumerate}
\end{proposition}
    \begin{proof}
\ref{p2 : (i)}$\Rightarrow$\ref{p2 : (ii)}: By Lemma \ref{LLY1}, $B_{\lambda}$ admits an irregular vector $f\in X$. By the proof of Proposition \ref{prop:Fernanda17}, we have $c_{00}(V)\subset\NS(B_\lambda)$. Therefore, applying Lemma \ref{l01} with the set $X_0:=c_{00}(V)$ and the constant sequence $(a_n)_{n\in\N}:=(f)_{n\in\N}$, we conclude that item \ref{p2 : (ii)} holds.

\ref{p2 : (ii)}$\Rightarrow$\ref{p2 : (iii)}: By Lemma \ref{LLY1}, $B_{\lambda}$ admits an irregular vector $f\in X$. Suppose that the family $\{B_\lambda^n\}_{n\geq 1}$ is equicontinuous. Thus, for each $m\in\N$ there exist $\ell\in\N$ and $C>0$ such that
$$\|B_\lambda^n f\|_m\leq C\|f\|_\ell\quad \text{for all }n\in\N.$$
Then the orbit of $f\in X$ is bounded, which is a contradiction.This proves the implication.

\ref{p2 : (iii)}$\Rightarrow$\ref{p2 : (iv)}: This implication follows from the definition of equicontinuity.

\ref{p2 : (iv)}$\Rightarrow$\ref{p2 : (v)}: This implication follows from the Banach-Steinhaus theorem.

\ref{p2 : (v)}$\Rightarrow$\ref{p2 : (vi)}: By Proposition \ref{prop:Fernanda17}, $B_{\lambda}$ is Li--Yorke chaotic. Thus, we conclude that condition \ref{teo:carac:li-yorke:B}  of Theorem \ref{teo:carac:li-yorke} holds. Hence, using condition \ref{frechet:condition- B1} of Definition \ref{defAd} and the continuity of $B_{\lambda}$ and of the seminorms, it follows that condition \ref{teo:carac:li-yorke:B}  of Theorem \ref{teo:carac:li-yorke} implies item \ref{p2 : (vi)}.

\ref{p2 : (vi)}$\Rightarrow$\ref{p2 : (i)}: Consider $X_0:= c_{00}(V)$. By the proof of Proposition \ref{prop:Fernanda17}, $c_{00}(V)\subset \NS(B_{\lambda})$. Thus, condition \ref{LLY1:iv:a} of the Li--Yorke chaos criterion holds. Defining 
$$g_n:= \begin{cases}
\displaystyle \frac{\sum_{v\in F_k} c_{v}e_v}{\left\| \sum_{v\in F_k} c_{v}e_v \right\|_{\rho(k)}},& \text{if }n=N_k\text{ for some }k\in \mathbb{N},\\
0,& \text{otherwise,}
\end{cases}$$
and arguing as in the proof of implication $(\Leftarrow)$ of Theorem \ref{teo:carac:li-yorke}, we conclude that condition \ref{LLY1:iv:b} of the Li--Yorke chaos criterion, holds. Therefore, $B_{\lambda}$ is Li--Yorke chaotic.
    \end{proof}
    If $V$ is rooted, then condition \ref{cond:(C):li-yorke} is automatically satisfied. Hence, we obtain the following corollary.
\begin{corollary}
    Let $V$ be a rooted directed tree. Let $X$ be an admissible Fréchet sequence space over $V$. Suppose that the weighted backward shift $B_\lambda$, with nonzero weights $\lambda := (\lambda_v)_{v \in V}$, is well-defined and continuous on $X$. Then assertions \ref{p2 : (i)}-\ref{p2 : (vi)} from Proposition \ref{p2} are equivalent.
\end{corollary}

Our results allows us to recover the following theorem from the recent literature.    
\begin{corollary}[{\cite[Theorem 22]{Fernanda}}] \label{cor:carac:li-yorke}
Let $X$ be a Fréchet sequence space over $\Z$ in which $(e_n)_{n\in\Z}$ is a Schauder basis. Suppose that the weighted backward shift $B_w$, with nonzero weights $w := (w_n)_{n \in \Z}$, is well-defined and continuous on $X$. Then $B_w$ is Li--Yorke chaotic if and only if the following conditions hold: 
\begin{enumerate}[label=$(\alph*)$]
    \item \label{cor:carac:li-yorke:(a)} the sequence $(w_{-n}\cdots w_{-1} e_{-n})_{n\geq 1}$ has a subsequence converging to $0\in X$;
    \item\label{cor:carac:li-yorke:(b)} there is a vector $x=(x_n)_{n\in \Z} \in X$ such that 
    \[(w_{j}\cdots w_{j+n-1}x_{j+n})_{j\in\Z}\not\rightarrow 0\in X \text{ as }n\to\infty.\]
\end{enumerate}
Moreover, whenever these conditions hold, $B_w$ is densely Li--Yorke chaotic.
\end{corollary}
\begin{proof}
Recall that, by Example \ref{ExFNZ}, we may regard $B_w$ as the weighted backward shift $B_\lambda$ introduced in \eqref{eqq0} on the admissible Fréchet sequence space over $V=\Z$. For this, it suffices to define $\lambda_{n+1}:=w_n$ for every $n\in\Z$.

($\Rightarrow$): Since $B_{\lambda}$ is Li--Yorke chaotic, by Lemma \ref{LLY1}, $B_{\lambda}$ admits an irregular vector and, hence, condition \ref{cor:carac:li-yorke:(b)} holds. On the other hand, by Theorem \ref{teo:carac:li-yorke}, there exist $n\in\Z$ and $c\in \K\setminus\{0\}$ such that $\liminf_{k\to \infty}B_{\lambda}^k(ce_n)=0$. Therefore, item \ref{cor:carac:li-yorke:(a)} holds.

($\Leftarrow$): Condition \ref{cor:carac:li-yorke:(a)} implies condition \ref{cond:(C):li-yorke}. Condition \ref{cor:carac:li-yorke:(b)} implies that there exists $f\in X$ such that $B_\lambda^n f\not\rightarrow 0$. Therefore, by Proposition \ref{p2}, $B_{\lambda}$ is Li--Yorke chaotic.
\end{proof}
   Assuming condition \ref{cond:(C):li-yorke} in the context of Köthe sequence spaces, we obtain the following characterization.
\begin{corollary} \label{cor:kothe-tree:outro}
Let $V$ be a directed tree and let $X:=\lambda_p(A,V)$, where $A=(a_{v,k})_{v\in V,\,k\in\N}$ is a Köthe matrix with non-zero entries and $p\in \{0\}\cup [1,\infty)$. Suppose that the weighted backward shift $B_\lambda$, with nonzero weights $\lambda := (\lambda_v)_{v \in V}$, is well-defined and continuous on $X$. Suppose that there is a subsequence $(n_j)_{j\in \N}\subset \N$ such that, for some $v_0\in V$,
\begin{equation}\label{koethec1}
    \left|\lambda\bigl(\Par^{n_j}(v_0)\to v_0\bigr)\right|
a_{\Par^{n_j}(v_0),k}\xrightarrow{j\to\infty}0
\quad\text{for every }k\in\N.
\end{equation}
Then the following assertions are equivalent:
\begin{enumerate}[label=$(\roman*)$]
    \item \label{cor:kothe-tree:outro i} $B_\lambda$ is Li--Yorke chaotic;
    \item\label{cor:kothe-tree:outro ii} $B_\lambda$ is densely Li--Yorke chaotic;
    \item\label{cor:kothe-tree:outro iii} there exists $m\in\N$ such that, for every $\ell\in\N$,
    \begin{equation}\label{eqcorKo}
    \begin{cases}
\displaystyle \sup_{\substack{n\geq 1\\ v\in V}} \sum_{u\in \operatorname{Chi}^n(v)} \left| \frac{a_{v,m}\lambda(v\to u)}{a_{u,\ell}} \right|=\infty, &\text{if }\, p=0,\\
\displaystyle \sup_{\substack{n\geq 1\\ v\in V}} \sup_{u\in \operatorname{Chi}^n(v)} \left| \frac{a_{v,m}\lambda(v\to u)}{a_{u,\ell}} \right|=\infty, &\text{if }\,  p=1,\\
\displaystyle \sup_{\substack{n\geq 1\\ v\in V}} \left( \sum_{u\in \operatorname{Chi}^n(v)} \left| \frac{a_{v,m}\lambda(v\to u)}{a_{u,\ell}} \right|^{p^*} \right)^{1/p^*}=\infty, &\text{if }\,  1<p<\infty.
\end{cases}
\end{equation}
\end{enumerate}
\end{corollary}
\begin{proof}
First, note that \eqref{koethec1} is equivalent to condition \ref{cond:(C):li-yorke}.  Using Hölder's inequality in the case $p\in [1,\infty)$ and the triangle inequality in the case $p=0$, we obtain
$$\sup_{\substack{n\geq1 \\ \|f\|_\ell\leq1}}\|B_\lambda^n f\|_m$$
is equal to the expression in (\ref{eqcorKo}). Therefore, the equivalences come from Proposition \ref{p2}.
\end{proof}
Let $V$ be a directed tree. For a fixed \(p \in [1,\infty)\), recall that, by Example \ref{exlpc0}, \(\ell^p(V)=\lambda_p(A,V)\), where \(A=(a_{v,k})_{v\in V,\,k\in\N}\) is given by \(a_{v,k}=1\) for all \(v\in V\) and \(k\in\N\). Likewise, \(c_0(V)=\lambda_0(A,V)\). Therefore, as we shall see in the following result, by Corollary~\ref{cor:kothe-tree:outro}, we get a more general characterization than the one obtained in \cite[Corollary 3.9]{Ko} on Li--Yorke chaos for weighted shifts on \(c_0(V)\) and \(\ell^p(V)\) that satisfy condition \ref{cond:(C):li-yorke}.
\begin{corollary}\label{CoKo}
    Let $V$ be a directed tree and let $X:=\ell^p(V)$ $(1\le p<\infty)$ or $X:=c_0(V)$. Suppose that the weighted backward shift $B_\lambda$, with nonzero weights $\lambda := (\lambda_v)_{v \in V}$, is well-defined and continuous on $X$. Moreover, suppose that 
    $$\liminf_{n\to \infty} |\lambda(\operatorname{par}^n(v)\to v)|=0\quad \text{for some }v\in V.$$
    Then the following assertions are equivalent:
    \begin{enumerate}[label=$(\roman*)$]
        \item \label{CoKo i} $B_\lambda$ is Li--Yorke chaotic on $X$;
        \item \label{CoKo ii} $B_\lambda$ is densely Li--Yorke chaotic on $X$;
        \item \label{CoKo iii} we have 
        $$\begin{cases}
\displaystyle \sup_{n\in\N,\,v\in V}\abs{\lambda(\Par^n(v)\to v)}=\infty, & \text{if } X=\ell^1(V),\\[6pt]
\displaystyle \sup_{n\in\N,\,v\in V}\sum_{u\in\Chi^n(v)}\abs{\lambda(v\to u)}^{p^*}=\infty, & \text{if } X=\ell^p(V)\ (1<p<\infty),\\[6pt]
\displaystyle \sup_{n\in\N,\,v\in V}\sum_{u\in\Chi^n(v)}\abs{\lambda(v\to u)}=\infty, & \text{if } X=c_0(V).
\end{cases}$$
    \end{enumerate}
    
\end{corollary}
In the next corollary, using our previous results, we can recover the characterization in \cite[Theorem 24]{Fernanda} in a much simpler way than in the original proof.
\begin{corollary}[{\cite[Theorem 24]{Fernanda}}] \label{cor:Nilson-Fernanda24}
Let $X:=\lambda_p(A,\Z)$, where $A=(a_{j,k})_{j\in\Z,\,k\in\N}$ is a Köthe matrix with nonzero entries and $p\in\{0\}\cup[1,\infty)$. Let $w=(w_j)_{j\in\Z}$ be a sequence of nonzero scalars such that the bilateral weighted backward shift $B_w$ is well-defined and continuous on $X$. Then $B_w$ is Li--Yorke chaotic if and only if the following conditions hold.
\begin{enumerate}[label=$(\alph*)$]
\item \label{cor:Nilson-Fernanda24:a} There exists a strictly increasing sequence $(n_k)_{k\in\N}$ of positive integers such that
\[
\lim_{k\to\infty} a_{-n_k,\ell}w_{-n_k}\cdots w_{-1}=0 \quad\text{for every }\ell\in\N.
\]
\item \label{cor:Nilson-Fernanda24:b} There exists $m\in\N$ such that, for every $\ell\in\N$,
\[
\sup
\left\{\frac{a_{i,m}|w_i\cdots w_j|}{a_{j+1,\ell}} : i,j\in\Z, i\leq j\right\}=\infty.
\]
\end{enumerate}
\end{corollary}

\begin{proof}
As observed in Example \ref{ExFNZ}, we may regard $B_w$ as the weighted backward shift $B_\lambda$ introduced in \eqref{eqq0} on $X$, by taking $\lambda_{n+1}:=w_n$ for every $n\in\Z$. Moreover, as noted in the proof of Corollary \ref{cor:kothe-tree:outro}, we have
\begin{equation} \label{eq1:cor:NF24}
\sup_{\substack{\|f\|_\ell\leq1 \\ n\geq1}}\|B_w^nf\|_k=\sup_{\substack{i,j\in\Z\\i\leq j}}\frac{a_{i,k}|w_i\cdots w_j|}{a_{j+1,\ell}}\quad \text{for all }k,\ell\in\N.
\end{equation}

($\Rightarrow$): By Corollary \ref{cor:particular:Kothe}, there exist $n\in\Z$, $c\in \K\setminus\{0\}$ and a subsequence $(m_k)_{k\in \N}\subset \N$ such that
\[
\lim_{k\to\infty}\left|cw_{n-m_k}\cdots w_{n-1} a_{n-m_k,\ell}\right|=\lim_{k\to\infty}\left| c\lambda(n-m_k\to n)a_{n-m_k,\ell} \right|=0\quad \text{ for all }\ell\in\N.
\]
Thus, condition \ref{cor:Nilson-Fernanda24:a} holds. Moreover, since condition \ref{cor:Nilson-Fernanda24:a} is equivalent to condition \ref{cond:(C):li-yorke} and $B_w$ is Li--Yorke chaotic, by Proposition \ref{p2} and \eqref{eq1:cor:NF24}, condition \ref{cor:Nilson-Fernanda24:b} follows.
		
($\Leftarrow$): As noted previously, condition \ref{cor:Nilson-Fernanda24:a} is equivalent to condition \ref{cond:(C):li-yorke}. Therefore, by Proposition \ref{p2}, condition \ref{cor:Nilson-Fernanda24:b} and \eqref{eq1:cor:NF24}, $B_w$ is Li--Yorke chaotic.
\end{proof}
For the next example, we need the following lemma, which characterizes hypercyclicity in the context of shifts on spaces of sequences over rooted trees. Recall that an operator $T:X \to X$ is said to be {\it hypercyclic} if there exists $x \in X$ whose orbit is dense, i.e., $\overline{\left\{ T^nx:n\in \N_0 \right\}}= X$. 
\begin{lemma}\cite[Theorem 4.4]{Gro1}\label{legro}
    Let $V$ be a rooted directed tree. Let $X=\ell^p(V)$, $1\leq p<\infty$, or $X=c_0(V)$, and suppose that the weighted backward shift $B_\lambda$, with nonzero weights $\lambda=(\lambda_v)_{v\in V}$, is well-defined and continuous on $X$. Then $B_\lambda$ is hypercyclic if and only if there is an increasing sequence of positive integers $(n_k)_{k\in \N}$ such that, for every $v\in V$,
\[
\begin{cases}
\displaystyle \sup_{u\in\Chi^{n_k}(v)} \abs{\lambda(v\to u)}\xrightarrow{k\to\infty}\infty, & \text{if } X=\ell^1(V),\\[6pt]
\displaystyle \sum_{u\in\Chi^{n_k}(v)} \abs{\lambda(v\to u)}^{p^*}\xrightarrow{k\to\infty}\infty, & \text{if } X=\ell^p(V),\ 1<p<\infty,\\[6pt]
\displaystyle \sum_{u\in\Chi^{n_k}(v)} \abs{\lambda(v\to u)}\xrightarrow{k\to\infty}\infty, & \text{if } X=c_0(V).
\end{cases}
\]
\end{lemma}
Examples of shifts that are Li--Yorke chaotic but not hypercyclic are already known; see, for instance, \cite[Page 16]{Fernanda}. In what follows,  we will give an example of a weighted backward shift on spaces of the form $\ell^p(V)$, $1\leq p<\infty$, or $c_0(V)$, where $V$ is a directed tree, that is Li--Yorke chaotic, but is not hypercyclic. The reader may observe how much simpler our example is than the one given in \cite{Fernanda}. This illustrates the power of studying shifts on sequences indexed by trees.

\begin{figure}[H]
        \centering
        \includegraphics[height=0.19\linewidth]{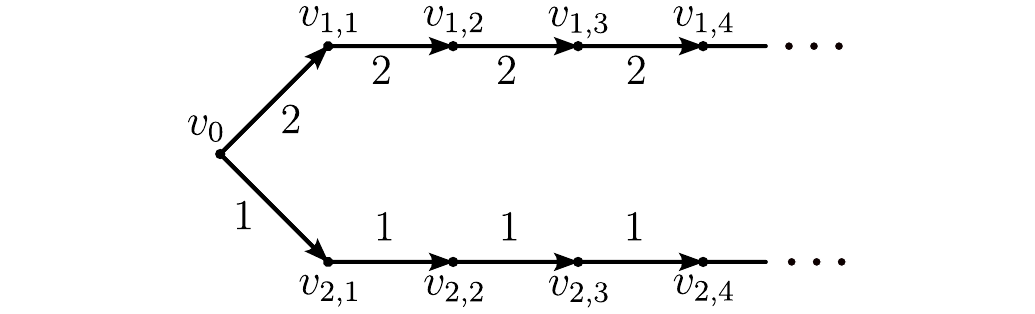}
        \caption{Directed tree in Example \ref{exNHcY}}
        \label{fig:2}
\end{figure}
    \begin{example}\label{exNHcY}
        Consider the directed tree $V$ with root $v_0$, with $\Chi(v_0)=\{v_{1,1},v_{2,1}\}$ and $\Chi(v_{i,j})=\{v_{i,j+1}\}$ for all $i=1,2$ and $j\in \N$. See Figure \ref{fig:2} for a representation of $V$. Let $\lambda = (\lambda_v)_{v\in V}$ be the weight sequence represented over the edges of $V$ in the figure.
     Consider $X=\ell^p(V)$, $1\leq p<\infty$, or $X=c_0(V)$, and note that the weighted backward shift $B_{\lambda}$ is well-defined and continuous on $X$. Note that, for all $n\in \N$, we have 
     $$\begin{cases}
\displaystyle \sup_{u\in\Chi^{n}(v_{2,1})} \abs{\lambda(v_{2,1}\to u)}=1, & \text{if } X=\ell^1(V),\\[6pt]
\displaystyle \sum_{u\in\Chi^{n}(v_{2,1})} \abs{\lambda(v_{2,1}\to u)}^{p^*}=1, & \text{if } X=\ell^p(V),\ 1<p<\infty,\\[6pt]
\displaystyle \sum_{u\in\Chi^{n}(v_{2,1})} \abs{\lambda(v_{2,1}\to u)}=1, & \text{if } X=c_0(V).
\end{cases}$$
Therefore, by Lemma \ref{legro}, $B_{\lambda}$ is not hypercyclic. On the other hand, since
$$\abs{\lambda(\operatorname{par}^n(v_{1,n+1})\to v_{1,n+1})}=2^n\quad \text{ for every }n\in \N,$$
by Corollary \ref{CoKo}, we conclude that $B_{\lambda}$ is Li--Yorke chaotic.
    \end{example}
    
\section{Mean Li--Yorke chaos} \label{sec:meanLY}

The formal concept of mean Li--Yorke chaos was introduced in the setting of topological dynamical systems by Huang, Li and Ye in \cite{Hu}, who investigated this form of chaos in positive entropy systems. We now state this definition.
\begin{definition}\label{lmdef}
		Let $(X,d)$ be a metric space and let $f:X\to X$ be a continuous map.
		A pair $(x,y)\in X\times X$, $x\neq y$, is called a \emph{mean Li--Yorke pair}
		for $f$ if
		\[
		\liminf_{n\to\infty} \frac{1}{n} \sum_{k=1}^{n} d\big(f^k(x), f^k(y)\big) = 0
		\quad\text{and}\quad
		\limsup_{n\to\infty} \frac{1}{n} \sum_{k=1}^{n} d\big(f^k(x), f^k(y)\big) > 0.
		\]
		The function $f$ is said to be ({\it densely}) \emph{mean Li--Yorke chaotic}
		if there exists an uncountable (dense) set $S\subset X$ such that every
		pair $(x,y)\in S\times S$ of distinct points is a mean Li--Yorke pair.
	\end{definition}
    
In the setting of the preceding definition, we recall some important concepts related to mean Li--Yorke chaos.
	\begin{itemize}
			\item We say that $f$ is \emph{mean sensitive} if there exists $\delta>0$ such that, for every $x\in X$ and every $\varepsilon>0$, there exists $y\in X$ with $d(x,y)<\varepsilon$ and
			\begin{equation}\label{meansens}
			    \limsup_{n\to\infty}\frac{1}{n}\sum_{k=1}^{n} d(f^k( x),f^k( y))\ge \delta.
			\end{equation}
			\item For $\delta>0$, a pair $(x,y)\in X\times X$ is called a \emph{mean Li--Yorke $\delta$-chaotic pair} if
			\[
			\liminf_{n\to\infty}\frac{1}{n}\sum_{k=1}^{n} d(f^k (x),f^k( y))=0
			\quad \text{and} \quad
			\limsup_{n\to\infty}\frac{1}{n}\sum_{k=1}^{n} d(f^k( x),f^k (y))\ge \delta.
			\]
			We say that $f$ is \emph{mean Li--Yorke sensitive} if there exists $\delta>0$ such that, for every $x\in X$ and every $\varepsilon>0$, there exists $y\in X$ with $d(x,y)<\varepsilon$ such that $(x,y)$ is a mean Li--Yorke $\delta$-chaotic pair. 
		\end{itemize}
    Observe that every mean Li--Yorke sensitive system is mean sensitive. From this point on, given a Fréchet space $\big(X,(\|\cdot\|_n)_{n\in\mathbb{N}}\big)$, we will always consider the \emph{compatible metric} $d$ to be the one given in \eqref{d}. In this section, we will use the theory of mean Li--Yorke chaos for linear operators on Fréchet spaces developed in \cite{JiangLi} to establish our characterizations. Before stating the results from \cite{JiangLi} that we will use, we need the following definitions, for which we fix a Fréchet space $X$ with compatible metric $d$ and a continuous linear operator $T$ on $X$.
    \begin{itemize}
        \item A vector $x\in X$ is called \textit{absolutely mean} \textit{semi-irregular} (or {\it mean semi-irregular point}) if
		$$\liminf_{n\to \infty}\frac{1}{n}\sum_{k=1}^{n}d(T^kx,0)=0\quad \text{and}\quad \limsup_{n\to \infty}\frac{1}{n}\sum_{k=1}^{n}d(T^kx,0)>0.$$
        \item The \textit{mean asymptotic cell} and the \textit{mean proximal cell} of $0$ are defined respectively by 
		\begin{equation}\label{MAsym}
		    \text{MAsym}(T,0):=\left\{ x \in X: \lim_{n\to \infty} \frac{1}{n}\sum_{k=1}^{n}d(T^kx,0)=0 \right\}\quad \text{and}
		\end{equation}
		\begin{equation}\label{MProx}
		    \text{MProx}(T,0):=\left\{ x \in X: \liminf_{n\to \infty} \frac{1}{n}\sum_{k=1}^{n}d(T^kx,0)=0 \right\}.
		\end{equation}
    \end{itemize}
    \begin{remark}\label{remark1}
 In \cite{JiangLi}, it was observed that $\text{MProx}(T,0)$ is a $G_\delta$ subset of $X$. Furthermore, by \cite[Lemma~4.7]{JiangLi}, if $\text{MAsym}(T,0)$ is residual in $X$, then $\text{MAsym}(T,0)=X$.
	\end{remark}
We now state as lemmas a series of results from \cite{JiangLi} that will be used in the this section.
    \begin{lemma}\cite[Proposition 4.11]{JiangLi}\label{LMLY1}
		Let $X$ be a Fr\'echet space with compatible metric $d$ and $T:X\to X$ be a continuous linear operator. Then the following assertions are equivalent:
		\begin{enumerate}[label=$(\alph*)$]
			\item \label{LMLY1 a} $T$ is mean sensitive;
			
			\item \label{LMLY1 b} there exist a sequence $(y_k)_{k\in \N}$ in $X$ and an increasing sequence $(N_k)_{k\in \N}$ in $\mathbb{N}$ such that $\lim_{k\to\infty} y_k=0$ and
			\[
			\inf_{k\in\mathbb{N}} \frac{1}{N_k}\sum_{i=1}^{N_k} d(T^i y_k,0)>0.
			\]
		\end{enumerate}
	\end{lemma}
	
	\begin{lemma}\cite[Theorem 4.15]{JiangLi}\label{LMLY2}
		Let $T:X \to X$ be a continuous linear operator on a Fréchet space $X$. Then the following statements are equivalent:
		\begin{enumerate}[label=$(\alph*)$]
			\item\label{LMLY2 a} $T$ is mean Li--Yorke chaotic; 
			\item\label{LMLY2 b} $T$ is mean Li--Yorke sensitive;
			\item\label{LMLY2 c} $T$ admits an absolutely mean semi-irregular vector.
		\end{enumerate}
	\end{lemma}
	\begin{lemma}\cite[Theorem 4.13 and Theorem 4.27]{JiangLi}\label{LMLY3}
		Let $T:X \to X$ be a continuous linear operator on a Fréchet space $X$ with compatible metric $d$. Then the following statements are equivalent:
		\begin{enumerate}[label=$(\alph*)$]
			\item\label{LMLY3 a} $T$ admits a dense set of absolutely mean semi-irregular vectors;
            \item\label{LMLY3 a'} $T$ admits a residual set of absolutely mean semi-irregular vectors;
			\item\label{LMLY3 b} $\text{MProx}(T,0)$ is dense in $X$ and there exists $x \in X$ such that
			\[
			\limsup_{n \to \infty} \frac{1}{n} \sum_{i=1}^{n} d(T^{i}x,0) > 0.
			\]
		\end{enumerate}
	\end{lemma}
    In the next theorem, we establish a characterization of mean Li–Yorke chaos for weighted backward shifts on admissible Fréchet sequence spaces over directed trees.
\begin{theorem} \label{thm:mean-LY-tree}
Let $V$ be a directed tree with generations enumerated according to some fixed vertex $v_0$ and let $X$ be an admissible Fréchet sequence space over $V$, with compatible metric $d$. Suppose that the weighted backward shift $B_\lambda$, with nonzero weights $\lambda := (\lambda_v)_{v \in V}$, is well-defined and continuous on $X$. Then $B_\lambda$ is mean Li--Yorke chaotic if and only if there exist a nonempty set $S\subset \{n\in\mathbb Z:\Gen_n\neq\varnothing\}$ and, for each $n\in S$, a vector $f_n \in X$ supported in $\Gen_n$ such that, putting $E:=\Span\{B_\lambda^if_n : n\in S, i\in\mathbb N_0\},$ the following conditions hold:
\begin{enumerate}[label=$(\Alph*)$]
\item \label{teo:mean:cond:A} there is a strictly increasing sequence $(M_k)_{k\in \N}$ of positive integers such that
\[
\lim_{k\to\infty} \frac{1}{M_k}\sum_{i=1}^{M_k} d(B_\lambda^if_n,0)=0 \quad \text{for all } \, n\in S;
\]
\item \label{teo:mean:cond:B} there exist $\varepsilon>0$, a strictly increasing sequence $(N_k)_{k\in\N}$ of positive integers, and a sequence $(h_k)_{k\in\N}\subset E$ such that
\[
d(h_k,0)<\frac{1}{k} \quad \text{and} \quad \frac{1}{N_k}\sum_{i=1}^{N_k} d(B_\lambda^i h_k,0)>\varepsilon \quad \text{for all } \, k\in\mathbb N.
\]
\end{enumerate}
\end{theorem}
\begin{proof}
($\Rightarrow$): Suppose that $B_\lambda$ is mean Li--Yorke chaotic. By Lemma \ref{LMLY2}, it admits an absolutely mean semi-irregular vector $f\in X$. Then we can take a subsequence $(M_k)_{k\in \N}\subset \N$ such that
\begin{equation}\label{cajou}
\lim_{k\to\infty}\frac{1}{M_k}\sum_{j=1}^{M_k}d(B_\lambda^jf,0)=0.
\end{equation}
Put $S:=\{n\in\mathbb Z:\Gen_n\neq \emptyset\}$. Fix $\eta>0$. By Lemma \ref{LM1}, the family of operators $(\Phi_{n})_{n\in \Z}$ given by 
        $$\Phi_{n}: g\in X \to g\chi_{\operatorname{Gen}_n}\in X, \quad \text{for } n\in \Z,$$
        is equicontinuous. Therefore, there exists $\delta \in (0,1)$ such that, for all $n\in \Z$ and $f\in X$,
\begin{equation}\label{eq:pitomba}
d(f,0)<\delta \implies d(f\chi_{\Gen_n},0)<\frac{\eta}{2}.
\end{equation}
By (\ref{cajou}), there exists $k_0\in \N$ such that
\begin{equation}\label{eq:pitomba2}
    \frac{1}{M_k}\sum_{j=1}^{M_k}d(B_\lambda^jf,0)< \frac{\eta\delta}{2}\quad \text{for all }k>k_0.
\end{equation}
Fix $k>k_0$ and consider the sets
$$A:=\left\{ 1\le j\le M_k: d(B_\lambda^jf,0)<\delta \right\}\quad \text{and}\quad B:=\left\{ 1\le j\le M_k: d(B_\lambda^jf,0)\geq\delta \right\}.$$
By (\ref{eq:pitomba}), (\ref{eq:pitomba2}), and Lemma \ref{LM2} we have, for all $n\in S$,
\begin{align*}
\frac{1}{M_k}\sum_{j=1}^{M_k}d(B_\lambda^j(f\chi_{\Gen_n}),0) 
    &= \frac{1}{M_k} \sum_{j\in A}^{} d(\chi_{\Gen_{n-j}}B_\lambda^jf,0)+ \frac{1}{M_k}\sum_{j\in B}^{} d(\chi_{\Gen_{n-j}}B_\lambda^jf,0)\\
    &\leq \frac{1}{M_k}M_k\frac{\eta}{2} + \frac{1}{M_k}\card (B)\\
    &\leq \frac{\eta}{2}+ \frac{1}{\delta M_k}\sum_{j=1}^{M_k}d(B_\lambda^jf,0)\\
& <\frac{\eta}{2}+\frac{\eta\delta}{2\delta}=\eta.
\end{align*}
Hence, defining $f_n:=f\chi_{\Gen_n}$ for all $n\in S$, condition \ref{teo:mean:cond:A} holds.

Define $E:= \Span\{ B_\lambda^if_n: n\in S,\ i\in\mathbb N_0\}$ and $Y:=\overline E$. The space $Y$ is a closed $B_\lambda$-invariant subspace of $X$. Notice that $f\in Y$. Consequently, $f$ is an absolutely mean semi-irregular vector for the restriction $B_\lambda|_Y$. By Lemma \ref{LMLY2}, $B_\lambda|_Y$ is mean Li--Yorke sensitive and, in particular, mean sensitive. By Lemma \ref{LMLY1}, there exist a sequence $(g_k)_{k\in\N}\subset Y$ such that $\lim_{k\to\infty}g_k=0$ and a strictly increasing sequence $(N_k)_{k\in \N}$ of positive integers such that
\begin{equation}\label{eq:ananas}
c:=
\inf_{k\in\mathbb N}
\frac{1}{N_k}\sum_{j=1}^{N_k}d(B_\lambda^jg_k,0)>0.
\end{equation}
By passing to a subsequence if necessary, we may assume that
\begin{equation}\label{eq:abricot}
d(g_k,0)<\frac{1}{2k}.
\end{equation}
By the density of $E$ in $Y$ and the continuity of $B_{\lambda}$, for each $k \in \N$, we can take $h_k\in E$ such that 
\begin{equation}\label{eq:peche}
d\left( B_{\lambda}^jg_k,B_{\lambda}^jh_k \right)< \min\left\{ \frac{1}{2k},\frac{c}{2} \right\}, \quad\text{for all } 0\leq j\leq N_k.
\end{equation}
Thus, equations \eqref{eq:ananas}-\eqref{eq:peche} imply that, for each $k\in \N$, $d(h_k,0)<\frac{1}{k}$ and
\begin{align*}
\frac{1}{N_k}\sum_{j=1}^{N_k}d(B_\lambda^jh_k,0) &\geq \frac{1}{N_k}\sum_{j=1}^{N_k}d(B_\lambda^jg_k,0) - \frac{1}{N_k}\sum_{j=1}^{N_k}d\bigl(B_\lambda^jh_k,B_\lambda^jg_k\bigr)  > c-\frac{1}{N_k} N_k \frac{c}{2}=\frac{c}{2}.
\end{align*}
Therefore, condition \ref{teo:mean:cond:B} holds with $\varepsilon=c/2$.

 ($\Leftarrow$): Assume that conditions \ref{teo:mean:cond:A} and \ref{teo:mean:cond:B} hold, and define again $Y:=\overline E$. Once more, $Y$ is a closed $B_\lambda$-invariant subspace. Let $f\in E$. Then there exist $r\in \N$, scalars $\alpha_1,\cdots, \alpha_r \in \K$, and $(m_1,n_1),\cdots, (m_r,n_r)\in S\times \N_0$ such that $f=\sum_{j=1}^{r}\alpha_j B_{\lambda}^{n_j}f_{m_j}$. Thus, by condition \ref{teo:mean:cond:A} and properties \ref{metric:P1} and \ref{metric:P2} of the metric $d$, we obtain 
 $$\frac{1}{M_k}\sum_{j=1}^{M_k}d(B_\lambda^{j}f,0)\le \sum_{i=1}^{r}\frac{(1+|\alpha_i|)}{M_k}\sum_{j=n_i+1}^{M_k+n_i}d(B_\lambda^{j}f_{m_i},0)\xrightarrow{k\to\infty}0.$$
 Therefore, $E \subset \MProx(B_\lambda|_Y,0)$ and, hence, $\MProx(B_\lambda|_Y,0)$ is dense in $Y$. On the other hand, condition \ref{teo:mean:cond:B} and Lemma \ref{LMLY1} show that $B_\lambda|_Y$ is mean sensitive. Applying the mean sensitivity at zero, there exists $z\in Y$ such that
\begin{equation}\label{eql12}
    \limsup_{N\to\infty} \frac{1}{N}\sum_{j=1}^{N}d(B_\lambda^jz,0)>0.
\end{equation}
Since $\MProx(B_\lambda|_Y,0)$ is dense in $Y$, combining \eqref{eql12} with Lemma \ref{LMLY3}, we conclude that $B_\lambda|_Y$ admits a dense set of absolutely mean semi-irregular vectors. By Lemma \ref{LMLY2}, $B_\lambda|_Y$ is mean Li--Yorke chaotic and, hence, $B_\lambda$ is mean Li--Yorke chaotic. This completes the proof.
\end{proof}
The following corollary improves \cite[Theorem 30]{joaojmaa} by removing an extra condition that the author assumes on the seminorms that define the Fréchet topology of the space.
\begin{corollary}
Let $X$ be a Fréchet sequence space over $\mathbb Z$ with compatible metric $d$ and suppose that $(e_n)_{n\in\mathbb Z}$ is a Schauder basis of $X$. Let $B_w$ be a well-defined continuous weighted backward shift on $X$ with nonzero weights $w=(w_n)_{n\in\mathbb Z}$. Then $B_w$ is mean Li--Yorke chaotic if and only if the following conditions hold:
\begin{enumerate}[label=$(\Alph*)$]
\item \label{corol:joao30:A}
$\liminf_{N\to\infty}\frac{1}{N} \sum_{i=1}^{N}d(w_{-i}\cdots w_{-1}e_{-i},0)=0;$
\item \label{corol:joao30:B} there exist $\eps>0$ and an increasing sequence $(N_k)_{k\in \N}$ of positive integers such that, for every $k\in\N$, there are $r=r(k)\in\N$ and $b_{-r,k},\ldots,b_{r,k}\in\mathbb K$ satisfying $d\big(\sum_{j=-r}^{r}b_{j,k}e_j,0\big)<\frac{1}{k}$ and
\[
\frac{1}{N_k}\sum_{i=1}^{N_k} d\left( \sum_{j=-r}^{r} b_{j,k} w_{j-i}\cdots w_{j-1}e_{j-i}, 0 \right) \geq \eps.
\]
\end{enumerate}
\end{corollary}
\begin{proof}
Recall that, by Example \ref{ExFNZ}, we may regard $B_w$ as the weighted backward shift $B_\lambda$ introduced in \eqref{eqq0} on the admissible Fréchet sequence space over $V=\Z$. For this, it suffices to define $\lambda_{n+1}:=w_n$ for every $n\in\Z$. As usual, we consider $\Gen_n:=\{n\}$ for $n\in \Z$. We shall use the elementary fact that, for every $m\in\Z$ and every increasing sequence $(M_k)_{k\in\N}$ of positive integers,
\begin{equation} \label{eq:mean-li-deslocado}
\lim_{k\to\infty}\frac{1}{M_k}\sum_{i=1}^{M_k}d(B_w^ie_m,0)=0 \quad \Longleftrightarrow \quad \lim_{k\to\infty} \frac{1}{M_k}\sum_{i=1}^{M_k} d(B_w^ie_0,0)=0.
\end{equation}
Indeed, the two orbit sequences differ by a fixed translation of the indices and multiplication by a nonzero scalar.

($\Rightarrow$): Suppose that $B_w$ is mean Li--Yorke chaotic. By Theorem \ref{thm:mean-LY-tree}, there exist $S\subset\Z$, vectors $f_n$ supported in $\Gen_n$ for $n\in S$ and $E=\Span\{B_w^if_n:n\in S,\ i\in\N_0\}$ satisfying conditions \ref{teo:mean:cond:A} and \ref{teo:mean:cond:B} of Theorem \ref{thm:mean-LY-tree}. Fix $n_0\in S$ such that $f_{n_0}\ne 0$. Since $\Gen_{n_0}=\{n_0\}$, we have $f_{n_0}=a_{n_0}e_{n_0}$ for some $a_{n_0}\neq0$. By condition \ref{teo:mean:cond:A} of Theorem \ref{thm:mean-LY-tree}, there exists a subsequence $(M_k)_{k\in\N}\subset \N$ with
\[
\lim_{k\to\infty}\frac{1}{M_k}\sum_{i=1}^{M_k}d(B_w^ie_{n_0},0)=0.
\]
By \eqref{eq:mean-li-deslocado}, the same holds for $e_0$. Therefore, condition \ref{corol:joao30:A} follows. Finally, every generator $B_w^if_n$ of $E$ is a scalar multiple of a canonical vector. Hence every $h_k\in E$ given by condition \ref{teo:mean:cond:B} of Theorem \ref{thm:mean-LY-tree} is finitely supported. We may therefore write
\[
h_k=\sum_{j=-r(k)}^{r(k)}b_{j,k}e_j
\]
for some $r(k)\in\N$ and some scalars $b_{j,k}$. Since $B_w^ie_j=w_{j-i}\cdots w_{j-1}e_{j-i}$, condition \ref{teo:mean:cond:B} of Theorem \ref{thm:mean-LY-tree} corresponds exactly to condition \ref{corol:joao30:B} in the statement.

($\Leftarrow$): By \ref{corol:joao30:A}, there is an increasing sequence $(M_k)$ such that
\[
\lim_{k\to\infty} \frac{1}{M_k} \sum_{i=1}^{M_k} d(B_w^ie_0,0)=0.
\]
By \eqref{eq:mean-li-deslocado}, the same conclusion holds with $e_n$ in place of $e_0$ for every $n\in\Z$. In Theorem \ref{thm:mean-LY-tree}, take $S=\Z$ and $f_n=e_n$ for every $n\in\Z$. Then condition \ref{teo:mean:cond:A} of Theorem \ref{thm:mean-LY-tree} holds and
\[
E=\Span\{B_w^ie_n : n\in \Z, \, i\in \N_0\} = c_{00}(\Z).
\]
For each $k\in\N$, let $h_k=\sum_{j=-r(k)}^{r(k)}b_{j,k}e_j\in E$. Thus, condition \ref{corol:joao30:B} in the statement gives condition \ref{teo:mean:cond:B} of Theorem \ref{thm:mean-LY-tree} with $\eps/2$. Therefore, $B_w$ is mean Li--Yorke chaotic.
\end{proof}
Let $V$ be a directed tree, $X$ be an admissible Fréchet sequence space over $V$ and $\lambda=(\lambda_v)_{v\in V}$ be a weight sequence such that the weighted backward shift $B_\lambda$ is well defined and continuous on $X$. We say that $\lambda$ satisfies condition \eqref{cond:C2} when there is $v_0\in V$ satisfying
\begin{equation}\label{cond:C2} \tag{$C2$}
    \liminf_{N\to\infty}\frac{1}{N} \sum_{j=1}^{N}d\left(\lambda\bigl(\operatorname{par}^j(v_0)\to v_0\bigr)e_{\operatorname{par}^j(v_0)},0\right)=0.
\end{equation}
\begin{proposition} \label{prop:mean-LY-dicotomia}
    Let $V$ be a directed tree and $X$ an admissible Fréchet sequence space over $V$ with compatible metric $d$. Suppose that the weighted backward shift $B_\lambda$, with nonzero weights $\lambda := (\lambda_v)_{v \in V}$, is well-defined and continuous on $X$. Suppose that condition \eqref{cond:C2} holds. Then $\MProx(B_\lambda,0)$ is dense in $X$ and exactly one of the following holds:
    \begin{enumerate}[label=$(\alph*)$]
        \item \label{prop:mean-LY-dicotomia:a} $\MAsym(B_\lambda,0)=X$;
        \item \label{prop:mean-LY-dicotomia:b} the absolutely mean semi-irregular vectors form a residual subset of $X$.
    \end{enumerate}
\end{proposition}
\begin{proof}
We shall prove that $c_{00}(V)\subset\MProx(B_\lambda,0)$. If $V$ is rooted, the inclusion is trivial. Suppose now that $V$ is unrooted. By condition \eqref{cond:C2}, there are $v_0\in V$ and an increasing sequence $(M_k)$ of positive integers such that
\begin{equation}\label{propC2}
    \lim_{k\to\infty}\frac{1}{M_k}\sum_{j=1}^{M_k}d(B_\lambda^je_{v_0},0)=0.
\end{equation}
Let $f\in c_{00}(V)$, say $f=\sum_{v\in F}f(v)e_v$ with $F\subset V$ finite. For each $v\in F$, let $p_v,q_v\in\N_0$ such that $\Par^{p_v}(v_0)=\Par^{q_v}(v)$. Note that, for all $j\ge q_v+1$, we can write
\begin{align*}
B_\lambda^je_v&=B_\lambda^{j-q_v}B_\lambda^{q_v}e_v\\
&=\lambda(\operatorname{par}^{q_v}(v)\to v)B_\lambda^{j-q_v}e_{\operatorname{par}^{q_v}(v)}\\
&=\overset{=:c_v}{\overbrace{\frac{\lambda(\operatorname{par}^{q_v}(v)\to v)}{\lambda(\operatorname{par}^{p_v}(v_0)\to v_0)}}}B_\lambda^{j-q_v}B_\lambda^{p_v}e_{v_0}\\
&=c_v B_{\lambda}^{j-q_v+p_v}e_{v_0}.
\end{align*}
Fix $N_0\in\N$ such that $N_0\geq  p_v-q_v$ for all $v\in F$. Thus, for all $k$ large enough and for all $v\in F$, we have
\begin{align*}
\frac{1}{M_k} \sum_{j=1}^{M_k} d(B_\lambda^je_v,0)&=\frac{1}{M_k}\sum_{j=1}^{q_v} d(B_\lambda^je_v,0)+
\frac{1}{M_k} \sum_{j=q_v+1}^{M_k} d(B_\lambda^je_v,0)\\
&\le \frac{q_v}{M_k} +\frac{1}{M_k} \sum_{j=q_v+1}^{M_k} d(c_vB_\lambda^{j-q_v+p_v}e_{v_0},0)\\
&\le \frac{q_v}{M_k}+\frac{(1+|c_v|)}{M_k}\sum_{i=1}^{M_{k}+N_0}d(B_\lambda^{i}e_{v_0},0).
\end{align*}
Using \eqref{propC2} and properties \ref{metric:P1} and \ref{metric:P2} of the metric $d$, it follows that the right-hand side of the last inequality tends to zero as $k\to\infty$. We thus conclude that $f\in \MProx(B_\lambda,0).$ This shows that $\MProx(B_\lambda,0)$ is dense in $X$. It remains to prove the final dichotomy. Of course \ref{prop:mean-LY-dicotomia:a} and \ref{prop:mean-LY-dicotomia:b} are mutually exclusive. Suppose that $\MAsym(B_\lambda,0)\neq X$. Then there is $f\in X$ such that 
\[
\limsup_{N\to\infty} \frac{1}{N} \sum_{j=1}^{N} d(B_\lambda^jf,0) > 0.
\]
Since $\MProx(B_\lambda,0)$ is dense in $X$, it follows from Lemma \ref{LMLY3} that the absolutely mean semi-irregular vectors form a residual subset of $X$. 
\end{proof}
For the next proposition, we will need the following lemma. To this end, we recall the following classical topological result.
\begin{itemize}
    \item (Mycielski theorem) If $X$ is a non-trivial separable Fréchet space and $R$ is residual in $X\times X$, then there exists a dense {\it Mycielski set} $\Gamma\subset X$, that is, a countable union of Cantor sets, such that $\Gamma\times\Gamma\setminus\Delta\subset R$, where $\Delta:=\{(x,x):x\in X\}$.
\end{itemize}
\begin{lemma}\label{LemmaPrep9}
       Let $T:X\to X$ be a continuous linear operator on a non-trivial separable Fréchet space $X$ with compatible metric $d$. Then the following assertions are equivalent:
\begin{enumerate}[label=$(\roman*)$]
    \item \label{LemmaPrep9 - i} $T$ is densely mean Li--Yorke chaotic;
    \item \label{LemmaPrep9 - ii}$\MProx(T,0)$ is dense in $X$ and there exists $x\in X$ such that
    \begin{equation}\label{eqp80}
        \limsup_{N\to\infty} \frac{1}{N} \sum_{j=1}^{N} d(T^jx,0) > 0.
    \end{equation}
\end{enumerate}
\end{lemma}
\begin{proof}
\ref{LemmaPrep9 - i}$\Rightarrow$\ref{LemmaPrep9 - ii}: Let $\Gamma\subset X$ be an uncountable dense set such that every pair of distinct elements of $\Gamma$ is a mean Li--Yorke pair. Fix $x_0\in \Gamma$ and consider the set $D:=\{x-x_0: x\in\Gamma\setminus\{x_0\}\}$. Note that $D$ is dense. By the invariance of the metric $d$, it is easy to see that $D\subset \MProx(T,0)$ and each $z\in D$ satisfies \eqref{eqp80}.

\ref{LemmaPrep9 - ii}$\Rightarrow$\ref{LemmaPrep9 - i}: Let
$$A:=\left\{ x\in X: \limsup_{N\to\infty} \frac{1}{N} \sum_{j=1}^{N} d(T^jx,0) > 0 \right\}.$$
By Lemma \ref{LMLY3}, $A$ is residual in $X$. From Remark \ref{remark1}, $\MProx(T,0)$ is also residual. Therefore, $R:= \MProx(T,0)\cap A$ is residual. Since the map $\Phi:X\times X\to X$, $\Phi(x,y)=x-y$, is continuous, open and surjective, the set
    $$\Phi^{-1}(R)=\{(x,y)\in X\times X:x-y\in R\}$$ is residual in $X\times X$. By Mycielski's theorem, there exists a dense Mycielski set $\Gamma\subset X$ such that $x-y\in R$ for every distinct $x,y\in\Gamma$. Therefore, $\Gamma$ is a dense set such that each pair of distinct points in $\Gamma$ is a mean Li--Yorke pair.
\end{proof}
The next proposition establishes a characterization of Li–Yorke chaos for weighted backward
shifts on admissible Fréchet sequence spaces over directed trees, under the assumption of condition \eqref{cond:C2}.
\begin{proposition} \label{prop:mean-LY-equiv}
    Let $V$ be a directed tree and $X$ an admissible Fréchet sequence space over $V$ with compatible metric $d$. Suppose that the weighted backward shift $B_\lambda$, with nonzero weights $\lambda := (\lambda_v)_{v \in V}$, is well-defined and continuous on $X$. Suppose that condition \eqref{cond:C2} holds. Then the following assertions are equivalent:
    \begin{enumerate}[label=$(\roman*)$]
        \item \label{prop:mean-LY-equi (i)} $B_\lambda$ is mean Li--Yorke chaotic;
        \item \label{prop:mean-LY-equi (ii)} $B_\lambda$ is densely mean Li--Yorke chaotic;
        \item \label{prop:mean-LY-equi (iii)} $B_\lambda$ is mean sensitive;
        \item \label{prop:mean-LY-equi (iv)} there exists $f\in X$ such that $\limsup_{N\to\infty}\frac{1}{N}\sum_{j=1}^Nd(B_\lambda^jf,0)>0$;
        \item \label{prop:mean-LY-equi (v)} there are $\varepsilon>0$, an increasing sequence $(N_k)\subset \N$ and, for each $k\in\N$, a finite set $F_k\subset V$, and $h_k\in X$ with support in $F_k$, such that
        \[d(h_k,0)<\frac{1}{k}\quad \text{and}\quad \frac{1}{N_k}\sum_{j=1}^{N_k}d(B_\lambda^jh_k,0)>\varepsilon \quad \text{for all } \, k\in\N.\]
    \end{enumerate}
\end{proposition}
\begin{proof}
    By proposition \ref{prop:mean-LY-dicotomia} and its proof, $c_{00}(V)\subset\MProx(B_\lambda,0)$ and therefore is dense in $X$. 

    \ref{prop:mean-LY-equi (i)}$\Rightarrow$\ref{prop:mean-LY-equi (ii)}: Since $B_{\lambda}$ is mean Li--Yorke chaotic, by Lemma \ref{LMLY2}, we conclude that $B_{\lambda}$ admits an absolutely mean semi-irregular vector. Since $\MProx(B_\lambda,0)$ is dense, by Lemma \ref{LemmaPrep9}, we conclude that $B_\lambda$ is densely mean Li--Yorke chaotic.

    \ref{prop:mean-LY-equi (ii)}$\Rightarrow$\ref{prop:mean-LY-equi (iii)}: Since $B_{\lambda}$ is densely mean Li--Yorke chaotic, by Lemma \ref{LMLY2}, $B_{\lambda}$ is mean Li--Yorke sensitive and, in particular, mean sensitive.  

    \ref{prop:mean-LY-equi (iii)}$\Rightarrow$\ref{prop:mean-LY-equi (iv)}: It suffices to apply the definition of mean sensitivity at the origin.

    \ref{prop:mean-LY-equi (iv)}$\Rightarrow$\ref{prop:mean-LY-equi (v)}: Since $\MProx(B_\lambda,0)$ is dense, by item \ref{prop:mean-LY-equi (iv)} and Lemma \ref{LemmaPrep9}, $B_{\lambda}$ is densely mean Li--Yorke chaotic. In particular, item \ref{teo:mean:cond:B} of Theorem \ref{thm:mean-LY-tree} holds. Therefore, by condition \ref{frechet:condition- B1} of Definition \ref{defAd} and the continuity of $B_{\lambda}$ and of the metric $d$, condition  \ref{teo:mean:cond:B} of Theorem \ref{thm:mean-LY-tree} implies item \ref{prop:mean-LY-equi (v)}.

    \ref{prop:mean-LY-equi (v)}$\Rightarrow$\ref{prop:mean-LY-equi (i)}: Item \ref{prop:mean-LY-equi (v)} implies item \ref{LMLY1 b} of Lemma \ref{LMLY1}. Thus, $B_{\lambda}$ is mean sensitive and, in particular, there exists $f\in X$ such that 
    $$\limsup_{N\to\infty}\frac{1}{N}\sum_{j=1}^Nd(B_\lambda^jf,0)>0.$$
    Since $\MProx(B_\lambda,0)$ is dense, by Lemma \ref{LemmaPrep9}, we conclude that $B_{\lambda}$ is mean Li--Yorke chaotic.
\end{proof}
If $V$ is rooted, then condition \eqref{cond:C2} is automatically satisfied. Hence, we obtain the following corollary.
\begin{corollary}
    Let $V$ be a rooted directed tree. Let $X$ be an admissible Fréchet sequence space over $V$. Suppose that the weighted backward shift $B_\lambda$, with nonzero weights $\lambda := (\lambda_v)_{v \in V}$, is well-defined and continuous on $X$. Then the assertions \ref{prop:mean-LY-equi (i)}-\ref{prop:mean-LY-equi (v)} from Proposition \ref{prop:mean-LY-equiv} are equivalent.
\end{corollary}

For the next corollary, we will need the following definition and a result from \cite{JiangLi}.

\begin{definition}
Let $(X,d)$ be a metric space and let $f:X\to X$ be a continuous map. We say that $f$ is \emph{mean equicontinuous} if, for every $\eps > 0$, there exists $\delta>0$ such that 
\[
x,y\in X, \, d(x,y) < \delta \implies \limsup_{N\to\infty} \frac{1}{N}\sum_{j=1}^N d(f^j(x),f^j(y))<\varepsilon.
\]
\end{definition}

\begin{lemma}[{\cite[Theorem 4.9]{JiangLi}}] \label{lemma:mean-eq:mean-sens}
Let $X$ be a Fréchet space and $T:X\to X$ be a continuous linear operator. Then exactly one of the following holds:
\begin{enumerate}[label=$(\roman*)$]
\item \label{lemma:mean-eq:mean-sens - i} $T$ is mean equicontinuous;
\item \label{lemma:mean-eq:mean-sen - ii} $T$ is mean sensitive.
\end{enumerate}
\end{lemma}

\begin{corollary} \label{cor:mean-dicotomy}
Let $V$ be a directed tree and $X$ an admissible Fréchet sequence space over $V$. Suppose that the weighted backward shift $B_\lambda$, with nonzero weights $\lambda := (\lambda_v)_{v \in V}$, is well-defined and continuous on $X$. Suppose that condition \eqref{cond:C2} holds. Then exactly one of the following assertions holds:
\begin{enumerate}[label=$(\alph*)$]
\item \label{mean-dicotomy a)}$B_\lambda$ is mean equicontinuous;
\item \label{mean-dicotomy b)}$B_\lambda$ is densely mean Li--Yorke chaotic.
\end{enumerate}
In particular, $B_\lambda$ is mean Li--Yorke chaotic if and only if it is not mean equicontinuous.
\end{corollary}
\begin{proof}
  Since every mean Li--Yorke chaotic operator is mean sensitive, it follows from Lemma \ref{lemma:mean-eq:mean-sens} that assertions \ref{mean-dicotomy a)} and \ref{mean-dicotomy b)} are mutually exclusive. Moreover, combining Lemma \ref{lemma:mean-eq:mean-sens} and Proposition \ref{prop:mean-LY-equiv}, we conclude that one of these two assertions must hold.
\end{proof}
In our next result, we establish necessary conditions in terms of the seminorms of an admissible Fréchet sequence space for the weighted backward shift $B_{\lambda}$ to be mean Li--Yorke chaotic. To this end, we need the following definition and lemma.
\begin{definition}
		Let $X$ be a Fréchet space endowed with an increasing sequence
		$(\|\cdot\|_k)_{k\in\mathbb{N}}$ of seminorms inducing the compatible metric $d$. Let $T:X\to X$ be a continuous linear operator. A vector $x\in X$ is called \textit{absolutely mean} $m$-\textit{irregular}, for some $m\in\N$, if 
		$$\liminf_{n\to \infty}\frac{1}{n}\sum_{k=1}^{n}d(T^kx,0)=0\quad \text{and}\quad \limsup_{n\to \infty}\frac{1}{n}\sum_{k=1}^{n}\left\| T^k x \right\|_m=\infty.$$ 
	\end{definition}
	\begin{lemma}\cite[Proposition 4.29]{JiangLi}\label{L0}
		Let $T:X \to X$ be a continuous linear operator on a Fréchet space $X$. Then the set of absolutely 
		mean semi-irregular vectors is contained in the closure of the set of 
		absolutely mean $m$-irregular vectors for some $m \in \mathbb{N}$.
	\end{lemma}

\begin{proposition} \label{prop:mean-LY-tree}
Let $V$ be a directed tree with generations enumerated according to some fixed vertex $v_0$ and $\left( X, \left( \left\| \cdot  \right\|_k \right)_{k\in\N} \right)$ be an admissible Fréchet sequence space over $V$, with compatible metric $d$. Suppose that the weighted backward shift $B_\lambda$, with nonzero weights $\lambda := (\lambda_v)_{v \in V}$, is well-defined and continuous on $X$. If $B_\lambda$ is mean Li--Yorke chaotic, then there exist $S \subset \{n\in\mathbb Z:\Gen_n\neq\varnothing\}$ and, for each $n\in S$, a vector $f_n \in X$ supported in $\Gen_n$ such that, putting $E:=\Span\{B_\lambda^if_n : n\in S, i\in\mathbb N_0\},$ the following conditions hold:
\begin{enumerate}[label=$(\Alph*)$]
\item \label{teo:mean:cond:A1} there is a strictly increasing sequence $(M_k)_{k\in \N}$ of positive integers such that
\[
\lim_{k\to\infty} \frac{1}{M_k}\sum_{i=1}^{M_k} d(B_\lambda^if_n,0)=0 \quad \text{for all } \, n\in S;
\]
\item \label{teo:mean:cond:B1} there exist $m \in \N$, a strictly increasing sequence $(N_k)_{k\in\N}$ of positive integers and a sequence $(h_k)_{k\in\N}\subset E$ such that
\[
\|h_k\|_{p(k)}>0 \quad \text{and} \quad \frac{1}{N_k\|h_k\|_{p(k)}}\sum_{j=1}^{N_k} \|B_\lambda^j h_k\|_{m} \geq  k,
\]
where $p(k) = \max\{k, m \}$ for all $k \in \N$.
\end{enumerate}
\end{proposition}
\begin{proof}
Since $B_\lambda$ is mean Li--Yorke chaotic, Theorem \ref{thm:mean-LY-tree} provides a nonempty set $S\subset\{n\in\Z:\Gen_n\neq\varnothing\}$ and vectors $f_n\in X$, supported in $\Gen_n$, for which condition \ref{teo:mean:cond:A1} holds. Define $Y:=\overline E$. The proof of Theorem \ref{thm:mean-LY-tree} also shows that $Y$ is $B_\lambda$-invariant and that $B_\lambda|_Y$ is mean Li--Yorke chaotic. By Lemma \ref{LMLY2}, $B_\lambda|_Y$ admits an absolutely mean semi-irregular vector. Hence, by Lemma \ref{L0}, there exist $m\in\N$ and an absolutely mean $m$-irregular vector $f\in Y$. Increasing $m$ if necessary, we may assume that $\|f\|_m>0$. Therefore,
\[
\limsup_{N\to\infty}\frac{1}{N} \sum_{j=1}^N\|B_\lambda^jf\|_m=\infty.
\]
For each $k\in\N$, put $p(k)=\max\{k,m\}$. Choose a strictly increasing sequence $(N_k)_{k\in\N}$ of positive integers such that
\begin{equation}\label{prop:mean:1}
\frac{1}{N_k}\sum_{j=1}^{N_k}\|B_\lambda^jf\|_m>4k\|f\|_{p(k)} \quad\text{for every }k\in\N.
\end{equation}
Since $E$ is dense in $Y$ and the seminorms and $B_\lambda$ are continuous, for each $k\in\N$ we can choose $h_k\in E$ such that
\begin{equation}\label{prop:mean:2}
\|h_k-f\|_{p(k)}<\|f\|_{p(k)} \quad \text{and} \quad \frac{1}{N_k}\sum_{j=1}^{N_k}\|B_\lambda^j(h_k-f)\|_m<k\|f\|_{p(k)}.
\end{equation}
In particular, $\|h_k\|_{p(k)}\geq\|f\|_{p(k)}-\|h_k-f\|_{p(k)}>0$ and $\|h_k\|_{p(k)}<2\|f\|_{p(k)}$. Moreover, by \eqref{prop:mean:1} and \eqref{prop:mean:2},
\begin{equation*}
\frac{1}{N_k}\sum_{j=1}^{N_k}\|B_\lambda^jh_k\|_m
    \geq\frac{1}{N_k}\sum_{j=1}^{N_k}\|B_\lambda^jf\|_m-\frac{1}{N_k}\sum_{j=1}^{N_k}\|B_\lambda^j(h_k-f)\|_m >3k\|f\|_{p(k)} >k\|h_k\|_{p(k)}.
\end{equation*}
Consequently,
\[
\frac{1}{N_k\|h_k\|_{p(k)}}\sum_{j=1}^{N_k}\|B_\lambda^jh_k\|_m>k
\]
for every $k\in\N$, and condition \ref{teo:mean:cond:B1} follows.
\end{proof}

\begin{remark} \label{mean:Kothe}
The corresponding necessary condition for the Köthe sequence spaces $\lambda_p(A,V)$ is obtained directly from Proposition~\ref{prop:mean-LY-tree} by inserting
\[ 
\|B_\lambda^jg\|_m = 
\begin{cases} \displaystyle
    \left(\sum_{v\in V}a_{v,m}^p \left|\sum_{u\in\Chi^j(v)} \lambda(v\to u)g(u) \right|^p \right)^{1/p}, & \text{when }1\leq p<\infty \\ \displaystyle
    \sup_{v\in V}a_{v,m} \left| \sum_{u\in\Chi^j(v)} \lambda(v\to u)g(u) \right|, & \text{when }p=0.
\end{cases}
\]
\end{remark}

\begin{example} \label{ex:fraise}
Let $\alpha\in \K\setminus\{0\}$ and consider the directed tree $V$ given by Figure \ref{fig:example:mean-LY:alpha}.
\begin{figure}[H]
    \centering
    \includegraphics[height=0.19\linewidth]{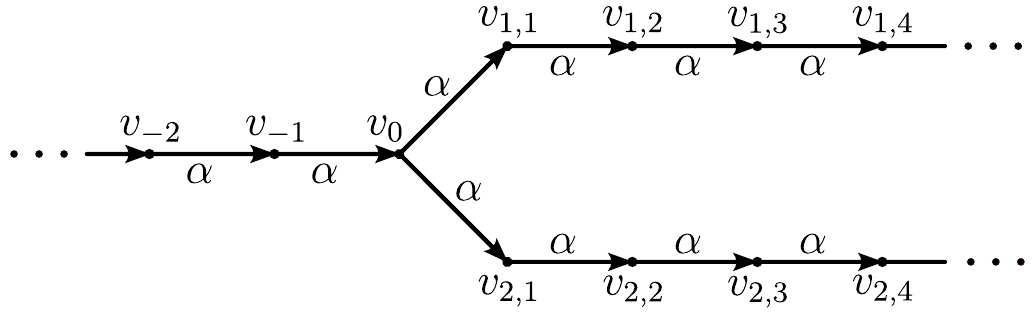}
    \caption{Set up of Example \ref{ex:fraise}}
    \label{fig:example:mean-LY:alpha}
\end{figure}    
Consider $(\lambda_v)_{v\in V}$ where $\lambda_v=\alpha$ for all $v\in V$. Fix $p\in\{0\}\cup[1,\infty)$ and consider $X:=\lambda_p(A,V)$, where the Köthe matrix $A=(a_{v,k})_{v\in V,\,k\in\N}$ is given by
\[
a_{v,k}=
\begin{cases}
j^k, & \text{if } v=v_{\ell,j}, \text{ for some } \ell\in\{1,2\},\ j\in\N,\\
1, & \text{if } v=v_{-j}, \text{ for some } j\in\N_0.
\end{cases}
\]
Note that $B_\lambda$ is well defined and continuous on $X$. We shall prove that the following assertions are equivalent:
\begin{enumerate}[label=$(\roman*)$]
    \item \label{example j^k i)} $B_{\lambda}$ is mean Li--Yorke chaotic;
    \item \label{example j^k ii)} $|\alpha|>1$.
\end{enumerate}
\ref{example j^k i)}$\Rightarrow$\ref{example j^k ii)}: Observe that, for every $n,k\in\N$ and $f\in X$,
\begin{equation}\label{eqExML1}
    \|B_\lambda^n f\|_k\leq
\begin{cases}
2^{1/p^*}|\alpha|^n\|f\|_k, & 1\leq p<\infty,\\
2|\alpha|^n\|f\|_k, & p=0,
\end{cases}
\end{equation}
where $p^*$ denotes the conjugate exponent of $p$, with the convention $1/p^*=0$ when $p=1$. Suppose that $|\alpha|\le 1$. By proposition \ref{prop:mean-LY-tree}, there exist $m \in \N$, a strictly increasing sequence $(N_k)_{k\in\N}$ of positive integers, and a sequence $(h_k)_{k\in\N}\subset E$ such that
\begin{equation}\label{eqExML2}
    \|h_k\|_{p(k)}>0 \quad \text{and} \quad \frac{1}{N_k\|h_k\|_{p(k)}}\sum_{j=1}^{N_k} \|B_\lambda^j h_k\|_{m} \geq  k,
\end{equation}
where $p(k) = \max\{k, m \}$, for all $k \in \N$. On the other hand, by \eqref{eqExML1}, for each $k\in\N$ we have
$$\frac{1}{N_k\|h_k\|_{p(k)}}\sum_{j=1}^{N_k} \|B_\lambda^j h_k\|_{m}\le 2,$$
contradicting \eqref{eqExML2}.

\ref{example j^k ii)}$\Rightarrow$\ref{example j^k i)}: Suppose that the generations are enumerated according to the vertex $v_0$. Take $S:=\N$ and, for each $n\in\N$, consider $f_n:= e_{v_{1,n}}-e_{v_{2,n}}$. Note that $B^n_{\lambda}f_n=0$ for every $n\in \N$. Therefore, condition \ref{teo:mean:cond:A} of Theorem \ref{thm:mean-LY-tree} holds. Now, set
\[
c:=\left\lceil\frac{2}{\log|\alpha|}\right\rceil,
\quad m_k:=c^2(k+1)^2,\quad N_k:=2m_k,
\quad\text{and}\quad
h_k:=\alpha^{-m_k}f_{N_k}, \quad\text{for all } k\in\N.
\]
Note that, for every $1\leq r\leq k$, we have
$$\|h_k\|_r\leq 2|\alpha|^{-m_k}(2m_k)^r\leq 2e^{-2c(k+1)^2+2ck(k+1)}=2e^{-2c(k+1)},$$
where we used $c\log|\alpha|\geq2$ and $\log(2c^2(k+1)^2)\leq2c(k+1)$. Consequently,
$$d(h_k,0)=\sum_{r=1}^{\infty}\frac{1}{2^r}\min\{1,\|h_k\|_r\}\leq 2e^{-2c(k+1)}+\sum_{r=k+1}^{\infty}\frac{1}{2^r}=2e^{-2c(k+1)}+2^{-k}<\frac1k.$$
On the other hand, for every $m_k\leq j<2m_k$, we have
\[
B_\lambda^j h_k=\alpha^{j-m_k}f_{2m_k-j},
\]
and hence
\[
\|B_\lambda^j h_k\|_r
=|\alpha|^{j-m_k}\|f_{2m_k-j}\|_r\geq1
\]
for every $r\in\N$. Thus,
\[
d(B_\lambda^j h_k,0)
=\sum_{r=1}^{\infty}\frac{1}{2^r}=1.
\]
Therefore,
$$\frac1{N_k}\sum_{j=1}^{N_k}d(B_\lambda^j h_k,0)\geq\frac1{2m_k}\sum_{j=m_k}^{2m_k-1}d(B_\lambda^j h_k,0)=\frac12.$$
Taking $\varepsilon=1/4$, condition \ref{teo:mean:cond:B} of Theorem \ref{thm:mean-LY-tree} holds. Therefore, $B_\lambda$ is mean Li--Yorke chaotic.
\end{example}

\subsection*{The case of admissible Banach sequence spaces}

Unlike in the sections on Li--Yorke chaos and distributional chaos, here we need to devote a separate subsection to the case of admissible Banach sequence spaces. The reason is that, even in the normed setting, mean Li--Yorke chaos with respect to the bounded metric $d$ is not, in general, equivalent to the usual notion of mean Li--Yorke chaos defined in terms of the norm. In what follows, we present the definitions in the normed context.

\begin{definition}
Let $T$ be a continuous linear operator on a Banach space $(X,\|\cdot\|)$. A pair $(x,y)\in X\times X$ is called a \emph{mean Li--Yorke pair} if
\[
\liminf_{N\to\infty}\frac{1}{N}\sum_{j=1}^N\|T^jx-T^jy\|=0 \quad \text{and}\quad 
\limsup_{N\to\infty}\frac{1}{N}\sum_{j=1}^N\|T^jx-T^jy\|>0.
\]
The operator $T$ is said to be ({\it densely}) \emph{mean Li--Yorke chaotic} if there exists an uncountable (dense) set $S\subset X$ such that every pair $(x,y)\in S\times S$ of distinct points is a mean Li--Yorke pair.
\end{definition}

\begin{definition}
Let $T$ be a continuous linear operator on a Banach space $(X,\|\cdot\|)$. A vector $x\in X$ is called \emph{absolutely mean semi-irregular} for $T$ if
\[
\liminf_{N\to\infty}\frac{1}{N}\sum_{j=1}^N\|T^jx\|=0
\quad\text{and}\quad
\limsup_{N\to\infty}\frac{1}{N}\sum_{j=1}^N\|T^jx\|>0.
\]
It is called \emph{absolutely mean irregular} for $T$ if
\[
\liminf_{N\to\infty}\frac{1}{N}\sum_{j=1}^N\|T^jx\|=0
\quad\text{and}\quad
\limsup_{N\to\infty}\frac{1}{N}\sum_{j=1}^N\|T^jx\|=\infty.
\]
\end{definition}

\begin{definition}
Let $T$ be a continuous linear operator on a Banach space $(X,\|\cdot\|)$. We say that $T$ satisfies the \emph{Mean Li--Yorke Chaos Criterion (MLYCC)} if there exists a subset $X_0\subset X$ with the following properties:
\begin{enumerate}[label=$(\roman*)$]
	\item\label{MLYCC i} $\displaystyle \liminf_{N\to\infty}\frac{1}{N}\sum_{j=1}^{N}\|T^{j}x\|=0$ for every $x\in X_0$;
	\item\label{MLYCC ii}  there exist a sequence $(y_k)_{k\in \N}$ in $\overline{\operatorname{span}(X_0)}$ and a subsequence $(N_k)_{k\in\N}\subset \N$ such that
	\[
	\frac{1}{N_k}\sum_{j=1}^{N_k}\|T^{j}y_k\|>k\|y_k\| \quad \text{for every } k\in\mathbb{N}.
	\]
\end{enumerate}
\end{definition}

The following lemma will used in the proof of our next theorem.

\begin{lemma}[{\cite[Theorems 5 and 9]{BBP}}] \label{LL}
Let $(X,\|\cdot\|)$ be a Banach space and $T:X\to X$ a continuous linear operator. Then, the following assertions are equivalent:
\begin{enumerate}[label=$(\alph*)$]
	\item\label{LL a} $T$ is mean Li--Yorke chaotic;
	\item\label{LL b} $T$ admits an \emph{absolutely mean irregular vector} $x \in X$;
	\item\label{LL c} $T$ satisfies the MLYCC.
\end{enumerate}
\end{lemma}

In the next theorem, we obtain a characterization analogous to Theorem \ref{thm:mean-LY-tree}, but in the normed case.

\begin{theorem} \label{teo:mean-LY-Banach}
Let $V$ be a directed tree with generations enumerated according to some fixed vertex $v_0$ and let $(X,\|\cdot\|)$ be an admissible Banach sequence space over $V$. Suppose that the weighted backward shift $B_\lambda$, with nonzero weights $\lambda := (\lambda_v)_{v \in V}$, is well-defined and continuous on $X$. Then $B_\lambda$ is mean Li--Yorke chaotic if and only if there exist a nonempty set $S\subset\{n\in\Z:\Gen_n\neq\varnothing\}$ and, for each $n\in S$, a vector $f_n\in X$ supported in $\Gen_n$ such that the following conditions hold:
\begin{enumerate}[label=$(\Alph*)$]
\item \label{teo:mean-LY-Banach:A} There exists an increasing sequence $(M_k)_{k\in\N}$ of positive integers such that
\[
\lim_{k\to\infty}\frac{1}{M_k}\sum_{j=1}^{M_k}\|B_\lambda^jf_n\|=0 \quad \text{for all } n\in S.
\]
\item \label{teo:mean-LY-Banach:B} There exists an increasing sequence $(N_k)_{k\in\N}$ of positive integers such that for each $k \in \N$ there is a finite subset $F_k\subset S$ and scalars $a_{j}\in \K$ for $j\in F_k$, with $\left\| \sum_{j\in F_k}a_jf_j\right\|$ $\ne 0$ and
\[
\frac{1}{N_k}\sum_{i=1}^{N_k}\left\| \sum_{j\in F_k}a_jB_\lambda^if_j \right\| > k\left\| \sum_{j\in F_k}a_jf_j \right\| \quad \text{for all } k\in\N.
\]
\end{enumerate}
\end{theorem}
\begin{proof} 
    ($\Rightarrow$): Suppose that $B_\lambda$ is mean Li--Yorke chaotic. By Lemma \ref{LL}, $B_{\lambda}$ admits an absolutely mean irregular vector $f\in X$, that is, 
    \begin{equation}\label{EqTeoNorm1}
        \liminf_{N\to\infty}\frac{1}{N}\sum_{j=1}^N\|B_{\lambda}^jf\|=0\quad\text{and}\quad\limsup_{N\to\infty}\frac{1}{N}\sum_{j=1}^N\|B_{\lambda}^jf\|=\infty.
    \end{equation}
    By \eqref{EqTeoNorm1}, there is an increasing sequence of positive integers $(M_k)_{k\in \N}$ such that
    \begin{equation} \label{Mean-Semi}
        \lim_{k \to \infty } \frac{1}{M_k} \sum^{M_k}_{j= 1} \|B_\lambda^j f \| = 0.
    \end{equation}
    Define $S:=\{n\in\mathbb Z:\Gen_n\neq \emptyset\}$. By Lemma \ref{LM1}, the family of operators $(\Phi_{n})_{n\in \Z}$ given by 
        $$\Phi_{n}: g\in X \to g\chi_{\operatorname{Gen}_n}\in X, \quad \text{for } n\in \Z,$$
        is equicontinuous. Therefore, there exists $C>0$ such that
\[
\|\Phi_n(g)\|\leq C\|g\|, \quad \text{for every }g\in X\text{ and }n\in\Z.
\]
Thus, by Lemma \ref{LM2} and \eqref{Mean-Semi}, we obtain
\[
\frac{1}{M_k}\sum_{j=1}^{M_k}\|B_\lambda^j(f\chi_{\Gen_n})\|=\frac{1}{M_k}\sum_{j=1}^{M_k}\|\Phi_{n-j}(B_\lambda^jf)\| \leq \frac{C}{M_k}\sum_{j=1}^{M_k}\|B_\lambda^jf\| \xrightarrow{k\to\infty} 0,
\]
for all $n\in S$. Thus, taking $f_n:=f\chi_{\Gen_n}$ for each $n\in S$, condition \ref{teo:mean-LY-Banach:A} holds. On the other hand, by \eqref{EqTeoNorm1}, for each $k\in \N$, there exists $N_k\in \N$ such that
    \begin{equation}\label{EqTeoNorm2}
        \frac{1}{N_k} \sum^{N_k}_{j= 1} \|B_\lambda^j f \|> k\left\| f \right\|+2.
    \end{equation}
     Take the integers $N_k$, so that $(N_k)_{k\in\N}$ forms a subsequence of $\N$. Since $f=\sum_{n\in S} f_n$, and $B_{\lambda}$ is continuous, for each $k\in \N$ there exists a finite subset $F_k \subset S$ such that 
\begin{equation}\label{EqTeoNorm3}
     \left| \left\|B_{\lambda}^{i}\left( \sum_{j\in F_k}^{}f_j \right) \right\| -\left\| B_{\lambda}^{i}f \right\|\right|<\frac{1}{2k},\quad \text{for all }i\in\{0,\cdots, N_k\}.
\end{equation}
Thus, by \eqref{EqTeoNorm2} and \eqref{EqTeoNorm3}, we obtain
$$\frac{1}{N_k}\sum_{i=1}^{N_k}\left\| \sum_{j\in F_k}B_\lambda^if_j \right\|>\frac{1}{N_k}\sum_{i=1}^{N_k}\left\| B_{\lambda}^{i}f \right\|-\frac{1}{2k}> k\left\| f \right\|+2-\frac{1}{2k}>k\left\| \sum_{j\in F_k}^{}f_j \right\|$$
for each $k\in\N$. Therefore, condition \ref{teo:mean-LY-Banach:B} holds.

($\Leftarrow$): Suppose that conditions \ref{teo:mean-LY-Banach:A} and \ref{teo:mean-LY-Banach:B} hold. Define $X_0 := \{ f_n : n\in S\}$. By condition \ref{teo:mean-LY-Banach:A}, item \ref{MLYCC i} of MLYCC holds. Now, taking 
$$y_k:=\sum_{j\in F_k}^{}a_jf_j \quad\text{for each }k\in\N,$$
we conclude that condition \ref{teo:mean-LY-Banach:B} implies item \ref{MLYCC ii} of MLYCC. Therefore, by Lemma \ref{LL}, $B_{\lambda}$ is mean Li--Yorke chaotic.
\end{proof}
    Next, we make the same observation, now concerning Theorem~\ref{teo:mean-LY-Banach}, as the one made in Remark~\ref{rmkLY}.
\begin{remark}\label{rmkD22}
    In Theorem \ref{teo:mean-LY-Banach}, for each $n\in S$, we chose a single vector $f_n \in X$ supported in $\Gen_n$, such that the family $(f_n)_{n\in S}$ satisfies conditions \ref{teo:mean-LY-Banach:A} and \ref{teo:mean-LY-Banach:B}. From the proof, we see that this was not necessary, that is, we could restate the theorem by taking, for each $n\in S$, not only a single vector, but rather a sequence $(f_{n,k})_{k\in \N}\subset X$ of vectors supported in $\Gen_n$, such that the family of vectors $\left\{ f_{n,k}:n\in S\text{ and }k\in \N \right\}$ satisfies conditions \ref{teo:mean-LY-Banach:A} and \ref{teo:mean-LY-Banach:B} of the theorem. 
\end{remark}

The following corollary improves \cite[Theorem 39]{joaojmaa} by requiring no additional basis assumption other than being a Schauder basis.

\begin{corollary} \label{cor:joao-mean-LY}
Let $(X,\|\cdot\|)$ be a Banach sequence space over $\Z$ in which $(e_n)_{n\in\Z}$ is a Schauder basis. Suppose that the bilateral weighted backward shift $B_w$, with nonzero weights $w:=(w_n)_{n\in\Z}$, is well-defined and continuous on $X$. Then $B_w$ is mean Li--Yorke chaotic if and only if the following conditions hold:
\begin{enumerate}[label=$(\Alph*)$]
\item \label{cor:joao:A} $\liminf_{N\to\infty}\frac{1}{N}\sum_{j=1}^N\|w_{-j}\cdots w_{-1}e_{-j}\|=0$;
\item \label{cor:joao:B} there exists an increasing sequence $(N_k)_{k\in\N}$ of positive integers such that, for each $k\in\N$, there are $r:=r(k)\in\N$ and scalars $b_{-r,k},\dots, b_{r,k}\in\K$ for which $\sum_{s=-r}^{r}b_{s,k}e_s \neq 0$ and
\[
\frac{1}{N_k\left\|\sum_{s=-r}^{r}b_{s,k}e_s\right\|}\sum_{j=1}^{N_k}\left\|\sum_{s=-r}^{r}b_{s,k}w_{s-j}\cdots w_{s-1}e_{s-j}\right\|> k.
\]
\end{enumerate}
\end{corollary}
\begin{proof}
Recall that, by Example \ref{ExFNZ}, we may regard $B_w$ as the weighted backward shift $B_\lambda$ introduced in \eqref{eqq0} on the admissible Fréchet sequence space over $V=\Z$. For this, it suffices to define $\lambda_{n+1}:=w_n$ for every $n\in\Z$. As usual, we consider $\Gen_n:=\{n\}$ for $n\in \Z$. We shall use the elementary fact that, for every $m\in\Z$ and every increasing sequence $(M_k)_{k\in\N}$ of positive integers,
\begin{equation} \label{eq:em:e0}
\lim_{k\to\infty}\frac{1}{M_k} \sum_{j=1}^{M_k} \|B_w^je_m\| = 0
\quad \Longleftrightarrow \quad \lim_{k\to\infty} \frac{1}{M_k} \sum_{j=1}^{M_k} \|B_w^je_0\| = 0.
\end{equation}

($\Rightarrow$): Suppose that $B_w$ is mean Li--Yorke chaotic. By Theorem \ref{teo:mean-LY-Banach}, there exist $n\in \Z$ and a nonzero vector $f_n$ supported in $\Gen_n=\{n\}$ satisfying condition \ref{teo:mean-LY-Banach:A} of Theorem \ref{teo:mean-LY-Banach}. Since $f_n$ is a nonzero scalar multiple of $e_n$, \eqref{eq:em:e0} yields \ref{cor:joao:A}. Moreover, it is easy to see that condition \ref{teo:mean-LY-Banach:B} of Theorem \ref{teo:mean-LY-Banach} is precisely \ref{cor:joao:B}.

($\Leftarrow$): Suppose that \ref{cor:joao:A} and \ref{cor:joao:B} hold. By \ref{cor:joao:A}, there is an increasing sequence $(M_k)_{k\in\N}$ of positive integers such that
\begin{equation} \label{eq:em:e1}
    \lim_{k\to\infty} \frac{1}{M_k} \sum_{j=1}^{M_k} \|B_w^je_0\| = 0.
\end{equation}
Take $S=\Z$ and define $f_n:=e_n$ for all $n\in S$. Thus, condition \ref{teo:mean-LY-Banach:A} of Theorem \ref{teo:mean-LY-Banach} follows from \eqref{eq:em:e0} and \eqref{eq:em:e1}. On the other hand, taking $F_k:=\{-r(k),\cdots, r(k)\}$ and $a_j:=b_{j,k}$ for $j\in F_k$, item \ref{cor:joao:B} implies condition \ref{teo:mean-LY-Banach:B} of Theorem \ref{teo:mean-LY-Banach}. Therefore, $B_w$ is mean Li--Yorke chaotic.
\end{proof}
The following corollary is a direct application of Theorem \ref{teo:mean-LY-Banach} to spaces of the form $c_0(\nu,V)$ and $\ell^p(\nu,V)$, $1\leq p<\infty$. We observe that the next corollary can also be stated in the form mentioned in Remark~\ref{rmkD22}.
\begin{corollary}\label{cor:mean-LY-lp-c0}
Let $V$ be a directed tree and let $\nu:=(\nu_v)_{v\in V}\subset (0,\infty)$, and consider $X:=c_0(\nu,V)$ or $X:=\ell^p(\nu,V)$ for some $1\leq p<\infty$. Suppose that the weighted backward shift $B_\lambda$, with nonzero weights $\lambda := (\lambda_v)_{v \in V}$, is well-defined and continuous on $X$. Then $B_\lambda$ is mean Li--Yorke chaotic if and only if there exist a nonempty set $S\subset\{n\in\Z:\Gen_n\neq\varnothing\}$ and scalars $c_v\in\K$, for $v\in\Gen_n$ and $n\in S$, such that $\sum_{v\in\Gen_n}c_ve_v\in X$ for all $n\in S$, and the following conditions hold:
\begin{enumerate}[label=$(\Alph*)$]
\item \label{cor:mean-LY-lp-c0 A} there exists an increasing sequence $(M_k)_{k\in\N}$ of positive integers such that, for every $n\in S$,
$$\begin{cases}
 \displaystyle \lim_{k\to\infty} \frac{1}{M_k} \sum_{j=1}^{M_k} \sup_{v\in\Gen_{n-j}} \left|\sum_{u\in\Chi^j(v)} \nu_v\lambda(v\to u)c_u\right| = 0,   & \text{if }\,X=c_0(\nu,V), \\
 \displaystyle \lim_{k\to\infty} \frac{1}{M_k} \sum_{j=1}^{M_k} \left(\sum_{v\in\Gen_{n-j}} \left|\sum_{u\in\Chi^j(v)}\nu_v \lambda(v\to u)c_u\right|^p\right)^{1/p}=0,   & \text{if }\,X=\ell^p(\nu,V).
\end{cases}$$

\item \label{cor:mean-LY-lp-c0 B} there exists an increasing sequence $(N_k)_{k\in \N}$ of positive integers such that for each $k \in \N$ there is a finite subset $F_k\subset S$ and scalars $a_{j}\in \K$ for $j\in F_k$, with $g_k:=\sum_{n\in F_k}a_n\sum_{v\in\Gen_n}c_ve_v\ne 0$ and
$$\begin{cases}
 \displaystyle  \frac{1}{N_k} \sum_{j=1}^{N_k}\max_{n\in F_k} \sup_{v\in\Gen_{n-j}} \left|\sum_{u\in\Chi^j(v)}\nu_va_n \lambda(v\to u)c_u\right| > k\left\| g_k \right\|,  & \text{if }\,X=c_0(\nu,V), \\
 \displaystyle  \frac{1}{N_k} \sum_{j=1}^{N_k} \left(\sum_{n\in F_k}\sum_{v\in\Gen_{n-j}} \left|\sum_{u\in\Chi^j(v)} \nu_va_n\lambda(v\to u)c_u\right|^p\right)^{1/p}> k\left\| g_k \right\|,   & \text{if }\,X=\ell^p(\nu,V).
\end{cases}$$
\end{enumerate}
\end{corollary}

\begin{example}
        Fix $N\in \N$ such that $N\ge 2$. Consider the unrooted directed $N$-adic tree $V$ centered at the vertex $v_0$. Consider $\nu:= (\nu_v)_{v\in V}\subset (0,\infty)$ such that there exists a sequence $(b_n)_{n\in\Z}\subset (0,\infty)$ where, for each $n\in\Z$,
        $$\nu_v=b_n\quad \text{for all }v\in \Gen_n.$$
        Consider $X:=c_0(\nu,V)$ or $X:=\ell^p(\nu,V)$ for some $1\leq p<\infty$. Fix $\alpha\in\K\setminus\{0\}$ and consider the weighted backward shift $B_\lambda$ with weights $\lambda_v=\alpha$ for every $v\in V$. We assume that
        $$\sup_{n\in \Z}\frac{b_{n-1}}{b_n}<\infty.$$
        Under this assumption, $B_\lambda$ is well-defined and continuous on $X$. Define
$$C:=
\begin{cases}
N|\alpha|, & \text{if }X=c_0(\nu,V),\\
|\alpha|, & \text{if }X=\ell^1(\nu,V),\\
N^{\frac{1}{p^*}}|\alpha|, & \text{if }X=\ell^p(\nu,V)\text{ for }p\in (1,\infty),
\end{cases}$$
where $p^*=\frac{p}{p-1}$. We will prove that the following assertions are equivalent:
\begin{enumerate}[label=$(\roman*)$]
    \item \label{example N alpha - i)} $B_{\lambda}$ is mean Li--Yorke chaotic;
    \item  \label{example N alpha - ii)} there exists a subsequence $(N_k)_{k\in\N}\subset \N$ such that for each $k \in \N$ there is a finite subset $F_k\subset \Z$ and  scalars $a_j>0$ for $j\in F_k$ where
    $$\begin{cases}
\displaystyle \frac{1}{N_k}\sum_{n=1}^{N_k}C^n\frac{\max\limits_{j\in F_k}a_{j}b_{j-n}}{\max\limits_{j\in F_k}a_{j}b_{j}}>k, & \text{if }\,X=c_0(\nu, V),\\
\displaystyle \frac{1}{N_k}\sum_{n=1}^{N_k}C^n\frac{\left( \sum\limits_{j\in F_k}a_{j}^pb_{j-n}^p \right)^{1/p}}{\left( \sum\limits_{j\in F_k}a_{j}^pb_{j}^p \right)^{1/p}}>k, & \text{if }\,X=\ell^p(\nu, V).
\end{cases}$$
\end{enumerate}
\ref{example N alpha - i)}$\Rightarrow$\ref{example N alpha - ii)}: By condition \ref{cor:mean-LY-lp-c0 B} of Corollary \ref{cor:mean-LY-lp-c0}, there exist $S\subset\Z$ and an increasing sequence $(N_k)_{k\in\N}$ of positive integers such that, for each $k\in\N$, there are a finite subset $F_k\subset S$ and vectors $f_j\in X$, supported on $\Gen_j$ for every $j\in F_k$, with $\sum_{j\in F_k}f_j\ne0$ and
\begin{equation}\label{eq:mean-ex-1}
\begin{cases}
\displaystyle \frac{1}{N_k}\sum_{n=1}^{N_k}\frac{\max\limits_{j\in F_k}\|B_\lambda^nf_j\|}{\max\limits_{j\in F_k}\|f_j\|}>k, &\text{if }\, X=c_0(\nu,V),\\
\displaystyle \frac{1}{N_k}\sum_{n=1}^{N_k}\frac{\left(\sum\limits_{j\in F_k}\|B_\lambda^nf_j\|^p\right)^{1/p}}{\left(\sum\limits_{j\in F_k}\|f_j\|^p\right)^{1/p}}>k, & \text{if }\,X=\ell^p(\nu,V).
\end{cases}
\end{equation}
Fix $k\in\N$. For each $j\in F_k$, define
$$
a_j:=
\begin{cases}
\displaystyle \sup_{u\in\Gen_j}|f_j(u)|, & \text{if }\,X=c_0(\nu,V),\\
\displaystyle \left(\sum_{u\in\Gen_j}|f_j(u)|^p\right)^{1/p}, & \text{if }\, X=\ell^p(\nu,V).
\end{cases}
$$
Removing the indices for which $a_j=0$, if necessary, we may assume that $a_j>0$ for every $j\in F_k$. Note that
\begin{equation}\label{eq:mean-ex-2}
\|f_j\|=a_jb_j
\end{equation}
for every $j\in F_k$. Moreover, using the triangle inequality if $X=c_0(\nu,V)$ and Hölder's inequality if $X=\ell^p(\nu,V)$, we get
\begin{equation}\label{eq:mean-ex-3}
\|B_\lambda^nf_j\|\le C^na_jb_{j-n}
\end{equation}
for every $j\in F_k$ and $n\in\N$. Therefore, combining \eqref{eq:mean-ex-1}--\eqref{eq:mean-ex-3}, we obtain item \ref{example N alpha - ii)}.

\ref{example N alpha - ii)}$\Rightarrow$\ref{example N alpha - i)}: For this implication, we use the version of Corollary \ref{cor:mean-LY-lp-c0} given in Remark \ref{rmkD22}. Suppose that item \ref{example N alpha - ii)} holds. For each $k\in\N$ and $j\in F_k$, choose $z_{k,j}\in\Gen_{j-N_k-1}$ and two distinct vertices $x_{k,j},y_{k,j}\in\Chi(z_{k,j})$. Define
$$
h_{k,j}:=
\begin{cases}
\displaystyle a_j\left(\sum_{u\in\Chi^{N_k}(x_{k,j})}e_u-\sum_{u\in\Chi^{N_k}(y_{k,j})}e_u\right), &\text{if }\, X=c_0(\nu,V),\\
\displaystyle \frac{a_j}{(2N^{N_k})^{1/p}}\left(\sum_{u\in\Chi^{N_k}(x_{k,j})}e_u-\sum_{u\in\Chi^{N_k}(y_{k,j})}e_u\right), &\text{if }\, X=\ell^p(\nu,V).
\end{cases}
$$
Then $h_{k,j}\in X$, $\supp(h_{k,j})\subset\Gen_j$, and $B_\lambda^{N_k+1}h_{k,j}=0$. Hence, condition \ref{cor:mean-LY-lp-c0 A} of Corollary \ref{cor:mean-LY-lp-c0} holds. Moreover, for every $1\le n\le N_k$, we have
$$
B_\lambda^nh_{k,j}=
\begin{cases}
\displaystyle a_j\alpha^nN^n\left(\sum_{u\in\Chi^{N_k-n}(x_{k,j})}e_u-\sum_{u\in\Chi^{N_k-n}(y_{k,j})}e_u\right), &\text{if }\, X=c_0(\nu,V),\\
\displaystyle \frac{a_j\alpha^nN^n}{(2N^{N_k})^{1/p}}\left(\sum_{u\in\Chi^{N_k-n}(x_{k,j})}e_u-\sum_{u\in\Chi^{N_k-n}(y_{k,j})}e_u\right), &\text{if }\, X=\ell^p(\nu,V).
\end{cases}
$$
Thus, for every $1\le n\le N_k$,
\begin{equation}\label{eq:mean-ex-norms}
\|h_{k,j}\|=a_jb_j,\quad \|B_\lambda^nh_{k,j}\|=C^na_jb_{j-n}.
\end{equation}
Hence, by item \ref{example N alpha - ii)} and \eqref{eq:mean-ex-norms}, we have 
\begin{equation}\label{eq:mean-ex-ratio}
\begin{cases}
\displaystyle \frac{1}{N_k}\sum_{n=1}^{N_k}\frac{\max\limits_{j\in F_k}\|B_\lambda^nh_{k,j}\|}{\max\limits_{j\in F_k}\|h_{k,j}\|}
> k, &\text{if }\, X=c_0(\nu,V),\\
\displaystyle \frac{1}{N_k}\sum_{n=1}^{N_k}\frac{\left(\sum\limits_{j\in F_k}\|B_\lambda^nh_{k,j}\|^p\right)^{1/p}}{\left(\sum\limits_{j\in F_k}\|h_{k,j}\|^p\right)^{1/p}}> k, &\text{if }\, X=\ell^p(\nu,V).
\end{cases}
\end{equation}
Therefore, condition \ref{cor:mean-LY-lp-c0 B} of Corollary \ref{cor:mean-LY-lp-c0} also holds, what implies that $B_w$ is mean Li--Yorke chaotic.
\end{example}

Given a continuous linear operator $T$ on a Banach space $(X,\|\cdot\|)$, throughout this subsection, the sets $\mathrm{MAsym}(T,0)$ and $\mathrm{MProx}(T,0)$ are defined as in \eqref{MAsym} and \eqref{MProx}, respectively, but with the norm $\|\cdot\|$ in place of the metric. Similarly, replacing the metric with the norm in \eqref{meansens}, we obtain the definition of mean sensitivity that we will use throughout this subsection. Next, we introduce a concept closely related to the definition of mean Li--Yorke chaos in the Banach space setting.
\begin{definition}
Let $T:X\to X$ be a continuous linear operator on a Banach space $(X,\|\cdot\|)$. We say that $T$ is \emph{absolutely Cesàro bounded} if there exists $C>0$ such that
\[
\frac{1}{N} \sum_{j=1}^N \|T^jx\| \leq C\|x\| \quad\text{for every $x\in X$ and $N\in\N$}.
\]
\end{definition}
Let $V$ be a directed tree, let $X$ be an admissible Banach sequence space over $V$, and let $B_\lambda$, with nonzero weights $\lambda:=(\lambda_v)_{v\in V}$, be a well-defined continuous weighted backward shift on $X$. We say that $\lambda$ satisfies condition \eqref{cond:C2-norm} if there exists $v_0\in V$ such that
\begin{equation}\tag{$C2_{\|\cdot\|}$} \label{cond:C2-norm}
\liminf_{N\to\infty}\frac{1}{N} \sum_{j=1}^N \|\lambda\bigl(\Par^j(v_0)\to v_0\bigr)e_{\Par^j(v_0)}\| = 0.
\end{equation}
For our next results, we will need the following lemmas.
\begin{lemma}[{\cite[Theorem 22]{BBP}}]\label{Th: BBP-22}
    Let $(X,\|\cdot\|)$ be a Banach space. If $T:X\to X$ is a continuous linear operator and  $\MProx(T, 0)$ is dense in $X$, then the following assertions are equivalent:
    \begin{enumerate}[label=$(\roman*)$]
        \item \label{Th: BBP-22: (i)}$T$ is mean Li--Yorke chaotic;
        \item \label{Th: BBP-22: (ii)} $T$ has a residual set of absolutely mean irregular vectors;
        \item \label{Th: BBP-22: (iii)} $T$ is not absolutely Cesàro bounded;
        \item \label{Th: BBP-22: (iv)}$\displaystyle\limsup_{N \to \infty} \frac{1}{N} \sum^N_{j =1} \|T^jy_0\| > 0$ for some $y_0 \in X$.
    \end{enumerate}
\end{lemma}

The following lemma comes from a combination of Lemma \ref{Th: BBP-22} with \cite[Theorem 17 and Remark 23]{BBP}.

\begin{lemma}\label{BBP-rem-23}
Let $X$ be a separable Banach space and let $T$ be a continuous linear operator on $X$. Suppose that $\MProx(T,0)$ is dense in $X$. Then $T$ is mean Li--Yorke chaotic if and only if $T$ is densely mean Li--Yorke chaotic.
\end{lemma}

The following dichotomy result is the analogous of Proposition \ref{prop:mean-LY-dicotomia} for the normed case.

\begin{proposition} \label{prop:norma-dicotmia}
Let $V$ be a directed tree and let $X$ be an admissible Banach sequence space over $V$. Suppose that the weighted backward shift $B_\lambda$, with nonzero weights $\lambda := (\lambda_v)_{v \in V}$, is well-defined and continuous on $X$. Suppose that condition \eqref{cond:C2-norm} holds. Then $\MProx(B_\lambda,0)$ is dense in $X$, and exactly one of the following assertions holds:
\begin{enumerate}[label=$(\alph*)$]
\item \label{prop:norma-dicotmia:a} $\MAsym(B_\lambda,0)=X$;
\item \label{prop:norma-dicotmia:b} the absolutely mean irregular vectors form a residual subset of $X$.
\end{enumerate}
\end{proposition}
\begin{proof}
    We shall prove that $c_{00}(V)\subset\MProx(B_\lambda,0)$. If $V$ is rooted, the inclusion is trivial. Suppose that $V$ is unrooted. By condition \eqref{cond:C2-norm}, there exist $v_0\in V$ and an increasing sequence $(M_k)_{k\in \N}$ of positive integers such that 
    \begin{equation}\label{a1a2a3}
        \lim_{k \to \infty} \frac{1}{M_k} \sum_{j=1}^{M_k} \|B_{\lambda}^j e_{v_0}\| = 0.
    \end{equation}
    Take $f \in c_{00}(V)$. Then there exists a finite set $F \subset V$ such that $f = \sum_{v \in F} f(v) e_v$. For each $v \in F$, we can take $p_v, q_v \in \N_0$ such that $\Par^{p_v}(v_0) = \Par^{q_v}(v)$. As in the proof of Proposition \ref{prop:mean-LY-dicotomia}, for all $j\ge q_v+1$, we can write
    \begin{equation}\label{a1a2a31}
        B_{\lambda}^{j  }e_v = c_v B_{\lambda}^{j-q_v+p_v}e_{v_0},
    \end{equation}
    for some constant $c_v \not= 0$, which depends only on $v$. Fix $I\in\N$ such that $I\ge \max\left\{ q_v:v\in F \right\}$. Hence, by \eqref{a1a2a31}, for all $v\in F$,
    \begin{align*}
        \sum^{M_k + I}_{i=1} \|B_\lambda^i e_v\| &= \overset{=:C}{\overbrace{\sum^{ I}_{i=1} \|B_\lambda^i e_v\|}} + \sum^{M_k + I}_{i=1 + I} \|B_\lambda^i e_v\| \\
        &= C + \sum^{M_k}_{j=1} \|B_\lambda^{j + I} e_v\| \\
        & \leq C + |c_v| \|B_\lambda^{I + p_v - q_v}\|\sum^{M_k}_{j=1} \|B_\lambda^{j} e_{v_0}\|.
    \end{align*}
    Thus, by \eqref{a1a2a3}, we conclude that
    \[\lim_{k \to \infty } \frac{1}{M_k+I}\sum^{M_k + I}_{i=1} \|B_\lambda^i f\| = 0.\]
    Then $c_{00}(V)\subset\MProx(B_\lambda,0)$. Therefore, $\MProx(B_\lambda, 0)$ is dense in $X$. It remains to prove the final dichotomy. Of course \ref{prop:norma-dicotmia:a} and \ref{prop:norma-dicotmia:b} are mutually exclusive. Suppose that \ref{prop:norma-dicotmia:a} is false. Then, there exists $f \in X$ such that 
    
    \[\limsup_{N\to\infty}\frac{1}{N}\sum_{j=1}^N\|B_\lambda^jf\|>0.\]
    Thus, since $\MProx(B_\la,0)$ is dense, by Lemma \ref{Th: BBP-22}, the absolutely mean irregular vectors form a residual subset of $X$. This concludes the proof.
\end{proof}
For obtaining the analogous of Proposition \ref{prop:mean-LY-equiv} in the normed case, we will need the following result from \cite{JiangLi}.

\begin{lemma}[{\cite[Proposition 4.33]{JiangLi}}]\label{china-Prop-4.33}
    Let $X$ be a Banach space and let $T$ be a continuous linear operator on $X$. Then $T$ is mean sensitive if and only if it is not absolutely Cesàro bounded.
\end{lemma}

\begin{proposition} \label{prop:norm-mean-LY-equivalencias}
Let $V$ be a directed tree and let $X$ be an admissible Banach sequence space over $V$. Suppose that the weighted backward shift $B_\lambda$, with nonzero weights $\lambda := (\lambda_v)_{v \in V}$, is well-defined and continuous on $X$. Suppose that condition \eqref{cond:C2-norm} holds. Then the following assertions are equivalent:
\begin{enumerate}[label=$(\roman*)$]
\item \label{prop:norm-mean-LY-equivalencias:i} $B_\lambda$ is mean Li--Yorke chaotic;
\item \label{prop:norm-mean-LY-equivalencias:ii} $B_\lambda$ is densely mean Li--Yorke chaotic;
\item\label{prop:norm-mean-LY-equivalencias:iii} $B_\lambda$ is mean sensitive;
\item\label{prop:norm-mean-LY-equivalencias:iv} $B_\lambda$ is not absolutely Cesàro bounded;
\item\label{prop:norm-mean-LY-equivalencias:v} there exists $f\in X$ such that $\displaystyle \limsup_{N\to\infty}\frac{1}{N} \sum_{j=1}^N \|B_\lambda^jf\|>0$;
\item\label{prop:norm-mean-LY-equivalencias:vi} there are an increasing sequence $(N_k)_{k\in\N}$ of positive integers and a sequence $(h_k)_{k\in\N}\subset c_{00}(V)$ of nonzero vectors such that
\[
\frac{1}{N_k} \sum_{j=1}^{N_k} \|B_\lambda^jh_k\|>k\|h_k\|  \quad \text{for every }k\in\N.
\]
\end{enumerate}
Moreover, whenever one of these conditions holds, the absolutely mean irregular vectors for $B_\lambda$ form a residual subset of $X$.
\end{proposition}
\begin{proof}
Since \eqref{cond:C2-norm} holds, it follows from Proposition \ref{prop:norma-dicotmia} that $\MProx(B_\lambda,0)$ is dense in $X$.

\ref{prop:norm-mean-LY-equivalencias:i}$\Rightarrow$\ref{prop:norm-mean-LY-equivalencias:ii}: Since $\MProx(B_\lambda,0)$ is dense, by Lemma \ref{BBP-rem-23}, we conclude that $B_{\lambda}$ is densely mean Li--Yorke chaotic.

\ref{prop:norm-mean-LY-equivalencias:ii}$\Rightarrow$\ref{prop:norm-mean-LY-equivalencias:iii}: Since $\MProx(B_\lambda,0)$ is dense, by Lemma \ref{Th: BBP-22}, we conclude that $B_{\lambda}$ is not absolutely Cesàro bounded. Thus, by Lemma \ref{china-Prop-4.33}, $B_{\la}$ is mean sensitive.

\ref{prop:norm-mean-LY-equivalencias:iii}$\Rightarrow$\ref{prop:norm-mean-LY-equivalencias:iv}: This implication follows from Lemma \ref{china-Prop-4.33}.

\ref{prop:norm-mean-LY-equivalencias:iv}$\Rightarrow$\ref{prop:norm-mean-LY-equivalencias:v}: Since $\MProx(B_\lambda,0)$ is dense, by Lemma \ref{Th: BBP-22}, we obtain item \ref{prop:norm-mean-LY-equivalencias:v}.

\ref{prop:norm-mean-LY-equivalencias:v}$\Rightarrow$\ref{prop:norm-mean-LY-equivalencias:vi}: By hypothesis and Lemma \ref{Th: BBP-22}, we may choose $g\in X$ and an increasing sequence $(N_k)_{k\in\N}$ of positive integers such that
\[
\frac{1}{N_k} \sum_{j=1}^{N_k} \|B_\lambda^jg\| > 4k\|g\|
\quad\text{for every }k\in\N.
\]
Since $c_{00}(V)$ is dense in $X$ and $B_\lambda$ is continuous, for each $k\in\N$ we may choose $h_k\in c_{00}(V)$ such that
\[
\|h_k-g\| < \|g\|
\quad\text{and}\quad
\frac{1}{N_k} \sum_{j=1}^{N_k} \|B_\lambda^j(h_k-g)\| < k\|g\|.
\]
Then $0<\|h_k\|<2\|g\|$ and
\[
\frac{1}{N_k} \sum_{j=1}^{N_k} \|B_\lambda^jh_k\| > 3k\|g\| > k\|h_k\|,
\]
which gives \ref{prop:norm-mean-LY-equivalencias:vi}.

\ref{prop:norm-mean-LY-equivalencias:vi}$\Rightarrow$\ref{prop:norm-mean-LY-equivalencias:i}: This implication follows from MLYCC by taking $X_0:=c_{00}(V)$ and by item \ref{prop:norm-mean-LY-equivalencias:vi}.
\end{proof}

If $V$ is rooted, then condition \eqref{cond:C2-norm} is automatically satisfied. Hence, we obtain the following corollary.
\begin{corollary}
    Let $V$ be a rooted directed tree. Let $X$ be an admissible Banach sequence space over $V$. Suppose that the weighted backward shift $B_\lambda$, with nonzero weights $\lambda := (\lambda_v)_{v \in V}$, is well-defined and continuous on $X$. Then the assertions \ref{prop:norm-mean-LY-equivalencias:i}-\ref{prop:norm-mean-LY-equivalencias:vi} from Proposition \ref{prop:norm-mean-LY-equivalencias} are equivalent.
\end{corollary}

\section{Distributional chaos} \label{sec:dist}

Schweizer and Smítal \cite{23} introduced the following notion, which is a stronger version of the original Li--Yorke concept.

\begin{definition}
		Given a metric space $(M,d)$, a map $f : M \to M$ is said to be ({\it densely}) {\em distributionally chaotic} if there exist
		an uncountable (dense) set $\Gamma \subset M$ and $\varepsilon > 0$ such that each pair $(x,y)$ of distinct points in $\Gamma$
		is a {\em distributionally chaotic pair for $f$}, in the sense that
		\[
		\liminf_{k\to \infty}\frac{\text{card}\left( \{n \in \left\{ 1, \cdots, k \right\} : d(f^n(x),f^n(y)) < \varepsilon\} \right)}{k}=0
		\]
		and
		\[
		\limsup_{k\to \infty}\frac{\text{card}\left( \{n \in \left\{ 1, \cdots, k \right\} : d(f^n(x),f^n(y)) < \tau \} \right)}{k}=1, \ \text{ for all } \tau > 0,
		\]
		where $\text{card}(A)$ denotes the cardinality of the set $A\subset \N.$ In this case, $\Gamma$ is said to be a distributionally scrambled set.
\end{definition}
	We also use the following notation: let $A\subset \N$, then we define
	$$\udens(A):=\limsup_{N \to \infty}\frac{\text{card}(\left\{ 1, \cdots, N \right\}\cap A)}{N}\quad \text{and} \quad \ldens(A):=\liminf_{N \to \infty}\frac{\text{card}(\left\{ 1, \cdots, N \right\}\cap A)}{N}.$$
	To prove our results in this section, we need the following definition and lemma from \cite{Nilson0}. We emphasize that \cite{Nilson0} is a seminal work in the theory of linear dynamics, providing a general study of the notion of distributional chaos for operators on Fréchet spaces.
    
	\begin{definition}Let $T$ be a continuous linear operator on a Fréchet space $X$. We say that $T$ satisfies the \emph{Distributional Chaos Criterion} (DCC) if there exist sequences $(x_k)_{k\in\N}, (y_k)_{k\in\N}$ in $X$ such that:
		\begin{enumerate}[label=$(\alph*)$]
			\item \label{dcc a} there exists $A \subset \mathbb{N}$ with $\udens(A) = 1$ such that $\lim_{n \in A} T^n x_k = 0$, for all $k \in \N$;
			\item \label{dcc b} $y_k \in \overline{\Span\{x_n : n \in \mathbb{N}\}}$, $\lim_{k \to \infty} y_k = 0$ and there exist $\varepsilon > 0$ and an increasing sequence $(N_k)_{k\in \N}$ in $\mathbb{N}$ such that
			\[
			\card\{1 \le j \le N_k : d(T^j y_k, 0) > \varepsilon\} \ge N_k(1 - k^{-1}), \quad \text{for all }k \in \mathbb{N}.
			\]
		\end{enumerate}
	\end{definition}
	\begin{lemma}[{\cite[Theorem 12]{Nilson0}}] \label{T0}
		Let $T$ be a continuous linear operator on a Fréchet space $Y$. Then, the following statements are equivalent:
		\begin{enumerate}[label=$(\alph*)$]
			\item\label{T0 a} $T$ is distributionally chaotic;
			\item\label{T0 b} $T$ admits a {\em distributionally irregular vector}, that is, 
			a vector $y \in Y$ for which there are $m \in \mathbb{N}$ and $A, B \subset \mathbb{N}$ with $\udens(A) = \udens(B) = 1$ such that
			\[
			\lim_{n \in A} T^n y = 0
			\quad \text{and} \quad
			\lim_{n \in B} \|T^n y\|_m = \infty.
			\]
			\item\label{T0 c}  $T$ satisfies the DCC.
		\end{enumerate}
	\end{lemma}
In the next result, we prove our most general result concerning distributional chaos: a characterization of weighted backward shifts on admissible Fréchet sequence spaces over directed trees.
    \begin{theorem}\label{teo:carac:dist}
     Let $V$ be a directed tree with generations 
enumerated according to some fixed vertex $v_0$. Let $\left( X, (\| \cdot \|_n)_{n\in \N} \right)$ be an admissible Fréchet sequence space over $V$. Suppose that the weighted backward shift $B_\lambda$, with nonzero weights $\lambda := (\lambda_v)_{v \in V}$, is well-defined and continuous on $X$. Then $B_\lambda$ is distributionally chaotic if and only if there exist $D\subset \N$ with $\udens(D)=1$, $S\subset \{n\in\mathbb Z:\Gen_n\neq\varnothing\}$ and, for each $n\in S$, a vector $f_n \in X$ supported in $\Gen_n$ such that the following conditions hold:
\begin{enumerate}[label=$(\Alph*)$]
    \item \label{teo:carac:dist A} for all $n\in S$, we have $\lim_{k\in D}B_{\lambda}^kf_n=0$.
    \item \label{teo:carac:dist B} there exist $m \in \N$ and an increasing sequence $(N_k)_{k \in \N}$ of positive integers such that, for each $k \in \N$, there is a finite subset $F_k\subset S$ with $\left\| \sum_{j\in F_k}f_j\right\|_{\rho(k)}$ $\ne 0$ and
    $$\card\left\{ 1 \le n \le N_k : \frac{\left\| \sum_{j\in F_k}^{}B_{\lambda}^nf_{j} \right\|_m}{\left\| \sum_{j\in F_k}^{}f_{j}\right\|_{\rho(k)}}>k\right\}> \left( 1- \frac{1}{k} \right) N_k,$$
    where $\rho(k):= \max\{m,k\}$.
\end{enumerate}
    \end{theorem}
    \begin{proof}
        ($\Rightarrow$): Suppose that $B_\lambda$ is distributionally chaotic. Then, by Lemma \ref{T0}, $B_{\lambda}$ admits a distributionally irregular vector $f\in X$, that is, there are $m \in \mathbb{N}$ and $D, E \subset \mathbb{N}$ with $\udens(D) = \udens(E) = 1$ such that
			\begin{equation}\label{eqD1}
			    \lim_{j \in D} B_{\lambda}^j f = 0
			\quad \text{and} \quad
			\lim_{j \in E} \|B_{\lambda}^j f\|_m = \infty.
			\end{equation}
            We take $m\in\N$ in \eqref{eqD1} sufficiently large such that $\left\| f\right\|_m\ne 0$. Fix such an $m\in\N$ and consider the definition of $\rho(k)$ given in the statement for this $m$. Define $S:=\left\{ n\in \Z:\Gen_n \neq \emptyset \right\}.$ By Lemma \ref{LM1}, the family of operators $(\Phi_{n})_{n\in \Z}$ given by 
        $$\Phi_{n}: g\in X \to g\chi_{\operatorname{Gen}_n}\in X, \quad \text{for } n\in \Z,$$
        is equicontinuous. Thus, given a neighborhood $O\subset X$ of $0$, there exists a neighborhood $U\subset X$ of $0$  that satisfies
 \begin{equation}\label{eqD2}
     g\in U \implies g\chi_{\Gen_n}\in O, \quad \text{for all }n \in \Z.
 \end{equation}
 By (\ref{eqD1}) and (\ref{eqD2}), there exists $k_0\in \N$ such that 
    $$k\in D\,\text{ and }\,k\ge k_0 \implies B_{\lambda}^kf\in U \implies (B_{\lambda}^kf)\chi_{\text{Gen}_n}\in O,$$
     for all $n\in \Z$. By Lemma \ref{LM2}, $B_\lambda^{k}(f\chi_{\Gen_n})=\chi_{\Gen_{n-k}}B_\lambda^{k}f$  for all $k\in \N$ and $n\in S$. Therefore, item \ref{teo:carac:dist A} holds with $f_n:= f\chi_{\Gen_{n}}$, $n\in S$. On the other hand, since for each $k\in \N$
     $$\overline{\text{dens}}\left\{ n \in \N : \left\| B_{\lambda}^nf\right\|_m > k(\left\| f \right\|_{\rho(k)} +1)\right\}= \overline{\text{dens}}(E)=1,$$
we can define an increasing sequence $(N_k)_{k\in\N}\subset \N$ such that 
\begin{equation}\label{eqD5}
    \text{card}\left\{1\le n\le N_k : \left\| B_{\lambda}^nf\right\|_m > k(\left\| f \right\|_{\rho(k)} +1)\right\}> N_k\left( 1-\frac{1}{2k} \right)
\end{equation}
for each $k\in\N$. Fix $k\in \N$ and define 
$$J_k:=\left\{1\le n\le N_k : \left\| B_{\lambda}^nf\right\|_m > k(\left\| f \right\|_{\rho(k)} +1)\right\}.$$
For each $N\in \N$, define $I_N:=\left\{ -N,\cdots, N \right\}$. By condition \ref{frechet:condition- B2} from Definition \ref{defAd}, we have 
$$ k(\left\| f \right\|_{\rho(k)} +1)<\left\| B_{\lambda}^nf\right\|_m =\lim_{N\to \infty}\left\| B_{\lambda}^{n}\left( \sum_{j\in S\cap I_N}^{}f_j \right) \right\|_m$$
     for each $n\in J_k$. Thus, there exists $N\in \N$ large enough such that, putting $F_k:=S\cap I_N$, the following inequalities hold
     \begin{equation}\label{eqD6}
         \left\| B_{\lambda}^{n}\left( \sum_{j\in F_k}^{}f_j \right) \right\|_m>k(\left\| f \right\|_{\rho(k)} +1)> k\left\| \sum_{j\in F_k}^{}f_j \right\|_{\rho(k)}>0\quad \text{for all }n\in J_k.
     \end{equation}
      Therefore, by (\ref{eqD5}) and (\ref{eqD6}), we obtain item \ref{teo:carac:dist B}.

     ($\Leftarrow$):  For this implication, we will use the Distributional Chaos Criterion. By item \ref{teo:carac:dist A}, item \ref{dcc a} of the DCC holds. Fix $m\in \N$ provided in item \ref{teo:carac:dist B}. For each $k\in \N$, consider 
 $$g_{k}:= \frac{\sum_{j\in F_k}^{}f_j}{k\left\| \sum_{j\in F_k}^{}f_j\right\|_{\rho(k)}},$$
 where the finite set $F_k \subset S$ comes from item \ref{teo:carac:dist B}. Fix $n\in \N$ and $\varepsilon>0$. Take $k_0 \in \N$ such that $k_0 > \max \{ m,n\}$ and $\frac{1}{k_0} < \varepsilon.$ Then, for all $k \ge k_0$, we have 
 $$\left\| g_k \right\|_n:= \frac{\left\| \sum_{j\in F_k}^{}f_j \right\|_n}{k\left\| \sum_{j\in F_k}^{}f_j\right\|_{\rho(k)}}\le \frac{1}{k}<\varepsilon.$$
 Thus, $\lim_{k\to \infty}g_k=0.$ By the definition of the metric $d$ in (\ref{d}), there exists $\delta >0$ such that
 \begin{equation}\label{eqD7}
     f\in X\text{ and }\left\| f \right\|_m > 1  \implies  d(f, 0)> \delta.
 \end{equation}
 Thus, by the definition of $g_k$, item \ref{teo:carac:dist B}, and (\ref{eqD7}), we have
 $$\text{card} \left\{ 1 \le n \le N_k :d(B_{\lambda}^ng_k,0)> \delta \right\}> (1- k^{-1}) N_k\quad \text{ for all }k\in\N.$$
Therefore, by the Distributional Chaos Criterion, we conclude that $B_\lambda$ is distributionally chaotic.
\end{proof}
Next, we make the same observation, now concerning Theorem~\ref{teo:carac:dist}, as that made in Remark~\ref{rmkLY}.
\begin{remark}\label{rmkD}
        In Theorem \ref{teo:carac:dist}, for each $n\in S$, we chose a single vector $f_n \in X$, supported in $\Gen_n$, such that the family $(f_n)_{n\in S}$ satisfies conditions \ref{teo:carac:dist A} and \ref{teo:carac:dist B}. From the proof, we see that this was not necessary, that is, we could restate the theorem by taking, for each $n\in S$, not only a single vector, but rather a sequence $(f_{n,k})_{k\in \N}\subset X$ of vectors supported in $\Gen_n$, such that the family of vectors $\left\{ f_{n,k}:n\in S\text{ and }k\in \N \right\}$ satisfies conditions \ref{teo:carac:dist A} and \ref{teo:carac:dist B} of the theorem. 
    \end{remark}
    Applying Theorem \ref{teo:carac:dist} to the case $V=\Z$, we recover the following characterization from \cite{joaojmaa}.
    \begin{corollary}\cite[Theorem 12]{joaojmaa}
          Let $X$ be a Fréchet sequence space over $\mathbb{Z}$, endowed with an increasing sequence $(\| \cdot \|_n)_{n \in \N}$ of seminorms, in which the sequence 
		$(e_n)_{n \in \mathbb{Z}}$ of canonical vectors is a Schauder basis. Suppose that the bilateral weighted backward shift $B_w$, with nonzero weights $w := (w_n)_{n \in \mathbb{Z}}$, is well-defined and continuous on $X$. Then $B_w$ is distributionally chaotic if and only if there exist $D \subset \N$ with $\udens(D)=1$ and $I \subset \Z$ such that the following conditions hold:
		\begin{enumerate}[label=$(\Alph*)$]
			\item \label{TThm12 JJMA A} for all $i \in I$, we have
			$$\lim_{n \in D} w_{i-n}\cdots w_{i-1}e_{i-n}=0.$$
			\item\label{TThm12 JJMA B} there exist $m \in \N$ and an increasing sequence $(N_k)_{k \in \N}$ of positive integers such that for each $k \in \N$ there are $r:=r(k) \in \N$, indices $i_{1,k}, \cdots, i_{r,k} \in I$ and scalars $b_{1,k}, \cdots, b_{r,k} \in \K \setminus\left\{ 0 \right\}$ with $\left\| \sum_{j=1}^{r} b_{j,k} e_{i_{j,k}} \right\|_{p(k)} \ne 0$ and
			$$\card\left\{ 1 \le n \le N_k :\frac{\left\| \sum_{j=1}^{r} b_{j,k} w_{i_{j,k}-n}\cdots w_{i_{j,k} -1}e_{i_{j,k}-n} \right\|_{m} }{\left\| \sum_{j=1}^{r} b_{j,k} e_{i_{j,k}} \right\|_{p(k)}}> k \right\}> (1- k^{-1}) N_k,$$
			where $p(k):= \max\{m,k\}$.
            \end{enumerate}
      \end{corollary}
    As a consequence of Theorem \ref{teo:carac:dist}, we obtain the following corollary, which characterizes distributional chaos in the context of weighted shifts on Köthe sequence spaces over directed trees. This corollary is obtained by directly applying Theorem \ref{teo:carac:dist} to the seminorms from Definition \ref{def1}. We observe that the next corollary can also be stated in the form mentioned in Remark~\ref{rmkD}.
  
 \begin{corollary} \label{cor:dist:Kothe}
       Let $V$ be a directed tree with generations enumerated according to some fixed vertex $v_0$. Consider a Köthe sequence space $X := \lambda_p(A,V)$, where $A := (a_{v,k})_{v \in V, k \in \N}$ is a Köthe matrix and $p \in \{0\} \cup [1,\infty)$. Suppose that the weighted backward shift $B_\lambda$, with nonzero weights $\lambda := (\lambda_v)_{v \in V}$, is well-defined and continuous on $X$. Then $B_\lambda$ is distributionally chaotic if and only if there exist $D\subset \N$ with $\udens(D)=1$, a subset $S\subset \{n\in\mathbb Z:\Gen_n\neq\varnothing\}$ and a family of scalars $\left\{ c_{v}\in \K : v \in \Gen_n \, \text{ for some } n\in S \right\}$ with $\sum_{v\in \operatorname{Gen}_n}^{}c_{v}e_{v}\in X$ for all $n\in S$, such that the following conditions hold:
        \begin{enumerate}[label=$(\Alph*)$]
            \item \label{cor:dist:Kothe A} for all $n \in S$, we have
            $$\begin{cases}
\displaystyle\lim_{k\in D}\underset{v \in \operatorname{Gen}_{n-k}}{\sup} \Bigg|\sum_{u\in \operatorname{Chi}^{k}(v)}^{} \lambda(v\to u)c_{u} a_{v, \ell} \Bigg|=0\quad \text{ for all }\ell\in \N, & \text{if }p=0,\\[4mm]
\displaystyle\lim_{k\in D}\underset{v \in \operatorname{Gen}_{n-k}}{\sum} \Bigg|\sum_{u\in \operatorname{Chi}^{k}(v)}^{} \lambda(v\to u)c_{u} a_{v, \ell} \Bigg|^p=0\quad \text{ for all }\ell\in \N, & \text{if }p\in[1,\infty).
\end{cases};$$
\item \label{cor:dist:Kothe B} there exist $m \in \N$ and an increasing sequence $(N_k)_{k \in \N}$ of positive integers such that, for each $k \in \N$, there is a finite subset $F_k\subset S$ with 
    $$\underset{j\in F_k}{\max} \underset{u \in\operatorname{Gen}_{j}}{\sup} \left| c_{u}a_{u,\rho(k)} \right|\ne 0\,\,\text{ if }\,p=0,\quad \text{or}\quad \underset{j\in F_k}{\sum}\underset{u \in \operatorname{Gen}_{j}}{\sum} \left| c_{u}a_{u,\rho(k)} \right|^p \ne 0\,\,\text{ if }\,p\in[1,\infty),$$
and
$$\begin{cases}
    \card\left\{1 \le n \le N_k : \frac{\underset{j\in F_k}{\max}\, \underset{v \in\operatorname{Gen}_{j-n}}{\sup}\Big| \underset{u\in \operatorname{Chi}^{n}(v)}{\sum}\lambda(v\to u) c_{u}a_{v,m}\Big|}{ \underset{j\in F_k}{\max} \underset{u \in\operatorname{Gen}_{j}}{\sup} | c_{u}a_{u,\rho(k)} |}>k \right\}> \left( 1- \frac{1}{k} \right) N_k, & \text{if }p=0,\\
    \card\left\{ 1 \le n \le N_k :\frac{ \underset{j\in F_k}{\sum}\,\underset{v \in \operatorname{Gen}_{j-n}}{\sum}\Big| \underset{u\in \operatorname{Chi}^{n}(v)}{\sum}\lambda(v\to u) c_{u}a_{v,m}\Big|^p}{\underset{j\in F_k}{\sum}\underset{u \in \operatorname{Gen}_{j}}{\sum} | c_{u}a_{u,\rho(k)} |^p}>k^p \right\}> \left( 1- \frac{1}{k} \right) N_k, & \text{if }p\in[1,\infty),
\end{cases}$$
         where $\rho(k):= \max\{k,m\}$.       
    \end{enumerate}
    \end{corollary}
    \begin{example}
       Fix $N\in\N$ such that $N\ge2$. Consider the unrooted directed $N$-adic tree $V$ centered at the vertex $v_0$. Let $A=(a_{v,k})_{v\in V,\,k\in\N}$ be a Köthe matrix that satisfies (\ref{eqndia1}). Let $p\in\{0\}\cup[1,\infty)$ and $X:=\lambda_p(A,V)$. Fix $\alpha\in\K\setminus\{0\}$ and consider the weighted backward shift $B_\lambda$ with weights $\lambda_v=\alpha$ for every $v\in V$, satisfying \eqref{eqndia2}. Consider the constant $C_p$ given by \eqref{eqndia3}. As observed in Example \ref{ex1}, $B_\lambda$ is well-defined and continuous on $X$. We will prove that the following assertions are equivalent:
\begin{enumerate}[label=$(\roman*)$]
    \item \label{dist:exe: i} $B_{\lambda}$ is distributionally chaotic;
    \item \label{dist:exe: ii} there exist $m\in \N$ and a subsequence $(N_k)_{k\in\N}\subset \N$ such that, for each $k\in\N$, there are $r_k\in \N$, distinct indices $i_{k,1}, \cdots ,i_{k,r_k}\in \Z$, and positive scalars $c_{k,1}\cdots, c_{k,r_k}$ with 
    $$\begin{cases}
\displaystyle \frac{1}{N_k}\card\left\{ 1\le n\le N_k:C_0^n\frac{\max\limits_{1\le j\le r_k}c_{k,j}b_{i_{k,j}-n,m}}{\max\limits_{1\le j\le r_k}c_{k,j}b_{i_{k,j},\rho(k)}}>k\right\}>1-\frac{1}{k}, & \text{if }p=0,\\
\displaystyle \frac{1}{N_k}\card\left\{1\le n\le N_k:C_p^{np}\frac{\sum\limits_{j=1}^{r_k}c_{k,j}b_{i_{k,j}-n,m}^p}{\sum\limits_{j=1}^{r_k}c_{k,j}b_{i_{k,j},\rho(k)}^p}>k^p\right\}>1-\frac{1}{k}, & \text{if }p\in[1,\infty).
\end{cases}$$
\end{enumerate}
\ref{dist:exe: i}$\Rightarrow$\ref{dist:exe: ii}: By condition \ref{cor:dist:Kothe B} of Corollary \ref{cor:dist:Kothe}, there exist $S\subset \Z$, $m\in \N$, $(f_n)_{n\in S}\subset X$, where $f_n$ is supported on $\text{Gen}_n$ for every $n\in S$, and an increasing sequence $(N_k)_{k \in \N}$ of positive integers such that for each $k \in \N$ there is a finite subset $F_k\subset S$ with $\sum_{j\in F_k}^{}f_j\ne 0$ and 
\begin{equation}\label{eqx21}
    \begin{cases}
\displaystyle\text{card}\left\{1\le n\le N_k : \frac{\max_{j\in F_k}\left\| B_{\lambda}^{n}f_j \right\|_m}{\max_{j\in F_k}\left\|f_j  \right\|_{\rho(k)}}>k \right\}>\left( 1- \frac{1}{k} \right) N_k,& \text{if }p=0,\\
\displaystyle\text{card}\left\{ 1\le n\le N_k :\frac{\sum_{j\in F_k}\left\| B_{\lambda}^{n}f_j \right\|_{m}^p}{\sum_{j\in F_k}\left\|f_j  \right\|_{\rho(k)}^p}>k^p \right\}>\left( 1- \frac{1}{k} \right) N_k,& \text{if }p\in[1,\infty).
\end{cases}
\end{equation}
Fix $k\in \N$ and write $F_k=\left\{ i_{k,1},\cdots, i_{k,r_k} \right\}$ for some $r_k\in\N$. Define, for each $j\in \{1, \cdots, r_k\}$,
$$c_{k,j}:=
\begin{cases}
\displaystyle\sup_{u\in\Gen_{i_{k,j}}}|f_{i_{k,j}}(u)|, & p=0,\\
\displaystyle\sum_{u\in\Gen_{i_{k,j}}}|f_{i_{k,j}}(u)|^p, & p\in[1,\infty).
\end{cases}$$
Note that, for all $j\in\{1,\cdots, r_k\}$ and $\ell\in \N$, we have 
\begin{equation}\label{eqx22}
    \|f_{i_{k,j}}\|_\ell=
\begin{cases}
\displaystyle c_{k,j}b_{i_{k,j},\ell}, & p=0,\\
\displaystyle c_{k,j}^{1/p}b_{i_{k,j},\ell}, & p\in[1,\infty).
\end{cases}
\end{equation}
Moreover, the triangle inequality, if $p=0$, and Hölder's inequality, if $p\in[1,\infty)$, imply that
\begin{equation}\label{eqx23}
    \|B_\lambda^nf_{i_{k,j}}\|_m
\le
\begin{cases}
\displaystyle C_0^nc_{k,j}b_{i_{k,j}-n,m}, & p=0,\\
\displaystyle C_p^nc_{k,j}^{1/p}b_{i_{k,j}-n,m}, & p\in[1,\infty),
\end{cases}
\end{equation}
for all $n\in \N$. Therefore, by \eqref{eqx21}-\eqref{eqx23}, we obtain item \ref{dist:exe: ii}.

\ref{dist:exe: ii}$\Rightarrow$\ref{dist:exe: i}: For this implication we use the version of Corollary \ref{cor:dist:Kothe} commented in Remark \ref{rmkD}. Suppose that item \ref{dist:exe: ii} holds. For each $k\in\N$ and $1\le j\le r_k$, choose $z_{k,j}\in\Gen_{i_{k,j}-N_k-1}$ and two distinct vertices $x_{k,j},y_{k,j}\in\Chi(z_{k,j})$. Define
$$
h_{k,j}:=
\begin{cases}
\displaystyle c_{k,j}\left(\sum_{u\in\Chi^{N_k}(x_{k,j})}e_u-\sum_{u\in\Chi^{N_k}(y_{k,j})}e_u\right), & p=0,\\
\displaystyle \left(\frac{c_{k,j}}{2N^{N_k}}\right)^{1/p}\left(\sum_{u\in\Chi^{N_k}(x_{k,j})}e_u-\sum_{u\in\Chi^{N_k}(y_{k,j})}e_u\right), & p\in[1,\infty).
\end{cases}
$$
Then $h_{k,j}\in X$, $\supp(h_{k,j})\subset\Gen_{i_{k,j}}$ and $B_{\lambda}^{N_k+1}(h_{k,j})=0$. Hence, item \ref{cor:dist:Kothe A} of Corollary \ref{cor:dist:Kothe} holds. Moreover, for $1\le n\le N_k$,
$$
B_\lambda^nh_{k,j}=
\begin{cases}
\displaystyle c_{k,j}\alpha^nN^n\left(\sum_{u\in\Chi^{N_k-n}(x_{k,j})}e_u-\sum_{u\in\Chi^{N_k-n}(y_{k,j})}e_u\right), & p=0,\\
\displaystyle \left(\frac{c_{k,j}}{2N^{N_k}}\right)^{1/p}\alpha^nN^n\left(\sum_{u\in\Chi^{N_k-n}(x_{k,j})}e_u-\sum_{u\in\Chi^{N_k-n}(y_{k,j})}e_u\right), & p\in[1,\infty).
\end{cases}
$$
Thus, for all $\ell\in \N$ and $1\le n\le N_k$, we have
\begin{equation}\label{eqx24}
    \|h_{k,j}\|_\ell=
\begin{cases}
\displaystyle c_{k,j}b_{i_{k,j},\ell}, & p=0,\\
\displaystyle c_{k,j}^{1/p}b_{i_{k,j},\ell}, & p\in[1,\infty),
\end{cases}
\quad
\|B_\lambda^nh_{k,j}\|_m=
\begin{cases}
\displaystyle C_0^nc_{k,j}b_{i_{k,j}-n,m}, & p=0,\\
\displaystyle C_p^nc_{k,j}^{1/p}b_{i_{k,j}-n,m}, & p\in[1,\infty).
\end{cases}
\end{equation}
Then, by \eqref{eqx24}, for each $k \in \N$ and $n\in \{1,\cdots, N_k\}$, we have 
\begin{equation}
    \begin{cases}\label{eqx25}
\displaystyle \frac{\max\limits_{1\le j\le r_k}\left\|B_\lambda^nh_{k,j}\right\|_m}{\max\limits_{1\le j\le r_k}\left\|h_{k,j}\right\|_{\rho(k)}}=C_0^n\frac{\max\limits_{1\le j\le r_k}c_{k,j}b_{i_{k,j}-n,m}}{\max\limits_{1\le j\le r_k}c_{k,j}b_{i_{k,j},\rho(k)}}, & p=0,\\
\displaystyle \frac{\sum\limits_{j=1}^{r_k}\left\|B_\lambda^nh_{k,j}\right\|_m^p}{\sum\limits_{j=1}^{r_k}\left\|h_{k,j}\right\|_{\rho(k)}^p}=C_p^{np}\frac{\sum\limits_{j=1}^{r_k}c_{k,j}b_{i_{k,j}-n,m}^p}{\sum\limits_{j=1}^{r_k}c_{k,j}b_{i_{k,j},\rho(k)}^p}, & p\in[1,\infty).
\end{cases}
\end{equation}
Therefore, by \eqref{eqx25} and item \ref{dist:exe: ii}, we conclude that item \ref{cor:dist:Kothe B} of Corollary \ref{cor:dist:Kothe} holds.
    \end{example}

For our next proposition, we will need the following results from \cite{Nilson0}.

\begin{lemma}[{\cite[Proposition 7]{Nilson0}}]\label{Nilson0:prop7}
    Let $T:X\to X$ be a continuous linear operator on a Fréchet space $X$, and $m\in \N$. The following assertions are equivalent:
    \begin{enumerate}[label=$(\roman*)$]
        \item \label{PProp12:N0:i} there exist $\varepsilon>0$, a sequence $(y_k)_{k\in\N}$ in $X$ and an increasing sequence $(N_k)_{k\in\N}$ in $\N$ such that $\lim_{k\to\infty}y_k=0$ and $\operatorname{card}\{1\leq j\leq N_k:\norm{T^jy_k}_m>\varepsilon\}\geq N_k(1-k^{-1})$ for all $k\in\N$.
        
        \item\label{PProp12:N0:ii} there exist $y\in X$ and $A\subset\N$ with $\udens(A)=1$ such that $\lim_{n\in A}\norm{T^ny}_m=\infty$.
        
        \item\label{PProp12:N0:iii} the set $\{y\in X:\exists A\subset\N\text{ with }\udens(A)=1\text{ such that }\lim_{n\in A}\norm{T^ny}_m=\infty\}$ is residual in $X$.
\end{enumerate}
\end{lemma}

\begin{lemma}[{\cite[Proposition 9]{Nilson0}}]\label{Nilson0:prop9}
    Let $T:X\to X$ be a continuous linear operator on a Fréchet space $X$. Consider the set
    \begin{equation}\label{seto}
        \mathcal{O}:=\{x\in X:\exists A\subset\N\text{ with }\udens(A)=1\text{ such that }\lim_{n\in A}T^nx=0\}.
    \end{equation}
    Suppose that there exists a dense subset $X_0\subset \mathcal{O}$. Then the set $\mathcal{O}$ is residual in $X$.
\end{lemma}

Before proving the proposition, let us show the following lemma, which is a generalization of \cite[Theorem 2.7]{HeYin}.

\begin{lemma}\label{lemmaP}
   Let $T:X\to X$ be a continuous linear operator on a non-trivial separable Fréchet space $X$. Then the following assertions are equivalent.
   \begin{enumerate}[label=$(\alph*)$]
       \item \label{lemmaP a} $T$ is densely distributionally chaotic.
       \item \label{lemmaP b} There exist $z\in X$ and a dense subset $X_0$ of $X$ such that:
       \begin{enumerate}[label=$(\roman*)$]
           \item \label{lemmaP i} there exist $m\in\N$ and a set $B\subset\N$ with $\udens(B)=1$ such that $\lim_{ n\in B}\norm{T^n z}_m=\infty$;

           \item \label{lemmaP ii} for every $x\in X_0$, there is a set $A_x\subset\N$ with $\udens(A_x)=1$ and $\lim_{\substack{n\in A_x}}T^n x=0$.  
       \end{enumerate}
   \end{enumerate}
\end{lemma}
\begin{proof}
    \ref{lemmaP a}$\Rightarrow$\ref{lemmaP b}: By definition, $T$ admits a dense distributionally scrambled set $\Gamma\subset X$. Fix $x_0\in \Gamma$ and consider the set $D:=\{x-x_0: x\in\Gamma\setminus\{x_0\}\}$. Note that $D$ is dense. Fix $z\in D$. Then 
    \begin{equation}\label{eqlemden}
        \udens\{n\in\N: d(T^{n}z,0)<\tau\}=1\quad \text{for all }\tau>0, 
    \end{equation}
    where $d$ is the usual metric in $X$. By \eqref{eqlemden}, using a diagonal argument, we obtain a set $A\subset\N$ such that $\udens(A)=1$ and $\lim_{n\in A}T^nz=0$. Thus, taking $X_0:=D$, we conclude item \ref{lemmaP ii}. Item \ref{lemmaP i} follows from the existence of a distributionally irregular vector, guaranteed by Lemma \ref{T0}.

    \ref{lemmaP b}$\Rightarrow$\ref{lemmaP a}: Let $Y$ be the set of all $x\in X$ for which there exists $B\subset\N$ with $\udens(B)=1$ such that $\lim_{ n\in B}\norm{T^n x}_m=\infty$. By Lemmas \ref{Nilson0:prop7} and \ref{Nilson0:prop9}, the sets $X_0$ and $Y$ are residual and, therefore, $R:=\mathcal{O}\cap Y$ is residual, where the set $\mathcal{O}$ is define as \eqref{seto}. Since the map $\Phi:X\times X\to X$, $\Phi(x,y)=x-y$, is continuous, open and surjective, the set
    $$\Phi^{-1}(R)=\{(x,y)\in X\times X:x-y\in R\}$$ is residual in $X\times X$. By Mycielski's theorem, there exists a dense Mycielski set $\Gamma\subset X$ such that $x-y\in R$ for every distinct $x,y\in\Gamma$. Therefore, $\Gamma$ is a dense distributionally scrambled set.
\end{proof}

Let $V$ be a directed tree, $X$ be an admissible Fréchet sequence space over $V$ and $\lambda=(\lambda_v)_{v\in V}$ be a weight sequence such that the weighted backward shift $B_\lambda$ is well defined and continuous on $X$. Consider the following condition:
\begin{enumerate}[label=$(C3)$]
    \item \label{cond:(C):dist} There are $v_0\in V$ and $D\subset \N$ with $\udens(D)=1$ such that
$$\lim_{n\in D} B_\lambda^{n}e_{v_0} =\lim_{n\in D}\lambda\bigl(\Par^{n}(v_0)\to v_0\bigr) e_{\Par^{n}(v_0)}=0.$$ 
\end{enumerate}
In the next proposition, we establish a characterization of distributional chaos for weighted backward
shifts on admissible Fréchet sequence spaces over directed trees, under the assumption of condition \ref{cond:(C):dist}.
\begin{proposition}\label{propD3}
Let $V$ be a directed tree. Let $\left( X, (\| \cdot \|_n)_{n\in \N} \right)$ be an admissible Fréchet sequence space over $V$. Suppose that the weighted backward shift $B_\lambda$, with nonzero weights $\lambda := (\lambda_v)_{v \in V}$, is well-defined and continuous on $X$. Suppose that condition \ref{cond:(C):dist} holds. Then the following assertions are equivalent:    
\begin{enumerate}[label=$(\roman*)$]
    \item \label{propD3 i} $B_{\lambda}$ is distributionally chaotic;
    \item \label{propD3 ii} $B_{\lambda}$ is densely distributionally chaotic;
    \item \label{propD3 iii} there exist $m\in \N$ and $f\in X$ such that $\lim_{n\in D}\left\| B_{\lambda}^nf \right\|_m=\infty$ for some $D\subset \N$ with $\udens(D)=1$.
    \item \label{propD3 iv} there exist $m \in \N$ and an increasing sequence $(N_k)_{k \in \N}$ of positive integers, such that for each $k \in \N$ there are a finite set $F_k\subset V$, and scalars $c_v$ for $v\in F_k$, with $\left\| \sum_{v\in F_k}^{} c_{v}e_v\right\|_{\rho(k)}$ $\ne 0$ and
    \begin{equation}\label{propD2}
        \card\left\{1\le n\le N_k: \frac{\left\| \sum_{v\in F_k}^{}\lambda(\Par^{n}(v)\to v) c_{v}e_{\Par^{n}(v)} \right\|_m}{\left\| \sum_{v\in F_k}^{} c_{v}e_v\right\|_{\rho(k)}}>k \right\}>(1-k^{-1})N_k,
    \end{equation}
    where $\rho(k):= \max\{m,k\}$.
\end{enumerate}
\end{proposition}
\begin{proof}
Let $\mathcal{O}$ be the set defined in \eqref{seto} relative to $B_\lambda$. Arguing as in the proof of Proposition \ref{prop:Fernanda17}, we conclude that $c_{00}(V)\subset\mathcal{O}$ and therefore $\mathcal{O}$ is dense.

    \ref{propD3 i}$\Rightarrow$\ref{propD3 ii}: Since $B_{\lambda}$ is distributionally chaotic, by Lemma \ref{T0}, we conclude that $B_{\lambda}$ admits a distributionally irregular vector. Since $\mathcal{O}$ is dense, by Lemma \ref{lemmaP}, $B_{\lambda}$ is densely distributionally chaotic.

    \ref{propD3 ii}$\Rightarrow$\ref{propD3 iii}: This implication is immediate from Lemma \ref{T0}. 

    \ref{propD3 iii}$\Rightarrow$\ref{propD3 iv}: By Lemma \ref{lemmaP}, $B_\lambda$ is densely distributionally chaotic and, in particular, distributionally chaotic. Thus, condition \ref{teo:carac:dist B} of Theorem \ref{teo:carac:dist} holds. Therefore, by condition \ref{frechet:condition- B1} of Definition \ref{defAd} and the continuity of $B_{\lambda}$ and of the seminorms, condition \ref{teo:carac:dist B} of Theorem \ref{teo:carac:dist} implies item \ref{propD3 iv}.

    \ref{propD3 iv}$\Rightarrow$\ref{propD3 i}: For each $k\in \N$ define 
    $$g_k:= \frac{\sum_{v\in F_k}^{} c_{v}e_v}{k\left\| \sum_{v\in F_k}^{} c_{v}e_v \right\|_{\rho(k)}}.$$
    Arguing as in the proof of implication $(\Leftarrow)$ of Theorem \ref{teo:carac:dist}, we have that \eqref{propD2} implies item \ref{PProp12:N0:i} of Lemma \ref{Nilson0:prop7}. Then there exist $f\in X$ and $A\subset\N$ with $\udens(A)=1$ such that $\lim_{n\in A}\norm{B_{\lambda}^nf}_m=\infty$. Since $\mathcal{O}$ is dense, by Lemma \ref{lemmaP}, we conclude that $B_\lambda$ is densely distributionally chaotic and, in particular, distributionally chaotic.
\end{proof}
 If $V$ is rooted, then condition \ref{cond:(C):dist} is automatically satisfied. Hence, we obtain the following corollary.
\begin{corollary}
    Let $V$ be a rooted directed tree. Let $X$ be an admissible Fréchet sequence space over $V$. Suppose that the weighted backward shift $B_\lambda$, with nonzero weights $\lambda := (\lambda_v)_{v \in V}$, is well-defined and continuous on $X$. Then assertions \ref{propD3 i}-\ref{propD3 iv} from Proposition \ref{propD3} are equivalent.
\end{corollary}
      In light of condition \ref{cond:(C):dist}, we obtain the following characterization for the spaces $\ell^p(V)$, $1\le p<\infty$, and $c_0(V)$.
      \begin{corollary}\label{corcasonor}
          Let $V$ be a directed tree and let $X:=\ell^p(V)$ $(1\le p<\infty)$ or $X:=c_0(V)$. Suppose that the weighted backward shift $B_\lambda$, with nonzero weights $\lambda := (\lambda_v)_{v \in V}$, is well-defined and continuous on $X$. Moreover, suppose that condition \ref{cond:(C):dist} holds. Then the following assertions are equivalent.
    \begin{enumerate}[label=$(\roman*)$] 
        \item \label{corcasonor:item:1} $B_\lambda$ is distributionally chaotic on $X$;
        \item \label{corcasonor:item:2} $B_\lambda$ is densely distributionally chaotic on $X$;
        \item \label{corcasonor:item:3}The following condition holds:
        \begin{enumerate}[label=$(\alph*)$] 
            \item \label{corcasonor:item:3.1} If $X=c_0(V)$, then for all $k\in \N$, we have 
        $$\limsup_{N\to\infty}\frac{1}{N}\operatorname{card}\left\{1\le n\le N:\sum_{\Chi^n(u)}\abs{\lambda(u\to v)}>k\,\text{ for some }u\in V\right\}=1.$$
            \item \label{corcasonor:item:3.2} If $X=\ell^p(V)$ for $1\le p<\infty$, then there exists an increasing sequence $(N_k)_{k \in \N}$ of positive integers, where for each $k \in \N$ there are a finite set $F_k\subset V$ and nonzero scalars $c_v$ for $v\in F_k$ such that
    $$\operatorname{card}\left\{1\leq n\leq N_k:\frac{\sum_{u\in\Par^n(F_k)}\left|\sum_{\substack{v\in F_k\\ \Par^n(v)=u}}\lambda(u\to v)c_v\right|^p}{\sum_{v\in F_k}|c_{v}|^p}>k^p\right\}>(1-k^{-1})N_k.$$
        \end{enumerate}
       
        \end{enumerate}
      \end{corollary}
      \begin{proof}
          By Proposition \ref{propD3}, it is sufficient to show that item \ref{corcasonor:item:3} is equivalent to item \ref{propD3 iv} of Proposition \ref{propD3}. Since, in the case $X=\ell^p(V)$ for $1\le p<\infty$, \ref{corcasonor:item:3.2} is precisely condition \ref{propD3 iv} of Proposition \ref{propD3}, it remains only to consider the case $X=c_0(V)$. We denote the norm in $c_0(V)$ by $\left\| \cdot  \right\|_{\infty}$.  Initially, note that we can take a family $(\alpha_v)_{v\in V}$ of scalars with $\abs{\alpha_v}=1$ such that $\lambda_v\alpha_v=\abs{\lambda_v}\alpha_{\Par(v)}$, whenever $\Par(v)$ is defined. Therefore, for every $v\in V$ and $n\in\mathbb{N}$ such that $\Par^n(v)$, 
          \begin{equation}\label{eqcor001}
              \lambda(\Par^n(v)\to v)\alpha_v=\abs{\lambda(\Par^n(v)\to v)}\alpha_{\Par^n(v)}.
          \end{equation}
          Now, suppose that item \ref{corcasonor:item:3.1} holds.  Thus, there exists a subsequence $(N_k)_{k\in\N}\subset \N$ such that, for each $k\in \N$,
          $$\frac{1}{N_k}\card(J_k)>1-k^{-1},\quad \text{where } J_k:=\left\{ 1\le n \le N_k : \sum_{\Chi^n(u)}\abs{\lambda(u\to v)}>k\,\text{ for some }u\in V \right\}.$$
          For each $n\in J_k$ there exists $u_{n,k}\in V$ and a finite set $F_{n,k}\subset\Chi^n(u_{n,k})$ such that
          $$\sum_{v\in F_{n,k}}\abs{\lambda(u_{n,k}\to v)}>k.$$
          Fix $k\in \N$ and define $F_k:=\bigcup_{n\in J_k}F_{n,k}$. Thus, by \eqref{eqcor001}, for each $n\in J_k$, we have 
          $$\frac{\left\| \sum_{v\in F_k}^{}\lambda(\Par^{n}(v)\to v) \alpha_{v}e_{\Par^{n}(v)} \right\|_\infty}{\left\| \sum_{v\in F_k}^{} \alpha_{v}e_v\right\|_{\infty}}\ge\sum_{v\in F_{n,k}}\left| \lambda(u_{n,k}\to v) \right|>k.$$
          Therefore, we obtain item \ref{propD3 iv} of Proposition \ref{propD3}. Conversely, suppose that item \ref{propD3 iv} of Proposition \ref{propD3} holds. Then there exists a subsequence $(M_k)_{k\in \N}\subset\N$ such that for each $k\in\N$ there are a finite set $G_k\subset V$ and scalars $(c_v)_{v\in G_k}\subset\K\setminus \{0\}$ with 
          $$\frac{1}{M_k}\card\overset{=:L_k}{\overbrace{\left\{1\le n\le M_k:\frac{\left\|\sum_{v\in G_k}\lambda(\Par^n(v)\to v)c_ve_{\Par^n(v)}\right\|_\infty}{\left\|\sum_{v\in G_k}c_ve_v\right\|_\infty}>k\right\}}}>1-k^{-1}.$$
          Fix $k\in\N$. For each $n\in L_{k+j}$ ($j\in\N_0$), there exists $u\in V$ such that
$$k<\frac{\left|\sum_{v\in G_{k+j}\cap\Chi^n(u)}\lambda(u\to v)c_v\right|}{\max_{v\in G_{k+j}}\abs{c_v}}\le\sum_{v\in\Chi^n(u)}\abs{\lambda(u\to v)}.$$
Hence, for all $\ell\ge k$ we have
$$\frac{1}{M_\ell}\card\left\{1\le n\le M_\ell:\sum_{v\in\Chi^n(u)}\abs{\lambda(u\to v)}>k\text{ for some }u\in V\right\}>1-\ell^{-1},$$
which proves item \ref{corcasonor:item:3.1}.
      \end{proof}
      Examples of shifts that are distributionally chaotic but not hypercyclic chaotic are already known; see, for instance, \cite[Example 20]{joaojmaa}. In what follows,  we will give an example of a weighted backward shift on spaces of the form $\ell^p(V)$, $1\leq p<\infty$, or $c_0(V)$, where $V$ is a directed tree, that is distributionally, but is not hypercyclic. The reader may observe how much simpler our example is than the one given in \cite{joaojmaa}. For an example of a hypercyclic shift that is not distributionally chaotic, see \cite[Example 18]{joaojmaa}. 
      \begin{example}
Consider the weighted backward shift $B_{\lambda}$ on $X=\ell^p(V)$, $1\leq p<\infty$, or $X=c_0(V)$, where the weights $(\lambda_v)_{v\in V}$ and the tree $V$ are given in Example \ref{exNHcY}. In this example, we proved that $B_{\lambda}$ is not hypercyclic and is Li--Yorke chaotic. Here we prove something stronger: $B_{\lambda}$ is densely distributionally chaotic. Indeed, first note that condition \ref{cond:(C):dist} holds, since $V$ is rooted. Therefore, it remains to prove that condition \ref{corcasonor:item:3} of Corollary \ref{corcasonor} is satisfied. We first consider the case $X=c_0(V)$. Fix $k\in\N$. Since $\Chi^n(v_{1,1})=\{v_{1,n+1}\}$ and $\abs{\lambda(v_{1,1}\to v_{1,n+1})}=2^n> k$ for every $n>\lfloor\log_2k\rfloor$, we have
$$\frac{1}{N}\card\left\{1\le n\le N:\sum_{v\in\Chi^n(u)}\abs{\lambda(u\to v)}>k\text{ for some }u\in V\right\}\ge 1-\frac{\lfloor\log_2k\rfloor}{N}$$
for all $N>\lfloor\log_2k\rfloor$. Taking $N\to\infty$, condition \ref{corcasonor:item:3.1} follows. Now, for the case $X=\ell^p(V)$ where $1\le p<\infty$, simple calculations show that condition \ref{corcasonor:item:3.2} holds by taking $N_k=k^2$, $F_k=\{v_{1,N_k}\}$, and $c_{v_{1,N_k}}=1$ for each $k\in\N$.
\end{example}

\section*{Funding} 

The first author was supported by the Coordenação de Aperfeiçoamento de Pessoal de Nível Superior (CAPES, Brazil). The second author was partially supported by the Conselho Nacional de Desenvolvimento Científico e Tecnológico (CNPq, Brazil) through grants 406457/2023-9, 403964/2024-5 and 304165/2025-5. The third author was supported by the Fundação de Amparo à Pesquisa do Estado do Rio de Janeiro (FAPERJ) through a \emph{FAPERJ nota 10} grant.

\section*{Acknowledgments}

This work was initiated during the XVIII ENAMA, held in Bahia, Brazil, in June 2025, where the authors met and discussed the initial results of this research. The first and second authors gratefully acknowledge the support provided by the Universidade Federal da Bahia (UFBA), which made their participation in the conference possible.

	\smallskip
	
	{\footnotesize
\bigskip\noindent
{\sc Mateus Azevedo}

\smallskip\noindent
Universidade Federal de Campina Grande - Campus I, Unidade Acadêmica de Matemática, Rua Aprígio Veloso, 882, Universitário, CEP 58429-900, Campina Grande, Brazil.\\
\textit{ e-mail address}: mateus.vinicius@academico.ufpb.br

\bigskip\noindent
{\sc Fernando Costa Jr.}

\smallskip\noindent
Universidade Federal da Paraíba - Campus I, Departamento de Matemática, Jardim Universitário, s/n, Bairro Castelo Branco, CEP 58051-900, João Pessoa, Brazil.\\
\textit{ e-mail address}: fernando@mat.ufpb.br

\bigskip\noindent
{\sc Jo\~ao V. A. Pinto}

\smallskip\noindent
Departamento de Matem\'atica Aplicada, Instituto de Matem\'atica, Universidade Federal do Rio de Janeiro, 
Caixa Postal 68530, RJ 21941-909, Brazil.\\
\textit{ e-mail address}: pinto.j.v.a@gmail.com

}

\end{document}